%% file: paper.tex
\documentclass[hidelinks,onefignum,onetabnum]{siamart220329}

\usepackage{lipsum}
\usepackage{amsmath,amssymb,amsfonts}
\usepackage{graphicx}
\usepackage{epstopdf}
\usepackage{algorithmic}
\ifpdf
\DeclareGraphicsExtensions{.eps,.pdf,.png,.jpg}
\else
\DeclareGraphicsExtensions{.eps}
\fi

\usepackage{cite}
\usepackage{xcolor}

\input{macros}

\usepackage{pgfplots} 
\pgfplotsset{width=7cm,compat=1.14}
\usepgfplotslibrary{statistics}
\usepackage{ifthen}
\usepackage{tikz} 
\usetikzlibrary{calc,patterns,shapes}
\usetikzlibrary{arrows.meta}
\tikzset{>={Latex[width=2mm,length=2mm]}}
\usetikzlibrary{decorations.pathreplacing}
\usepgfplotslibrary{fillbetween}
\usetikzlibrary{patterns}
\usepackage{booktabs}
\usepackage{siunitx}
\usepackage{verbatim}
\usepackage{pgfplotstable}
\usepackage[caption=false,font=footnotesize,labelfont=sf,textfont=sf]{subfig} 

\newsiamremark{remark}{Remark}
\newsiamremark{hypothesis}{Hypothesis}
\crefname{hypothesis}{Hypothesis}{Hypotheses}
\newsiamthm{claim}{Claim}

\def\mytitle{Algebraic Temporal Blocking for Sparse Iterative Solvers on Multi-Core CPUs}
\def\myshorttitle{Algebraic Temporal Blocking}
\def\myshortauthor{C.~Alappat et al.}

\headers{\myshorttitle}{\myshortauthor}

\title{\mytitle\thanks{Submitted to the editors DATE.
		\funding{This work was funded by the Fog Research Institute under contract no.~FRI-454.}}}

\author{Christie Alappat\thanks{Erlangen National High Performance Computing Center, Friedrich-Alexander-Universität Erlangen-Nürnberg. Email: \{\email{christie.alappat}, \email{georg.hager}, \email{holger.fehske}, \email{gerhard.wellein}\}\email{@fau.de}.}
\and Jonas Thies\thanks{Institute of Applied Mathematics, Delft University of Technology, Delft, The Netherlands. Email: \email{j.thies@tudelft.nl}.}
\and Georg Hager\footnotemark[2]
\and Holger Fehske\footnotemark[2]
\and Gerhard Wellein\footnotemark[2] \thanks{Department of Computer Science, Friedrich-Alexander-Universität Erlangen-Nürnberg.}}

\usepackage{amsopn}

\ifpdf
\hypersetup{
	pdftitle={\mytitle},
	pdfauthor={\myshortauthor}
}
\fi

\begin{document}

\maketitle

\begin{abstract}
Sparse linear iterative solvers are essential for many large-scale simulations.
Much of the runtime of these solvers is often spent in the implicit evaluation of matrix polynomials via a sequence of sparse matrix-vector products.
A variety of approaches has been proposed to make these polynomial evaluations explicit (i.e., fix the coefficients), e.g., polynomial preconditioners or $s$-step Krylov methods. Furthermore, it is nowadays a popular practice to approximate triangular solves by a matrix polynomial to increase parallelism.
Such algorithms allow to evaluate the polynomial using a so-called matrix power kernel (MPK), which computes the product between a power of a sparse matrix $A$ and a dense vector $x$, i.e., $A^px$, or a related operation.
Recently we have shown that using the level-based formulation of sparse matrix-vector multiplications in the Recursive Algebraic Coloring Engine (RACE) framework we can perform temporal cache blocking of MPK to increase its performance.
In this work, we demonstrate the application of this cache-blocking optimization in sparse iterative solvers.

By integrating the RACE library into the Trilinos framework, we demonstrate the speedups achieved in (preconditioned) $s$-step GMRES, polynomial preconditioners, and algebraic multigrid (AMG).
For MPK-dominated algorithms we achieve speedups of up to 3$\times$ on modern multi-core compute nodes.
For algorithms with moderate contributions from subspace orthogonalization, the gain reduces significantly, which is often caused by the insufficient quality of the orthogonalization
routines.
Finally, we showcase the application of RACE-accelerated solvers in a real-world wind turbine simulation (Nalu-Wind)  and highlight the new opportunities and perspectives opened up by RACE as a cache-blocking technique for MPK-enabled sparse solvers.
\end{abstract}



\section{Introduction and related work}

The solution of linear systems involving large sparse matrices 
is at the core of many computational workflows. Apart from application-specific approaches
like domain decomposition methods and geometric multigrid, the most popular classes of solvers  are Krylov subspace methods (often combined with preconditioning)
or algebraic multigrid. These algorithms are key components in open-source parallel simulation frameworks like Trilinos~\cite{trilinos-website}.

Krylov subspace methods perform a sequence of sparse matrix-vector multiplications (SpMV), vector updates (axpy) and inner products
to construct some basis of the Krylov subspace $\mathcal{K}_k(A,v)=\{v, Av, A^2v, \dots, A^{k-1}v\}$, and then extract an approximate solution by solving a
much smaller problem involving $A$ projected onto that subspace. In general, maintaining some orthogonality property of the basis is essential for stability, which leads to other vector operations being required in between SpMVs. 
If a preconditioner is used, the SpMVs may also be alternated with other operators, e.g., 
approximations of $A^{-1}$ or triangular factors $A^{-1}\approx U^{-1}L^{-1}$. Preconditioning is a broad field of research; for an overview of methods, see~\cite{wathen_2015}.
%
%
For sufficiently large matrices $A$, the SpMVs (including preconditioning) typically dominate the runtime, and it is known that these operations are main-memory bound for appropriately chosen sparse data layouts and may achieve high spatial locality when accessing the data elements of the matrix~\cite{KreutzerSELLC14}.

In the early days of parallel computing, $s$-step methods were developed to improve data locality (i.e., reduce communication) in Krylov methods~\cite{s_step_gmres_0, s_step_cg, s_step_gmres_1, s_step_gmres_2}.
They break up the data dependency by first computing a sequence of SpMVs and then using a sequence of scalar/vector operations to approximate the 
basis produced by, e.g., a Conjugate Gradient~(CG, \cite{s_step_cg}) or Generalized Minimum Residual~(GMRES~\cite{s_step_gmres_0, s_step_gmres_1, s_step_gmres_2}) method. 
These variants have recently received attention 
as they may use fast `kernels'  like the `Matrix-Power Kernel' (MPK) and the `Tall-Skinny QR' (TSQR), see~\cite{Demmel2008avoiding,hoemmen_thesis}. 
Recent work focuses on distributed-memory systems, i.e., reducing the number of messages and synchronization points in MPI implementations (e.g.,~\cite{Ichi_precon,7914608,Ichi_ca_pipe_gmres}).

However, the performance potential of the MPK for modern cache-based multicore architectures has not been exploited so far in any solver frameworks.
MPK involves the successive application of SpMV with the same matrix and offers the opportunity to exploit temporal locality by reusing the matrix elements from cache instead of repeatedly loading them from main memory.
For regular stencil algorithms it is well known how to improve temporal locality by temporal blocking~\cite{kaushik_2009};
on the other hand, for irregular sparse matrices such geometrical blocking approaches are generally not applicable.
Instead, an algebraic formulation of the problem needs to be considered to control the data dependencies and cache-access locality between successive SpMVs.
In~\cite{RACE_MPK} we have shown that this can be realized by a cache-aware traversal of the levels obtained from a breadth-first search (BFS) on the graph underlying the matrix.
Our implementation of the MPK achieves good scalability and high performance on modern multicore architectures for a broad range of matrices:
Compared to state-of-the-art implementations, RACE provides speedups in the range of 2--4$\times$.
We refer to~\cite{RACE_MPK} for an overview of related work on optimizing MPK.
Besides $s$-step Krylov algorithms there are other classes of methods
like polynomial preconditioning, smoothers in multigrid, Chebyshev time propagation, and power methods for eigenvalue solvers, which may also benefit from cache blocking of MPK.

\subsection*{Contributions}
In this paper we address the integration of the cache-blocked RACE MPK~\cite{RACE-git} into a number of representative iterative methods and evaluate the overall performance benefit on various solvers.
The cache-blocking strategy does not change the numerical behavior of the methods. 
Thus, the purpose of this paper is not to compare different iterative schemes or identify the most efficient preconditioners.
Instead, we focus on a broad range of numerical algorithms including several preconditioners and investigate the performance gains achieved through optimized MPK. 
Our specific contributions can be summarized as follows:
\begin{itemize}
\item Demonstration of the use of RACE MPK to accelerate $s$-step GMRES on modern multi-core CPUs,
\item incorporation of diagonal and triangular preconditioners into the MPK, where the triangular solves are approximated using Jacobi-Richardson iterations,
\item application of RACE MPK to GMRES polynomial preconditioning, demonstrating substantial performance improvements for high matrix powers,
\item introduction of strategies to accelerate algebraic multigrid (AMG) methods using RACE's cache blocking technique, and
\item showcasing the impact of highly efficient MPK on algorithmic choices using a case study from wind turbine simulation (Nalu-Wind~\cite{Nalu_Wind}).
\end{itemize}
In all cases a thorough performance analysis is conducted and the speedup obtained by RACE for different solvers is quantified.

\subsection*{Outline}
Throughout the paper we use the same representative hardware and matrices for demonstration purposes; these are introduced in Sec.~\ref{sec:testbed}.
We start by briefly recapitulating the idea of cache-blocking MPK using RACE in Sec.~\ref{sec:RACE_MPK}.
Section~\ref{sec:$s$-step_GMRES} discusses the hardware-efficient integration of RACE MPK into $s$-step GMRES methods. 
Section~\ref{sec:precon} addresses the integration of preconditioners into the MPK for $s$-step GMRES methods. We choose Jacobi and Gauss-Seidel sweeps as representative examples for diagonal and triangular preconditioners, where the triangular systems are solved approximately using Jacobi-Richardson iterations.
Further in the section we discuss the application of RACE to advanced polynomial and AMG preconditioners.
In Sec.~\ref{sec:case_study} we bring together the ideas developed in the paper to accelerate the solution of a momentum equation arising in the Nalu-Wind wind turbine simulation.
Finally, we summarize our findings in Sec.~\ref{sec:conclusion}.

\section{Hardware and software environment}
\label{sec:testbed}
\subsection{Hardware testbed}
\label{sec:hw_testbed}
The experiments presented in this paper were performed on single 
Intel Ice Lake (ICL) and AMD Epyc Rome (ROME) multicore processors.
These processors or similar ones are used in the majority of Top500~\cite{top500_list} systems today.
Key features of the chips are listed in Table~\ref{tab:testbed}.
Both architectures implement an x86 instruction set.
The 10\,nm ICL processor supports the AVX-512, 
while the 7\,nm ROME processor supports AVX2.
The systems are capable of sustaining more than 2\,\GHZ\ clock frequency and the turbo mode was active for all our experiments.
The AMD system has a higher core count (64 per socket) compared to its Intel counterpart (38 per socket).
Both systems have three levels of cache: private, inclusive
 L1 and L2 caches, and shared victim L3 cache.
The L3 cache on ICL is shared among all the cores within  a socket, 
while on ROME the L3 cache is shared within one core complex (CCX) unit comprising four cores.
Due to ROME's hierarchical ``chiplet'' design, the L3 cache is highly scalable 
and it can sustain an aggregate L3 load only bandwidth of 2700\,\GBS, 
while ICL achieves only 420\,\GBS.
The total L3 cache size of ROME is also much larger compared to ICL.
Both systems have eight-channel DDR memory and sustain  similar memory bandwidth.
Both are configured with one ccNUMA domain per socket configuration,
i.e., Sub-NUMA Clustering (SNC) was disabled on ICL and one NUMA node per socket (NPS1) mode was used on ROME.

\begin{table}[!tb]
	\begin{minipage}{0.47\linewidth}
		\centering
		\caption{Key specification of test bed machines.\label{tab:testbed}}
		\resizebox{\linewidth}{!}{%
			\begin{tabular}{l c c}
				Architecture   & ICL  & ROME  \\
				\midrule
				Chip Model         & Xeon Platinum 8368 & AMD EPYC 7662 \\
				Microarchitecture  &  Sunny Cove  & Zen-2 \\
				Cores per socket  &    38 & 64   \\
				Max. SIMD width      & 512\,\bits & 256\,\bits\\
				L1D cache capacity   & 38$\times$48\,\KiB   & 64$\times$32\,\KiB   \\
				L2 cache capacity    & 38$\times$1.25\,\MiB  & 64$\times$512\,\KiB     \\
				L3 cache capacity   &  57\,\MiB    & 16$\times$16\,\MiB \\
				L3 Bandwidth     & 420\,GB/s & 2700\,GB/s   \\
				Mem. Configuration  & 8 ch. DDR4-3200 &  8 ch. DDR4-3200   \\
				Mem. Bandwidth     & 170\,GB/s & 146\,GB/s   \\
			\end{tabular}
		} 
	\end{minipage}%
	\hfil
	\begin{minipage}{0.49\linewidth}
		\caption{Details of the benchmark matrices. See~\cite{UOF} for details. 
			\label{tab:matrices}}
		\resizebox{\linewidth}{!}{%
			\input{matrices.tex}

		}
	\end{minipage}
\end{table}

\subsection{Software environment}
The ICL system runs Red Hat Enterprise Linux (RHEL) version 8.4 while ROME runs Ubuntu 20.04.4 LTS.
For best performance, the OS setting ``Transparent Huge Pages'' (THP) was set to ``always'' on both the systems~\cite{Understanding}.
For compilation we used the Intel compiler version 2021.5.0 and 19.0.5 on ICL and ROME, respectively,
at the highest optimization level \texttt{-O3}.
Machine-specific code generation was employed via \texttt{-xHOST} on ICL and
\texttt{-march=core-avx2} \texttt{-mtune=core-avx2} on ROME.
All floating-point computations were performed in double precision, and integers were 32~\bits\ wide.
The linear solvers from the Trilinos framework~\cite{trilinos-website} used in this work were adapted to use the RACE MPK.
Both RACE and the modified Trilinos solvers are available through GitHub repositories; the exact versions used for the experiments are
available at \cite{RACE_git_commit} and \cite{Trilinos_git_commit}.
For BLAS computations,  Intel MKL~\cite{MKL} version 2022.0 was used on ICL.
On ROME we used MKL version 2020.0.4 unless otherwise stated.
It is well known that MKL sometimes exhibits low performance when it detects AMD hardware.
In order to make results comparable, we overwrite the \texttt{mkl\_serv\_intel\_cpu\_true} symbol with a function that always returns true (see~\cite{MKL_AMD_hack} for details).
On ROME we occasionally use the AOCL BLIS library~\cite{AOCL_BLIS, BLIS1} version 3.2 as an alternative. This is clearly indicated in the text.
Thread affinity was enforced by setting \texttt{OMP\_PLACES=cores} and \texttt{OMP\_PROC\_BIND=close}.
The run-to-run fluctuations in the experiments were less than 5\% and therefore we do not present any error bars.

\subsection{Benchmark matrices}
\label{sec:benchmark_matrices}
For our experiments we choose matrices from two prior publications~\cite{Jenifer_gmrespoly_precon_Trilinos,BergerVergiatGS2}
that are relevant to our work and are also available in the SuiteSparse Matrix Collection~\cite{UOF}.
We selected only square matrices with a memory footprint beyond 100\,\MB\, as our optimization targets big matrices that have to be loaded from main memory.
Table~\ref{tab:matrices} lists the matrices together with some relevant parameters:
number of rows (\NR), number of non-zeros (\NNZ), and average number of non-zeros per row (\NNZR).
In the following discussion we refer to the matrices by their IDs (first column of Table~\ref{tab:matrices}).
The compressed row storage (CRS) format was used throughout.


\section{Accelerating MPK using RACE}
\label{sec:RACE_MPK}
The central theme of this work revolves around speeding up various 
iterative solvers by using cache-blocked MPK.
MPK computes the application of powers of a sparse matrix to a dense vector. For a given sparse square  matrix $A$ and input vector $x$, 
MPK computes the matrix powers $A^px$ up to a maximum power of $p_m$, i.e., $p=1,\ldots,p_m$, and stores the result into $p_m$ vectors ($y_p=A^px$).
Usually this is done by performing back-to-back SpMVs as shown in Alg.~\ref{fig:baseline_alg}.
The input vector $x$ is stored in $y_0$ and each SpMV computation promotes the power by one.
Therefore, we reach $A^{p_m}x$ after $p_m$ SpMV computations.
\begin{algorithm}[tbp]
	\caption{Computing $A^{p_m}x$ using back-to-back SpMVs. 
		The arrays $val$, $col$, and $rowPtr$ hold the CRS data structure of $A$.
		The input and output vectors are stored in the $y$ matrix.}
	\label{fig:baseline_alg}
	\begin{algorithmic}
	        \STATE{$double$ $A.val[$\NNZ$]$} \textcolor{gray}{//store values of nonzeros in $A$}
		\STATE{$int$ $A.col[$\NNZ$]$} \textcolor{gray}{//column index of $A$}
		\STATE{$int$ $A.rowPtr[$\NR$+1]$} \textcolor{gray}{//row pointer of $A$}
		\STATE{$double:: y[$0:$p_m$, \NR$]$} \textcolor{gray}{//to store input and output vectors}
                \STATE{$y[0,:]=x[:]$} \textcolor{gray}{//starting vector $x$}
		\STATE \textcolor{gray}{//Perform $p_m$ SpMVs}
		\FOR{$p=1:p_m$}
		\STATE $y[p,:]$=SpMV($A$, $y[p-1,:]$) \textcolor{gray}{//$y_p = Ay_{p-1}$}
		\ENDFOR		
	\end{algorithmic}
\end{algorithm}

\begin{algorithm}[tbp]
	\caption{A prototype of SpMV callback function that can be passed to RACE for cache blocking MPK computation. The function is based on CRS data format.
	}
	\label{fig:callback_SpMV}
	\begin{algorithmic}
		\STATE \textbf{function} SpMV\_callback($int$ $row\_s$, $int$ $row\_e$, $int$ $p$, $arg\_type$ $kernel\_args$)
		\begin{ALC@g}
			\STATE{$A = kernel\_args.A$}
			\STATE{$y = kernel\_args.y$}
			\STATE \textcolor{gray}{//Loop over rows}
			\STATE{\textcolor{darkgray}{\#pragma omp parallel for schedule(static)}}
			\FOR{$row=row\_s:row\_e$}
			\STATE{$double$ $tmp=0$}
			\STATE \textcolor{gray}{//Loop over nonzeros in row}
			\FOR{$idx=(A.rowPtr[row]):(A.rowPtr[row+1]-1)$}
			\STATE{$tmp \mathrel{+}= A.val[idx]*y[p-1, A.col[idx]]$} 
			\ENDFOR
			\STATE{$y[p, row] = tmp$}
			\ENDFOR
		\end{ALC@g}
		\STATE \textbf{end function}	
	\end{algorithmic}
\end{algorithm}
SpMV being the central kernel in Alg.~\ref{fig:baseline_alg}, it is clear that its performance will be similar to that of SpMV.
For most of the matrices encountered in computational science and engineering, the latter is limited by main memory bandwidth on modern CPUs.
In Alg.~\ref{fig:baseline_alg}, if $A$ is larger than any cache it will be loaded $p_m$ times from memory. 
However, since MPK uses the same matrix $A$ for every SpMV, cache blocking of matrix accesses across successive matrix power calculations may reduce main memory traffic and thus improve performance.
This blocking is not straightforward due to dependencies among the SpMV computations, i.e., $y_p=A^{p}x=Ay_{p-1}$ depends on the results of the previous $y_{p-1}=A^{p-1}x=Ay_{p-2}$ computation; we denote this by $A^{p-1} \rightarrow A^{p}$.
The dependency, however, is not necessarily an all-to-all dependency, i.e., to calculate $A^{p}x$ on a subset of rows there is no need to finish the full $A^{p-1}x$ computation first.
The structure (or graph) of the matrix determines the dependencies between successive power calculations; 
the RACE library exploits it to enable cache blocking.
This is done by analyzing the graph of the matrix using a BFS traversal.
RACE  then uses the information from the level structure formed from the BFS 
and the cache size of the hardware to determine an execution order that enables cache blocking.
The routine to be blocked (here SpMV) is passed to RACE via a user-defined callback function.
The callback function takes the range of rows, the current power $p$, and any input required by the kernel as arguments.
In the execution phase, RACE supplies the values to these arguments and executes the
kernel in a cache-blocked manner according to an internally created execution order.
The user has to write a generic SpMV function computing $y_p=Ay_{p-1}$ on a range of rows.
Algorithm~\ref{fig:callback_SpMV} shows a prototype of such a function.
%
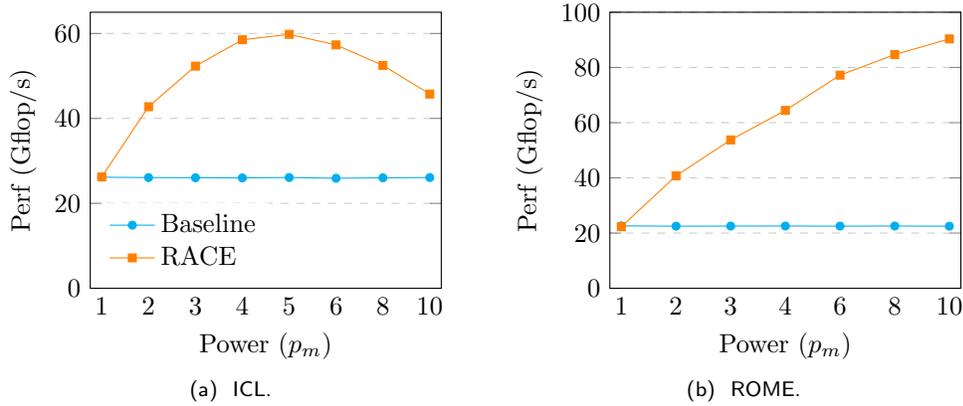
\begin{figure}[tbp]
	\centering
	\subfloat[ ICL.]{%
	\input{plots/perf_vs_power/plot_p_Flan_1565_ICL.tex}%
\label{fig:param_study:p:Flan:ICL}}
	\hfill
	\subfloat[ ROME.]{%
	\input{plots/perf_vs_power/plot_p_Flan_1565_ROME.tex}%
\label{fig:param_study:p:Flan:ROME}}
	\caption{\label{fig:perf_vs_power} Performance as a function of maximum power $p_m$
		for RACE and the baseline implementation of MPK.
		The experiment was conducted on \texttt{Flan\_1565} matrix
		and on ICL (a) and ROME (b).
		Figure reprinted from~\cite{RACE_MPK}.
	}
\end{figure}

More details on cache blocking via RACE can be found in~\cite{RACE_MPK}, where  
we have also shown that the level-based idea needs to be combined with a number of optimizations to achieve significant performance speedup (up to $5\times$) on MPK computations.
The maximum power $p_m$ has a substantial influence on the performance of the cache-blocked MPK.
Figure~\ref{fig:perf_vs_power}(a) shows the MPK performance of the \texttt{Flan\_1565} matrix as a function of $p_m$ on one socket of ICL compared with the baseline (non-blocked) kernel.
As a baseline for comparison we also show the performance of the naive kernel from Alg.~\ref{fig:baseline_alg}.
Both the RACE and baseline variants are parallelized using OpenMP~\cite{openmp}.
At $p_m=1$, both the variants are on par as expected. 
As $p_m$ increases, the RACE variant ideally needs to load the matrix only once from the 
main memory and the remaining $p_m-1$ accesses can be served from the caches.
Therefore, performance increases with $p_m$ (see Fig.~\ref{fig:perf_vs_power}(a)) until a maximum is reached.
Larger $p_m$ has a detrimental effect due to overhead from the blocking.
Hence, for maximum performance it is recommended to run the MPK 
kernel with the optimal power value, which we denote by $p_{\textrm{opt}}$.
If the required $p_m$ of an application is greater than $p_{\textrm{opt}}$
we execute multiple MPKs with $p_{\textrm{opt}}$ until $p_m$ is reached. 
Of course the last  MPK computation (the ``remainder loop'') might only operate up to a $p<p_{\textrm{opt}}$.
The $p_{\textrm{opt}}$ value depends on the matrix structure and the hardware and needs to be determined once for a given setting.
Comparing Fig.~\ref{fig:perf_vs_power}(a) and (b) demonstrates the qualitative impact of the hardware on $p_{\textrm{opt}}$. 
Due to its larger cache and massive cache bandwidth, ROME has a higher $p_{\textrm{opt}}$ value and achieves substantially better performance.

In the next sections we will discuss various applications of the MPK
and similar kernels in various iterative solvers.
Application-specific details and how to reformulate the algorithms 
to use RACE's cache blocking will be discussed.
Via thorough performance analysis we will observe the speedup achieved
by cache blocking and also detail some optimization strategies.


\section{$s$-step GMRES solver}
\label{sec:$s$-step_GMRES}
The GMRES method~\cite{GMRES} is a linear solver algorithm that computes an approximation $\hat{x}$ to the unknown solution $x$ of the linear system of equation $Ax=b$.
It is a Krylov subspace method and constructs $\hat{x}$ within the span of the Krylov subspace $\mathcal{K}_n(A,r) = \{r, Ar, A^2r, ..., A^{n-1}r\}$, where $r = b-Ax$ is the residual vector and $A^kr$ ($\forall k\in[0,n-1]$) are the Krylov vectors.
The size $n$ of the subspace is expanded iteratively to improve the approximation quality of $\hat{x}$.
Algorithm~\ref{alg:GMRES_pseudo}(b) shows the subspace generation routine of the GMRES solver.
Within each iteration, the algorithm generates a new vector $v[j+1]$ by performing a SpMV operation with the previous Krylov vector $v[j]$.
The newly generated vector is then orthonormalized against all previously generated Krylov basis vectors and added to the subspace.
Theoretically, the procedure can be repeated until the system converges.
However, within each iteration the memory and computational requirement grows as the subspace is expanded.
Therefore, the procedure is restarted every $m$ iterations. The parameter $m$ is commonly known as the restart length of the GMRES solver.
The pseudocode in Alg.~\ref{alg:GMRES_pseudo}(a) shows the wrapper around the subspace generation routine that restarts the GMRES solver after every $m$ iterations.


\begin{algorithm}[tbp]
	\centering
	\resizebox{0.49\linewidth}{!}{%
		\begin{minipage}{0.61\linewidth}
			\begin{algorithmic}[1]
				\STATE{$iter=0$}
				\WHILE{$iter<n$ \AND !converged}
				\STATE{$r$ = $b-Ax$}
				\STATE{$v[0:m-1]$ = \textcolor{ao(english)}{GenerateKrylovSupspace}($A$,$r$)}
				\STATE{$iter=iter+m$}
				\STATE{Find $x$ in the subspace $v[0:m-1]$}
				\ENDWHILE
			\end{algorithmic}
			\vspace{0.5em}
			\centering (a) ($s$-step)GMRES solver
			
			\vspace{0.5em}
			\begin{algorithmic}[1]
				\STATE \textbf{function} \textcolor{ao(english)}{GenerateKrylovSupspace}($A$,$r$)
				\begin{ALC@g}
					\STATE{$v[0]$ = $r/\lVert r\rVert$}
					\FOR{$j = 0:m-1$}			
					\STATE{$v[j+1]$ = SpMV($A$, $v[j]$)};
					\STATE{Orthonormalize $v[j+1]$ against $v[0:j]$}
					\STATE{\textcolor{gray}{//Check for convergence}}
					\STATE {converged = checkConvergence()}
					\IF{converged}
					\STATE{break}
					\ENDIF
					\ENDFOR
					\STATE {return $v[0:m-1]$}
				\end{ALC@g}
			\end{algorithmic}
			\vspace{0.5em}
			\centering (b) subspace generation routine of GMRES solver
		\end{minipage}
	}
	\resizebox{0.49\linewidth}{!}{%
		\begin{minipage}{0.61\linewidth}
			\begin{algorithmic}[1]
				\STATE \textbf{function} \textcolor{ao(english)}{GenerateKrylovSupspace}($A$,$r$)
				\begin{ALC@g}
					\STATE{$v[0]$ = $r/\lVert r\rVert$}
					\FOR{$j = 0:s:m-1$}
					\STATE{\textcolor{col3!50!black}{//MPK kernel}}
					\FOR{$p = 0:1:s-1$}
					\STATE{$v[j+p+1]$ = SpMV($A$, $v[j+p]$)};
					\ENDFOR
					\STATE{\textcolor{col2!50!black}{//BOrtho}}
					\STATE{Orthogonalize $v[j+1:j+s+1]$ against $v[0:j]$}
					\STATE{\textcolor{col1!50!black}{//TSQR}}
					\STATE{Orthonormalize vectors within $v[j+1:j+s+1]$}
					\STATE{\textcolor{gray}{//Check for convergence}}
					\STATE {converged = checkConvergence()}
					\IF{converged}
					\STATE{break}
					\ENDIF
					\ENDFOR
					\STATE {return $v[0:m-1]$}
				\end{ALC@g}
			\end{algorithmic}
			\vspace{0.5em}
			\centering (c) subspace generation routine of $s$-step GMRES solver
		\end{minipage}
	}
	\caption{Pseudocode of the GMRES and $s$-step GMRES solvers. (a) shows the general algorithmic structure of the solvers. 
		(b) and (c) show the specialized algorithms of the Krylov subspace generation routine for GMRES and $s$-step GMRES solvers, respectively. 
		\label{alg:GMRES_pseudo}}
\end{algorithm}

As the result of each SpMV operation is fed to the orthogonalization procedure,
this dependency prohibits the idea of calling a MPK, and thus the temporal blocking optimizations provided by RACE can not be applied to the subspace generation as presented in Alg.~\ref{alg:GMRES_pseudo}(b). 
However, the alternative $s$-step~\cite{s_step_gmres_1,s_step_gmres_2} formulation of GMRES
allows for MPK computations.
The basic structure of the $s$-step GMRES  solver remains the same as that of the standard GMRES solver (Alg.~\ref{alg:GMRES_pseudo}(a)). 
However, the Krylov subspace generation is modified to compute blocks of $s$ orthonormal vectors together (see  Alg.~\ref{alg:GMRES_pseudo}(c)). 
In the first step, the construction of the Krylov vectors in Alg.~\ref{alg:GMRES_pseudo}(c) is done as a sequence of $s$ back-to-back SpMV operations (line 5-7), which can be replaced by an MPK. 
The subsequent orthonormalization procedure is split into two routines: BOrtho and TSQR.
The BOrtho routine orthogonalizes the newly generated block of vectors with the previously
generated Krylov basis vectors, and TSQR orthonormalizes the vectors within the block.
The main advantage of the $s$-step variant is that it can use highly efficient BLAS kernels in the orthogonalization routines and
effectively reduce the frequency of MPI communications in distributed MPI-parallel setting due to the block-wise computations.
This results in a performance speedup over the standard GMRES solver\cite{Ichi_ca_GMRES_paper}.
Such implementations have therefore been called ``communication-avoiding GMRES'' (CA-GMRES) in the literature~\cite{mohiyudeen,hoemmen_thesis}.

\label{sec:$s$-step_GMRES_RACE_MPK}
We use the $s$-step GMRES method as 
implemented in the Belos package~\cite{Belos} of the Trilinos~\cite{trilinos-website} framework.
Belos performs back-to-back SpMVs as shown in Alg.~\ref{alg:GMRES_pseudo}(c) (lines 5--7) to generate the new Krylov vectors. 
This part is replaced with our cache-blocked MPK from RACE.
Note that, for stability reasons, the actual implementation of the $s$-step GMRES solver 
uses a Newton basis instead of the monomial basis~\cite{hoemmen_thesis}.
This means that the MPK routine computes $[v$, $(A-\lambda_1 I)v$, $(A-\lambda_2 I)^{2}v$, $(A-\lambda_3 I)^{3}v, \dots]$ instead of $[v$, $Av$, $A^2v$, $A^3v, \dots]$, where the $\lambda_i$ are just constant shifts. 
As the shifts only change the matrix diagonal, the RACE adaptation is straightforward and we pass the shifted SpMV callback function to RACE.
The $\lambda_i$ are computed from the eigenvalue information gathered by running a few steps of standard GMRES in Trilinos~\cite{Ichi_ca_GMRES_paper}.
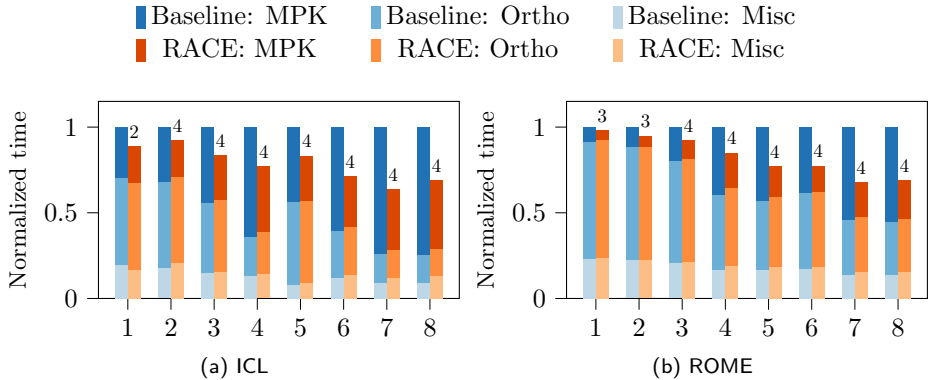
\begin{figure}[tbp]
	\centering
	%
	\input{plots//gmres_s_step_legend.tex}%

	\subfloat[ICL]{%
	\input{plots/horeka/gmres_s_step_none.tex}%
}
	\subfloat[ROME]{%
	\input{plots/tg097/gmres_s_step_none.tex}%
}
	\caption{Normalized execution time for 1000 iterations of the $s$-step GMRES method with (orange bars) and without (blue bars)  RACE MPK for the eight matrices (x-axis) shown in Table~\ref{tab:matrices} on ICL (a) and ROME (b).
	The absolute execution time is normalized to the baseline variant for each matrix separately.
	The stacked bar plot shows the time contributions of orthonormalization (Ortho) kernels, SpMV kernel 
		and other miscellaneous (Misc) routines.
	The numbers on top of the orange bars indicate the tuned power value $p_\textrm{opt}$ of RACE MPK operation.
		\label{fig:gmres_s_step_none}}
\end{figure}

Figure~\ref{fig:gmres_s_step_none} shows the performance advantage of 
the RACE-accelerated $s$-step GMRES solver on the ICL and ROME systems (see Sec.~\ref{sec:testbed} for details).
The numbers on the $x$-axis represent the matrix IDs from Table~\ref{tab:matrices}; the matrices are ordered by increasing size ($N_\mathrm{nz}$).
Unless mentioned otherwise, we set the cache size parameter $C$ of RACE to 85~\MB\ and 200~\MB\ for ICL and ROME, respectively.
The restart length $m$ of the solver was set  to 50.
Typically the step size $s$ of the $s$-step GMRES solver is kept under eight for stability reasons~\cite{hoemmen_thesis}. 
In our experiment we used $s=4$, which limits the maximum matrix power $p_m$.
The numbers above bars in Fig.~\ref{fig:gmres_s_step_none} denote the optimal power $p_\textrm{opt}$ at which RACE executed the kernel.
As this value is the same as $p_m$, i.e. $s$, for most matrices, an increase in $s$ will lead to higher speedups.
However, even with $s=4$ we manage to achieve a significant fraction of the maximum MPK speedup. 
This is in line with the discussion on Fig.~\ref{fig:perf_vs_power}, where we observe substantial performance gains already at low/moderate matrix power values.
On our test matrices, RACE accelerates the MPK computation by an average factor of 1.8$\times$ and 2.1$\times$ over the baseline method on ICL and ROME, respectively.
Note that this baseline uses the SpMV provided by the Trilinos package but modified by us to achieve a performance in line with the \rlm\  (see the discussion in Appendix~\ref{appendix:spmv_perf_trilinos} for details).
Otherwise the RACE MPK speedup would be even higher. 

Of course only the MPK routine is accelerated by RACE; the runtime of the other routines in the solver will be almost the same for both variants.
This reduces the average speedup of RACE for the complete solver to 1.3$\times$ and $1.2\times$ as seen in Figs.~\ref{fig:gmres_s_step_none}(a) and (b). 
On the other hand, the overall speedup stems purely from the performance gain, while numerically both $s$-step solver variants are identical. 
The reduced overall performance impact of RACE is mainly due to the significant cost of the orthogonalization procedure (Ortho).
In our experiments, one sweep of iterated classical Gram-Schmidt (ICGS) was performed in the BOrtho step of the Ortho routine (line 9 in Alg.~\ref{alg:GMRES_pseudo}(c)) and tall-skinny QR decomposition was used in the TSQR step (line 11).
Both of these routines are accounted for in the Ortho time. 
Note that for robustness it is advisable 
to perform two or more sweeps of ICGS, but for all examples in this paper the above choice of one sweep proved to be sufficiently stable due to the relatively short restart length.
The cost of Ortho is especially high for matrices with low \NNZR\ values because here the runtime complexity for both SpMV and Ortho approaches  $\mathcal{O(\NR)}$.

Although RACE MPK attains higher speedup on ROME compared to ICL, the total solver speedup on ROME is lower than on ICL.
Again the Ortho routine is the culprit as it takes substantially longer on ROME for the same matrices.
An in-depth performance analysis revealed that the BLAS calls associated with the Ortho routines performed poorly on ROME (see Apendix~\ref{appendix:ortho_perf_trilinos}).
Overall it could be said that the speedup of the $s$-step GMRES solver achieved by cache blocking is limited by the Ortho routines.
As a result, implementing cache blocking on Krylov methods with short recurrence such as the $s$-step conjugate gradients (CG), where the Ortho cost is minimal, might lead to a solver speedup that approaches the RACE MPK speedup.
However, Trilinos currently does not have an $s$-step CG implementation and we would perform this analysis in the future once the solver becomes available.

The runtime contribution of the remaining kernels (apart from MPK and Ortho) is minor (indicated as ``Misc'' in Fig.~\ref{fig:gmres_s_step_none}).
Typically this contribution is slightly higher in the RACE variant because it includes the pre-processing cost of RACE (usually 30--50 SpMVs).
However, the total number of solver iterations is very large in most of the applications, and the extra cost can be easily amortized.

\section{Preconditioners}
\label{sec:precon}
A GMRES solver is rarely used without a preconditioner because its long recurrence makes it infeasible if the number of iterations increases.
Restarting, on the other hand, may drastically increase the total number of iterations and even cause stagnation.
A preconditioner transforms the linear system to an equivalent system by formally multiplying the system matrix with another linear operator from the
left or right, or both. The goal of this transformation is to improve the condition number of the overall operator, or more specifically to decrease the number of iterations required to solve the system by GMRES.
Throughout the paper we will apply the preconditioner from the right, but other choices can be implemented analogously.
Hence, the system $Ax=b$ is transformed to $AP^{-1}y=b$, where $y=Px$ is solved for $x$ by applying the preconditioner a final time in the end.
The preconditioner $P^{-1}$ is chosen to be some approximation of $A^{-1}$ which is cheap to construct and to apply to a vector. 
In practice, we then perform the operation $AP^{-1}v$ instead of the SpMV routine computing $Av$.

In case of $s$-step GMRES, the introduction of a preconditioner requires us to compute 
the vectors $[v$, $AP^{-1}v$, $(AP^{-1})^{2}v$, $(AP^{-1})^{3}v$, $\dots]$ instead of $[v$, $Av$, $A^2v$, $A^3v, \dots]$ in the MPK.
The main challenge here is that the preconditioner results in an additional dependency between 
$P^{-1}$ and $A$, and, 
using the dependency notation introduced in Sec.~\ref{sec:RACE_MPK},
we can denote the MPK dependencies as  $P^{-1} \rightarrow  AP^{-1} \rightarrow P^{-1}AP^{-1} \rightarrow (AP^{-1})^2 \dots$.
We will show that, despite these additional dependencies, it is possible to
cache block the preconditioned $s$-step GMRES using RACE and achieve significant speedups on modern multicore CPUs.

\subsection{Relaxation preconditioners}
\label{sec:precon_relax}
Relaxation preconditioners use iterations from a stationary iterative splitting method.
In the following we will investigate two popular choices in this category, Jacobi and Gauss-Seidel, 
and demonstrate how RACE can be used to accelerate the preconditioned $s$-step GMRES solvers.

\subsubsection{Jacobi}
\label{sec:Jacobi}

The application of a Jacobi preconditioner to a vector $v$
follows the Jacobi iteration:
\begin{equation}
 \label{eq:Jacobi}
 z^{k+1} = D^{-1}v -D^{-1}(L+U)z^{k}\eos
\end{equation}
Here $z^{k+1}$ and $z^{k}$ denote the new and old iterate of $P^{-1}v$.
The matrices $L$ and $U$ are the strictly lower and upper triangular part of matrix $A$ and matrix $D$ is the diagonal.
In many use cases only one Jacobi iteration is applied for the preconditioner. 
If the initial guess $z^0$ is also chosen to be zero, the entire Jacobi preconditioner simplifies to $z^1 = D^{-1}v$,  which is just a diagonal scaling of the input vector $v$.
Consequently this type of Jacobi preconditioner is also commonly known as diagonal preconditioner.

\begin{figure}[tbp]
	\centering
	%
	\input{plots//gmres_s_step_legend_w_precon.tex}%

	\subfloat[ICL]{%
	\input{plots/horeka/gmres_s_step_Jacobi.tex}%
}
	\subfloat[ROME]{%
	\input{plots/tg097/gmres_s_step_Jacobi.tex}%
}
	\caption{Time taken by baseline and RACE accelerated variant of  $s$-step GMRES solver using Jacobi preconditioner. 
		The stacked bar plot displays the time contributions by orthonormalization (Ortho) kernels, SpMV kernel 
		and other small miscellaneous (Misc) routines.
		\label{fig:gmres_s_step_Jacobi}}
\end{figure}
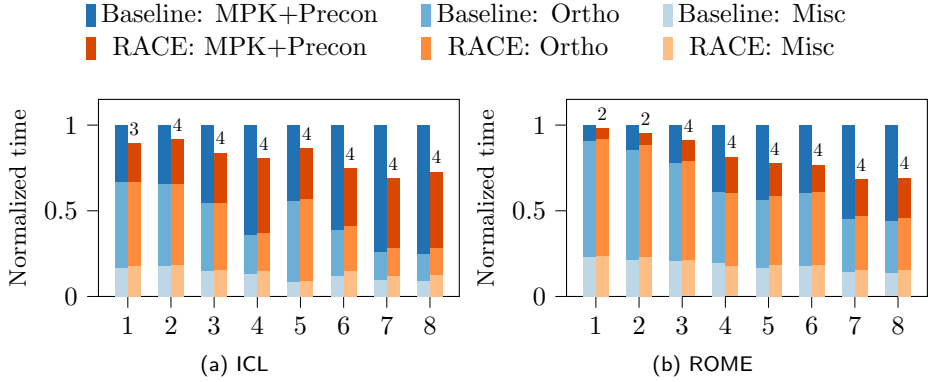

An advantage of the diagonal preconditioner is that the diagonal scaling $P^{-1} = D^{-1}$ does not introduce any additional dependency between $D^{-1}$ and $A$ 
and therefore $AD^{-1}$ can be fused to a single kernel.
Hence, the actual MPK dependency shown in Sec.~\ref{sec:precon} boils down to  
$AD^{-1} \rightarrow  (AD^{-1})^2 \rightarrow (AD^{-1})^3 \dots$
and is analogous to the one seen for the plain MPK routine without preconditioners in Sec.~\ref{sec:RACE_MPK}.
The diagonally preconditioned SpMV routine computing $AD^{-1}v$ is straightforward
and similar to SpMV as shown in Alg.~\ref{fig:callback_SpMV}, except that it requires an extra diagonal scaling along the columns.
The baseline $s$-step GMRES solver using a Jacobi preconditioner calls this routine $s$ times on the whole matrix to compute the MPK.
The RACE variant, however, blocks the matrices $A$ and $D^{-1}$ in cache across the $s$ iterations.
To achieve this, we pass the diagonally scaled SpMV callback routine to RACE, which then performs 
cache blocking based on the internally created execution order, similar to the plain unpreconditioned MPK computations seen in Sec.~\ref{sec:$s$-step_GMRES_RACE_MPK}.
Due to the similarities between the plain MPK and the Jacobi-preconditioned
MPK, the performance characteristics of $s$-step GMRES solver remain unchanged (compare Fig.~\ref{fig:gmres_s_step_none} and Fig.~\ref{fig:gmres_s_step_Jacobi}).
In this case, too, the RACE-accelerated solver achieves an average speedup of almost 1.25$\times$ (see Fig.~\ref{fig:gmres_s_step_Jacobi}). 

In practice, more than one Jacobi iteration is rarely used as a preconditioner as it requires additional SpMVs ($z^k\neq 0$ in \eqref{eq:Jacobi}). In fact, performing $k$ Jacobi iterations is equivalent to a simple matrix polynomial preconditioner based on the Neumann series $(I-B)^{-1} \approx \sum_{j=0}^k B^k$ for $B=-D^{-1}(L+U)$.
The cache-blocking approach in RACE MPK offers the opportunity to reduce the computational cost of these additional SpMVs to the point that this approach may be competitive. Cache-blocked polynomial preconditioners can even  accelerate a standard Krylov method,
as we will show in Sec.~\ref{sec:case_study}.
Combining it with an $s$-step method allows to cache block for higher powers, which may be even more efficient.

\subsubsection{Gauss-Seidel}
\label{sec:GS2_precon}
The Gauss-Seidel (GS) preconditioner is derived from the GS iteration:
\begin{equation}
	\label{eq:GS}
	(L+D)z^{k+1} = v - Uz^{k}\eos
\end{equation}
For many linear systems, GS is considered to be superior to
Jacobi since it uses the new iterate $z^{k+1}$ whenever available ($L$ is applied to $z^{k+1}$ in~\eqref{eq:GS}).
However, shared-memory parallelization is a challenge as
a thread does not know when other threads have updated their $z$ entries.
Two well-known solutions to this triangular solver problem are multicoloring~\cite{gs_mc_EVANS} and level scheduling~\cite{level_scheduling_SAAD}.
Reordering via multicoloring often degrades the data locality and the convergence rate, resulting in performance loss.
On the other hand, level scheduling maintains the convergence rate but often exhibits limited parallelism.
Another promising solution, especially for preconditioners, is
the two-stage Gauss-Seidel (GS2) iteration~\cite{GS2_original}.
Here a fixed number of Jacobi-Richardson iterations is used to solve~\eqref{eq:GS}.
This means that, within each iteration of GS, we have inner iterations of Jacobi-Richardson.
The benefit with this approach is that we can solve the system using simple SpMVs and BLAS-1 operations.
This technique has been used to increase parallelism for GPUs~\cite{Chow_ILU_with_JR, BergerVergiatGS2}, but not for cache blocking on multi-core CPUs.
Of course, in contrast to level scheduling, the system is not solved exactly with the Jacobi-Richardson iterations;
however, it was shown in~\cite{BergerVergiatGS2} that for preconditioners, where $A^{-1}$ is already approximated, 
this method produces similar convergence rates for many matrices.
\begin{algorithm}[tbp]
	\resizebox{0.49\linewidth}{!}{%
	\begin{minipage}{0.5\linewidth}
		\vspace{1em}
		\begin{algorithmic}[1]
			\STATE{\textcolor{white!50!black}{//$z^0$ is the initial guess}}
			\FOR{$k = 0:K$}
			\STATE{$g^{0}_k = D^{-1}(v - Uz^k)$}
			\FOR{$j = 1:\gamma$}
			\STATE{$g^{k}_{j} = g^{k}_0 - D^{-1}Lg^{k}_{j-1}$};
			\ENDFOR
			\STATE{\textcolor{white!50!black}{//update}}
			\STATE{$z^{k+1} = z^k + g^{k}_{\gamma}$}
			\ENDFOR
		\end{algorithmic}
	\centering
	\vspace{0.5em}
	\text{{\small(a) GS2 pseudocode}}
	\end{minipage}
}
	\resizebox{0.49\linewidth}{!}{%
	\begin{minipage}{0.5\linewidth}
	\begin{algorithmic}[1]
		\IF{$j$ == 0}
			\STATE{$g^{0}_0 = D^{-1}(v[p] - Uz^0)$}			
		\ELSIF{$j$ == 1}
			\STATE{$g^{0}_j = g^{0}_0 - D^{-1}Lg^{0}_{j-1}$}
		\ELSIF{$j$ == 2}
			\STATE{$g^{0}_j = g^{0}_0 - D^{-1}Lg^{0}_{j-1}$}
			\STATE{$z^{1} = z^{0} + g^{0}_j$}
		\ELSIF{$j$ == 3}
			\STATE{$v[p+1] = Az^{1}$}
		\ENDIF
	\end{algorithmic}
	\text{\small{(b) MPK with GS2 of $\gamma=2$}}
	\end{minipage}
	}
\caption{(a) Pseudocode of a two-stage Gauss-Seidel (GS2) iteration with $\gamma$ inner Jacobi-Richardson iterations. 
The $k$-loop is the outer Gauss-Seidel iteration and $j$ is the inner iteration.
(b) Unrolled implementation of MPK with GS2 preconditioner ($\gamma=2$) passed to RACE for cache blocking.
The implementation takes an input vector $v[p]$ and performs the computation $v[p+1] = AP^{-1}v[p]$, where $P^{-1}$ is the GS2 preconditioner.
\label{alg:GS2}}
\end{algorithm}

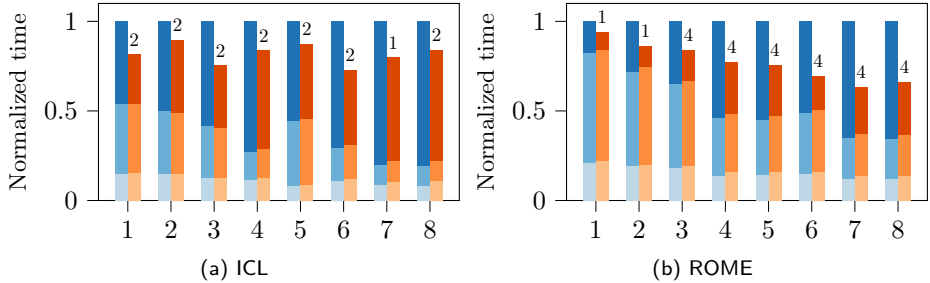
\begin{figure}[tbp]
	\centering
	\subfloat[ICL]{%
	\input{plots/horeka/gmres_s_step_GS2_1.tex}%
}
	\subfloat[ROME]{%
	\input{plots/tg097/gmres_s_step_GS2_1.tex}%
}
	\caption{Time taken by baseline and RACE-accelerated variant of the $s$-step GMRES solver using the GS2 preconditioner with one inner Jacobi-Richardson iteration.
	The legend is the same as in Fig.~\ref{fig:gmres_s_step_Jacobi}.
		\label{fig:gmres_s_step_GS2_1}}
\end{figure}
A GS2 iteration algorithm based on the non-compact form of the GS iteration (see~\cite{BergerVergiatGS2} for details) is shown in Alg.~\ref{alg:GS2}(a).
As in the Jacobi case, typically only one outer iteration of GS2 is employed as a preconditioner and frequently combined 
with one ($\gamma=1$) or two ($\gamma=2$) inner Jacobi-Richardson iterations.
The cache-blocking of GS2 preconditioner with RACE is more involved than
the Jacobi counterpart.
Here, the preconditioner $P^{-1}$ itself involves many interdependent steps.
The dependency within $P^{-1}$ of GS2 with $\gamma=2$
can be expressed as $U \rightarrow LU \rightarrow L^{2}U$.
Finally, when applying the preconditioner to the matrix in the $s$-step GMRES solver, there is an additional dependency  $P^{-1} \rightarrow A$.
This means that even in a single step ($s=1$) of the $s$-step GMRES solver we have a
dependency chain of length four.
To allow for an easy integration of RACE into the such preconditioners, we 
divide each power computation into a fixed number of \emph{sub-powers}.
In case of GS2 as shown above we would have four sub-powers within a power,
and the power loop in RACE will map to the power loop along the MPK computations of the $s$-step GMRES solver (e.g., line 5 in Alg.~\ref{alg:GMRES_pseudo}(c)).
Algorithm~\ref{alg:GS2}(b) shows the MPK routine 
with the GS2 ($\gamma=2$) preconditioner  that is passed to RACE for cache blocking.
Here we distinguish each stage of the dependency chain using the sub-power ($j$ in Alg.~\ref{alg:GS2}(b)),
which goes from zero to $\gamma+1$ (three in this case) for each power computation.
The first three sub-powers compute the application of the GS2 preconditioner on a vector $v$,
and the result is stored in vector $z^1$.
The last sub-power, i.e., at $j=3$ in case of Alg.~\ref{alg:GS2}(b), 
calculates $Az^1$.
Note that ``sub-power'' is just a convenient abstraction 
on top of the power loop in RACE; it satisfies all the BFS level dependencies mentioned in Sec.~\ref{sec:RACE_MPK}, ensuring that
 the dependencies within $P^{-1}$ (i.e., $U \rightarrow LU \rightarrow L^{2}U$) are met.
Of course in this case the callback function passed to RACE will have an extra input argument for the sub-power $j$, 
which will be imported by RACE internally during the execution phase.
The levels in RACE are still generated from matrix $A$ only, as
the matrices $U$ and $L$ have a subset of the sparsity pattern of $A$.

In contrast to the plain unpreconditioned or Jacobi-preconditioned $s$-step GMRES solver,
a GS2 preconditioned solver has the benefit that even with $s=1$ (normal GMRES solver) some performance advantage is possible since we reuse the matrices within $P^{-1}$ and between $P^{-1}$ and $A$.
For example, with two inner iterations ($\gamma=2$) a straightforward implementation 
requires to load the matrices $A$ and $U$ once and $L$ twice.
However, with the cache-blocked variant ideally we need to load the matrix $L$ only once.
We can further save the traffic from matrix $A$ if we perform the SpMV with a split form, i.e., $A=L+D+U$, leading to reuse in matrices $L$ and $U$. 
This optimization is applied in our RACE implementation.
For $s$-step GMRES solvers with $s>1$, this benefit adds to the advantage of blocking the matrices to higher powers.

The GS2 preconditioner is implemented in the Ifpack2 package~\cite{Ifpack2} of Trilinos.
In this section, we use this preconditioner for the $s$-step GMRES solver as a baseline for comparison.
Similar to the Jacobi preconditioner, it is common practice to choose the starting vector $z^0$ to be zero.
This allows for short-circuiting some computations in the GS2 iteration, i.e.,
line 3 in Alg.~\ref{alg:GS2}(a) simplifies to $g^{k}_0 = D^{-1}v$ and the update step (line 8) to $z^{k+1} = g^{k}_{\gamma}$.
The Ifpack2 implementation currently initializes the vector with zero by default but it does not short circuit the computations although the Kokkos backend supports this.
This mainly results in an additional overhead of half an SpMV ($Uz^{k}$).
In the interest of a fair comparison we modified the Ifpack2 code to allow for short-circuiting the unnecessary computations as well. 

Figure~\ref{fig:gmres_s_step_GS2_1} shows the performance of the GS2-preconditioned $s$-step GMRES solver
with $\gamma$=1  on ICL and ROME.
On ICL, the power value at which RACE operates is lower compared
to the previously discussed $s$-step solvers because
each power step ($p$) contains multiple sub-power computations.
Thus, the effective total power to which RACE applies cache blocking is the 
product of the two power computations and the maximum performance is achieved at lower $p_\textrm{opt}$ (see discussion of Fig.~\ref{fig:perf_vs_power}(a)).
However, on ROME most matrices reach the maximum power value of four ($p_\textrm{opt}=s$) due to its large cache size. 
Increasing the $\gamma$ from one to two improves the speedup slightly 
from 1.22$\times$ to 1.3$\times$ on ICL (not shown in figure), 
as more reuse can be applied within the inner iterations. 

The results for the GS2 preconditioner demonstrate the applicability of RACE to a wider range of preconditioners that require chaining of multiple routines.
This includes sparse approximate inverse preconditioners were the $P^{-1}$ is explicitly computed 
and can be chained effectively with the matrix $A$ to implement the $s$-step GMRES solver.
Similarly, factorization-based preconditioners like ILU can be implemented with RACE using the chaining idea and applying Jacobi-Richardson iterations to solve the triangular systems~\cite{Chow_ILU_with_JR}.
In this paper we do not investigate further on this class of preconditioners
but rather demonstrate the applicability of RACE to two different classes of preconditioners which show significant performance improvement even on standard ($s=1$) GMRES solvers.

%
%
%
%

\subsection{Polynomial preconditioners}
\label{sec:gmres_poly_precon}
A polynomial preconditioner has the form $P^{-1} = \cal{P}(A)$, where $\cal{P}$ is a polynomial.
Such preconditioners have been extensively studied in the context of Krylov-based solvers~\cite{cg_poly_precon,saad_least_square_poly_precon}.
Stability and 
setup costs of polynomial preconditioners, such as extreme eigenvalue calculations, were major concerns for a long time. 
However, many recent studies, e.g.~\cite{Jenifer_gmrespoly_precon,Saad_poly_precon}, have stimulated renewed interest in these methods.  
In particular, the lack of global communication in the evaluation of the polynomial make them attractive for large-scale computing.

In each preconditioning step, a polynomial of degree $d$ is applied. The application of $P^{-1}$ to a vector $x$ can be expressed as
\begin{equation}
P^{-1}x = {\cal P}(A)x = \lambda_0x + \lambda_1 Ax + \lambda_2 A^2x + \dots +  \lambda_d A^dx\cma
\label{eq:poly_precon}
\end{equation}
where the $\{\lambda_i\}$ are scalar coefficients that determine the type of polynomial, with Chebyshev~\cite{cg_poly_precon} and GMRES polynomials~\cite{gmres_precon_orig} being the most popular ones.
Here we focus on the GMRES polynomial,
where the scalar constants are generated by running $d$ iterations of a GMRES solver in a pre-processing step.

The optimal degree $d$ for GMRES polynomial preconditioners is rather high (in the range of $40$--$100$), thus applying the preconditioner requires many back-to-back SpMVs, making this operation frequently the dominant part of the solver.\footnote{Note that the high computational cost of the preconditioner is often amortized by a decrease in the total number of iterations, making the preconditioner effective.}
This is a very attractive scenario for the MPK of RACE as potentially large speedups are achievable and it directly applies to the hotspot of the solver.
Thus, we exclusively focus on accelerating the polynomial preconditioner, which is called in each iteration of the GMRES solver. In the following the parameters $C$ and $m$ remain similar to the previous experiments with the $s$-step GMRES solver  (see Sec.~\ref{sec:$s$-step_GMRES_RACE_MPK}).
We use the Belos package of Trilinos as the implementation baseline (see~\cite{Jenifer_gmrespoly_precon_Trilinos} for details)
and choose a polynomial of degree 80.

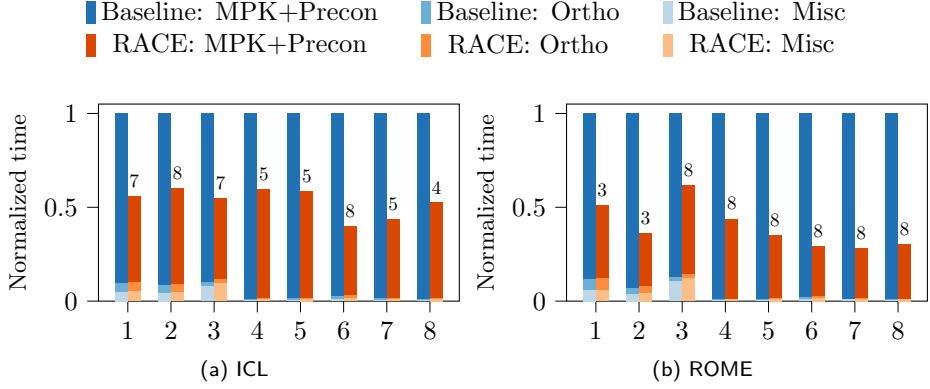
\begin{figure}[tbp]
	\centering
	%
	\input{plots//gmres_s_step_legend_w_precon.tex}%

	\subfloat[ICL]{%
	\input{plots/horeka/gmres_poly_none.tex}%
}
	\subfloat[ROME]{%
	\input{plots/tg097/gmres_poly_none.tex}%
}
	\caption{Time taken by the baseline and RACE-accelerated variants of the GMRES solver using a polynomial preconditioner of degree 80.
		\label{fig:gmres_poly_none}}
\end{figure}

Figure~\ref{fig:gmres_poly_none} shows the performance benefit when using RACE 
to accelerate the polynomial preconditioner in the GMRES solver.
The striking observation is that the speedup obtained by RACE
is significantly higher than previously observed with $s$-step GMRES solvers.
There are two reasons for this:
First, the polynomial application $P^{-1}x$ consumes a significant fraction of the entire solver runtime (more than 95\% for most matrices; see Fig.~\ref{fig:gmres_poly_none}).
Second, high powers in the MPK computations (80 in our case) allow RACE to block for higher power values and thus operate at high performance levels. 
This is also the reason why most of the matrices on ROME operates at $p=8$, where eight is the highest power value in our tuning space of $p \in [1:8]$. 
Higher $p$ did not prove to be significantly faster.
On ICL, the optimal power values for RACE are lower due to the smaller cache size
and therefore most of the matrices have $p_\textrm{opt}<8$.
Overall, RACE improves the MPK performance on ROME (ICL) by an average factor of almost 3$\times$ (2$\times$), which translates to an average 2.7$\times$ (1.9$\times$) speedup on the entire GMRES solver.

\begin{figure}[tbp]
	\centering
	\subfloat[ICL]{%
	\input{plots/horeka/gmres_poly_GS2_1.tex}%
}
	\subfloat[ROME]{%
	\input{plots/tg097/gmres_poly_GS2_1.tex}%
}
	\caption{Time taken by baseline and RACE-accelerated variants of the GMRES solver using GS2 with one inner iteration on top of the polynomial preconditioner of degree 80.
	Colors have the same meaning as in Fig.~\ref{fig:gmres_poly_none}.	
		\label{fig:gmres_poly_GS2_1}}
\end{figure}
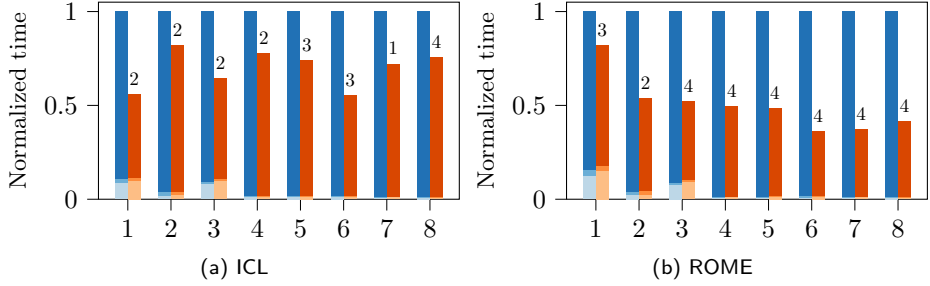

GMRES polynomial preconditioners can be further combined with other preconditioners. 
In such scenarios, the polynomial is formulated in terms of the combined matrix $AM^{-1}$, where $M^{-1}$ is another preconditioner.  
The application of the preconditioner then reads:
\begin{equation}
	P^{-1}x = {\cal P}(AM^{-1})x = \lambda_0x + \lambda_1 AM^{-1}x + \lambda_2 (AM^{-1})^2x + \dots 
	\label{eq:poly_precon_w_precon}
\end{equation} 
The dependencies caused by the new matrix $M^{-1}$ have to be taken into account in the cache blocking, and a similar approach to the one described in Sec.~\ref{sec:precon_relax} is in order.
In case of a Jacobi preconditioner, only matrix diagonal scaling is required and we attain almost the same speedup as with plain polynomial preconditioning as shown in Fig.~\ref{fig:gmres_poly_none}.
The speedups obtained by RACE when using the GS2 ($\gamma=1$) preconditioner
on top of the polynomial preconditioner are shown in Fig.\ref{fig:gmres_poly_GS2_1}.
Although the computational kernels remain similar to the
ones discussed in Sec.~\ref{sec:GS2_precon} above, the speedup
is much higher (1.5$\times$ and 2.1$\times$ on ICL and ROME) in the case of polynomial preconditioners due to the two reasons discussed in the previous paragraph.


\subsection{Algebraic multigrid preconditioners}
\label{sec:amg}
Algebraic Multigrid (AMG) preconditioners are among the most widely used
preconditioners for Krylov solvers~\cite{wathen_2015}.
AMG preconditioners are particularly effective for solving large 3D problems.
Similar to any multilevel or geometrical multigrid schemes~\cite{MG_ref}, the AMG preconditioner
uses a hierarchy of grids with various refinement (discretization) levels.
Inter-grid transfer operators are used to transfer information between the grid levels.
The restriction and prolongation operators are the two inter-grid transfer operators used to 
transfer information from a fine to coarse grid and vice versa.
Within each grid level a smoothing operator is applied to reduce the error within the level.
In contrast to geometrical multigrid, AMG algebraically determines the 
coarse grids and the inter-grid  transfer operators~\cite{AMG_falgout} based on the matrix entries.
As AMG does not require any explicit knowledge of the problem geometry,
it is particularly useful for problems having a complicated or even unknown geometry.

\begin{algorithm}
	\begin{algorithmic}[1]
	\STATE{\textcolor{darkgray}{//Solve $Az=v$}}
	\STATE{$z = 0$}
	\STATE{AMG($A$,$v$,$z$,0)} 	
	\STATE \textbf{function} \textcolor{ao(english)}{AMG}($A_k$,$b$,$x$,$k$)
	\begin{ALC@g}
		\STATE{$x = S^{\textrm{pre}}_k(A_k, b, x)$ \textcolor{darkgray}{//Pre-smoothing}}
		\IF{$k \ne max\_levels-1$}
			\STATE{$r_k = b-A_kx$ \textcolor{darkgray}{//Residual}}
			\STATE{$r_{k+1} = R_k(r_k)$ \textcolor{darkgray}{//Restriction on the residual}}
			\STATE{$c_{k+1} = 0$}
			\STATE{\textcolor{darkgray}{//Call AMG with next coarser matrix $A_{k+1}$}}
			\STATE{AMG($A_{k+1}$, $r_{k+1}$, $c_{k+1}$, $k+1$)}
			\STATE{$c_k = P_kc_{k+1}$ \textcolor{darkgray}{//Prolongation on the correction vector}}
			\STATE{$x=x+c_k$ \textcolor{darkgray}{//Add correction}}
			\STATE{$x = S^{\textrm{post}}_k(A_k, b, x$)  \textcolor{darkgray}{//Post-smoothing}}
		\ENDIF
	\end{ALC@g}
		
	\end{algorithmic}
	\caption{Pseudocode of a single AMG V-cycle, adapted from~\cite{AMG_NALU_Thomas}. The letter $k$ denotes the grid level.
		\label{alg:MG}}
\end{algorithm}

A single iteration of AMG starts with a smoothing operation (pre-smoothing) performed on
the finest grid level with an initial guess of zero.
The residual is then calculated on the finest level and transferred to the next coarser level
using the restriction operator.
Now the coarse level performs the smoothing and recursively applies the procedure till
the coarsest level is reached.
The linear system on the coarsest level is usually solved by a direct solver.
This solution then serves as a correction on the next finer level and is  transferred using the prolongation operator.
The smoothing operation is again performed (post-smoothing) on the finer level using the correction as the initial guess.
The grid then transfers the solution to next finer level and the process repeats until we reach the finest level.
Due to the manner in which the grids are traversed, i.e., finest to coarsest and then back to finest, this is called \emph{V-cycle AMG}.
Although there are many other types of cycles, we will concentrate on the V-cycle in this paper.
When using AMG as a preconditioner, a single V-cycle of AMG is typically used to compute $P^{-1}v$.
Algorithm~\ref{alg:MG} shows the corresponding high-level algorithm computing $z=P^{-1}v$, where $P^{-1}$ is an approximation to $A^{-1}$.

\begin{figure}[tbp]
	\centering
	%
	\input{plots//gmres_mg_legend.tex}%

	\subfloat[ICL, 2 outer sweeps]{%
	\input{plots/horeka/gmres_amg_gs2_1_2.tex}%
}
	\subfloat[ROME, 2 outer sweeps]{%
	\input{plots/tg097/gmres_amg_gs2_1_2.tex}%
}
	
	\subfloat[ICL, avg. speedup]{%
	\input{plots/horeka/gmres_amg_gs2_avg_speedup.tex}%
}
	\subfloat[ROME, avg. speedup]{%
	\input{plots/tg097/gmres_amg_gs2_avg_speedup.tex}%
}
	\caption{\label{fig:gmres_amg_gs2}
	(a), (b) Comparison of time taken by the baseline and the RACE-accelerated variant of the GMRES solver preconditioned by the algebraic multigrid preconditioner using the GS2 smoother with one inner iteration ($\gamma$=1) and two outer sweeps.
	(c), (d) Improvement in speedup (averaged across the eight benchmark matrices) as the number of outer sweeps is increased on ICL and ROME.
}
\end{figure}
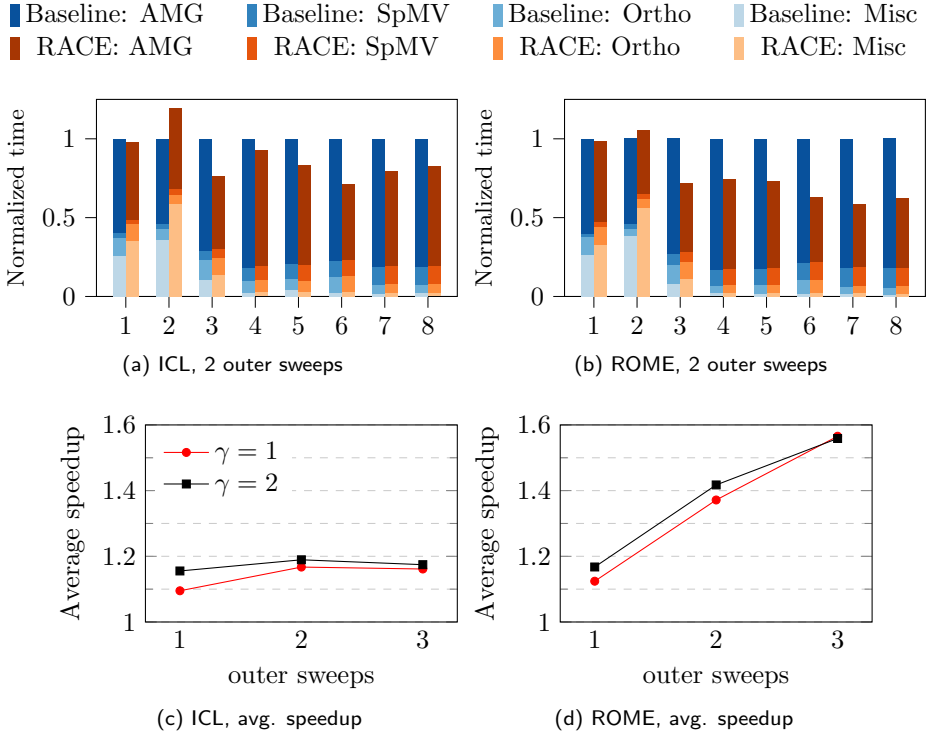

\begin{figure}[tbp]
	\centering
	\subfloat[ICL, degree 3]{%
	\input{plots/horeka/gmres_amg_cheb3.tex}%
}
	\subfloat[ROME, degree 3]{%
	\input{plots/tg097/gmres_amg_cheb3.tex}%
}
	
	\subfloat[time contribution of level 0]{%
	\input{plots/gmres_amg_cheb_sweep3_level0_contrib.tex}%
}
	\subfloat[avg. speedup]{%
	\input{plots/gmres_amg_cheb_avg_speedup.tex}%
}
	\caption{\label{fig:gmres_amg_cheb}
		(a), (b) Comparison of time taken by the baseline and the RACE-accelerated variant of the GMRES solver preconditioned by algebraic multigrid using Chebyshev smoothers of degree three (same legend as in Fig.~\ref{fig:gmres_amg_gs2}).
		(c) Fraction of AMG runtime spent on the finest level (i.e., level 0) for the different matrices.
		(d) Average speedup of the solver as a function of Chebyshev polynomial degree.
	}
\end{figure}
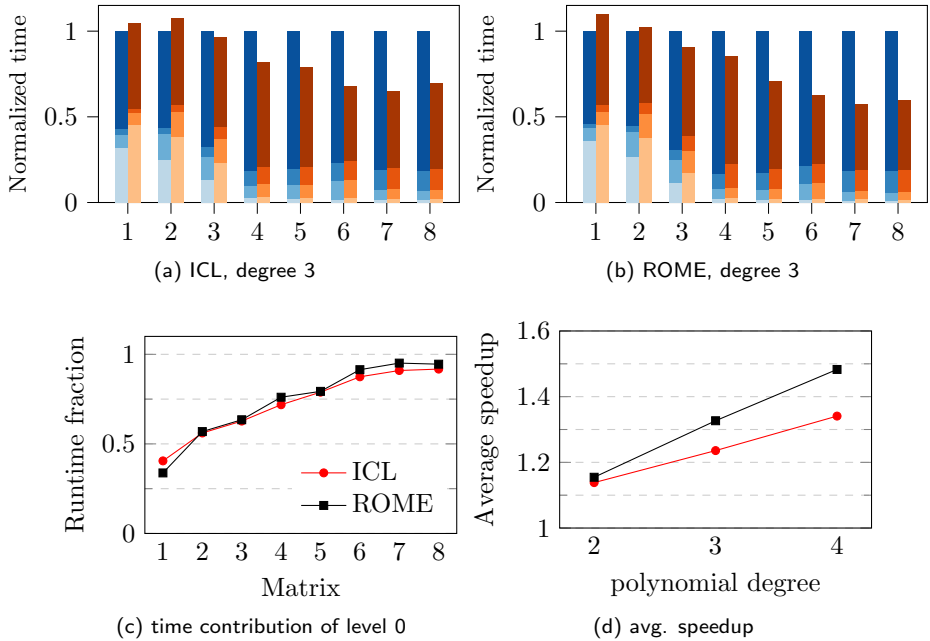

For large problems, most of the solver's runtime is spent in the smoothing operation on the finest grid level.
Typically, a few sweeps of simple iterative methods like Jacobi or GS are used for this.
As the matrix does not change between the sweeps, RACE can cache block the matrix entries.
For example, the GS2 sweeps introduced in Sec.~\ref{sec:GS2_precon} can be used as a smoother
and RACE can block both for inner Jacobi-Richardson iterations within GS2 and outer sweeps of the smoother.
Figures~\ref{fig:gmres_amg_gs2}(a) and (b) demonstrate the speedup attained by RACE over the baseline method when using AMG with the GS2 smoother ($\gamma=1$). 
The baseline performs a single  V-cycle of smoothed-aggregation AMG~\cite{SA_AMG_ref} provided by the MueLu package~\cite{MueLu} in Trilinos as the preconditioner to the GMRES solver
from Belos package.
The pre-smoothing employs two forward sweeps of GS2 while the post-smoothing uses two backward sweeps ($L$ and $U$ are exchanged in Alg.~\ref{alg:GS2}).
The GS2 smoother baseline is provided by the Ifpack2 package.
For the RACE variant we modified the MueLu code such that the cache-blocked variant of GS2 
is used as the smoother for the finest level.
The only difference from the previous implementation shown in Alg.~\ref{alg:GS2}(b) is that we do not need the computation of $Az^1$ in the last sub-power, i.e., at $j=3$ in Alg.~\ref{alg:GS2}(b).
In case of pre-smoothing, at the last sub-power we instead compute the residual that is required in the next step of 
AMG (line 7 in Alg.~\ref{alg:MG}).
This allows us to reuse the matrix $A_k$ from the cache when computing the residual.
Note that this reuse is on top of the reuse in the GS2 sweeps.
As the number of sweeps in the smoother is typically small (in the range of 1--4) we do not tune the power value at which RACE operates
but set it to account for all the sweeps and the residual computation. 
The post-smoother performance can also theoretically benefit by fusing it to the next kernel.
Post-smoothing is the last step of AMG preconditioner; the next step of the 
Krylov solver (in case of right preconditioning) involves an SpMV of the matrix with the preconditioned vector.
By integrating this SpMV as the last sub-power computation, the matrix can be served from the cache. 
However, in the current implementation we do not do this as the SpMV is performed
not by the MueLu package but by the Belos package and therefore would involve fusing kernels
from two different packages, which is possible but requires significant changes.

From Figures~\ref{fig:gmres_amg_gs2}(a) and (b) we see that the cache blocking of the GS2 smoother by RACE achieves a moderate overall speedup of 17\% and 37\% in the solver time compared to the baseline on ICL and ROME, respectively.
For the first time a slowdown of the RACE variant is encountered:  With the \texttt{thermal2} matrix (matrix ID=2) the overall RACE runtime is 19\% and 6\% higher than the baseline on ICL and ROME, respectively.
This is due to the extra cost of RACE's pre-processing (see miscellaneous contributions in Fig.~\ref{fig:gmres_amg_gs2}), which is typically in the range of 30--50 SpMVs (see~\cite{RACE_MPK} for more details).
As the number of solver iterations can be relatively small for AMG-preconditioned solvers (28 in case of the \texttt{thermal2} matrix), this cost cannot always be amortized.

Figures~\ref{fig:gmres_amg_gs2}(c) and (d) show how the number of outer sweeps for the GS2 smoother with one and two inner iterations~($\gamma$) influences the speedup.
Increasing these parameters improves the speedup as higher effective powers in RACE can be used.  
This effect is most pronounced on ROME where we achieve 40\% improvement when applying three outer sweeps instead of a single one.

Good smoothers have the general property that they dampen the error component orthogonal to the coarse grid correction step, which typically means damping high-frequency errors~\cite{cheb_smoother_freq}.
GS-based smoothers enjoy this property.
Another very attractive and commonly used smoother in this regard is the Chebyshev polynomial smoother.
The polynomial is tailored to dampen the high-frequency errors
and is constructed using spectral information and Chebyshev recursion.
Similar to the polynomial preconditioners described in Sec~\ref{sec:gmres_poly_precon}, the Chebyshev polynomial has the property that it can be solely implemented with MPKs as it takes the form shown in~\eqref{eq:poly_precon}.
In contrast to polynomial preconditioners, the degree $d$ of the polynomial smoother is typically low (less than ten).
When using Chebyshev polynomials, each smoothing step of AMG computes the application of a Chebyshev polynomial to a vector using the MPK, which is subject to cache blocking via RACE.
Note that in case of pre-smoothing we also cache-block the residual computation similar to the GS2 pre-smoothing seen above.
As a baseline for comparison  we use the Chebyshev smoother from the Ifpack2 package, which
implements specialized (with appropriate scales and shifts) SpMV-based MPK kernels in Kokkos;  its performance is impacted by the same dynamic scheduling problem for larger matrices as discussed in Sec.~\ref{sec:$s$-step_GMRES_RACE_MPK}. 
Again we modified the code to always use static scheduling in the baseline variant.

Figures~\ref{fig:gmres_amg_cheb}(a) and (b) compare the solver time of baseline and RACE variants using AMG with the Chebyshev smoother of degree three.
On average, the RACE variant achieves a speedup of 24\% and 32\% on ICL and ROME, respectively.
Interestingly, the performance benefit of RACE tends to increase with matrix size\footnote{Matrices are ordered according to increasing size; see Table~\ref{tab:matrices}.}, e.g., on ROME a $1.7\times$ speedup is attained for the largest matrix (\texttt{Flan\_1565}).
This is mainly due to three reasons:
First, for small problems with low iteration counts, the Misc contribution (including RACE preprocessing) to the runtime is significant.
Second, the small test matrices that we have considered (see Table~\ref{tab:matrices}) also tend to have a low \NNZR. 
This makes the smoothing operation 
less prominent compared to the BLAS-1 type (vector-only) operations.
Third, RACE cache blocking is currently only implemented on the finest grid level and the share of this level in overall AMG runtime increases as matrix size increases (see Fig.~\ref{fig:gmres_amg_cheb}(c)).
%
Another observation in line with our previous results is the positive correlation between the polynomial degree $d$ and RACE's speedup (see Fig.~\ref{fig:gmres_amg_cheb}(d)).

In summary, above experiments have demonstrated that RACE provides moderate speed\-ups on the AMG-preconditioned GMRES solver.
Two main factors that currently prevent larger speedups on some matrices are
the large time contribution from coarser grid levels and the low number of solver iterations.

\section{Case study: Momentum equation in the Nalu-Wind solver}
\label{sec:case_study}
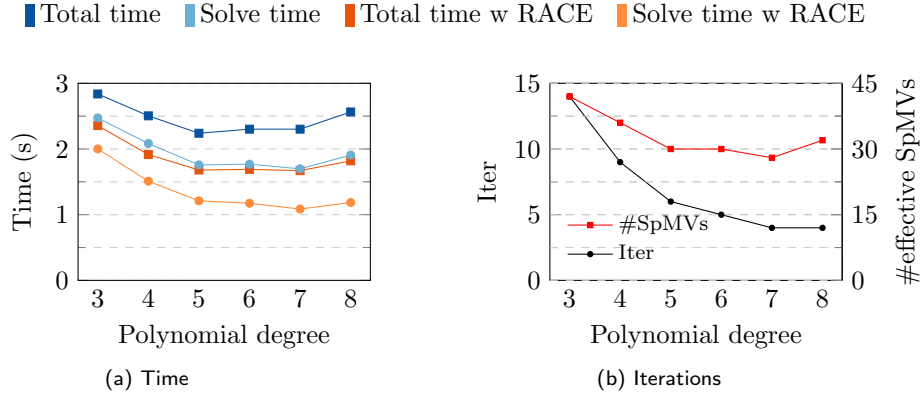
\begin{figure}[tbp]
	\centering
	%
	\input{plots/tg097/NALU_wind/NALU_wind_legend.tex}%

	\hspace{-2em}
	\subfloat[Time]{%
	\input{plots/tg097/NALU_wind/NALU_wind_gmres_poly.tex}%
}
	\hfill
	\subfloat[Iterations]{%
	\input{plots/tg097/NALU_wind/NALU_wind_gmres_poly_iter.tex}%
}
	\caption{\label{fig:NALU_wind_poly}
		(a) Time required to solve the momentum equation using different degrees of the polynomial+Jacobi preconditioner with and without RACE on ROME.
		(b) Number of solver iterations and the effective number of SpMVs performed.
	}
\end{figure}%
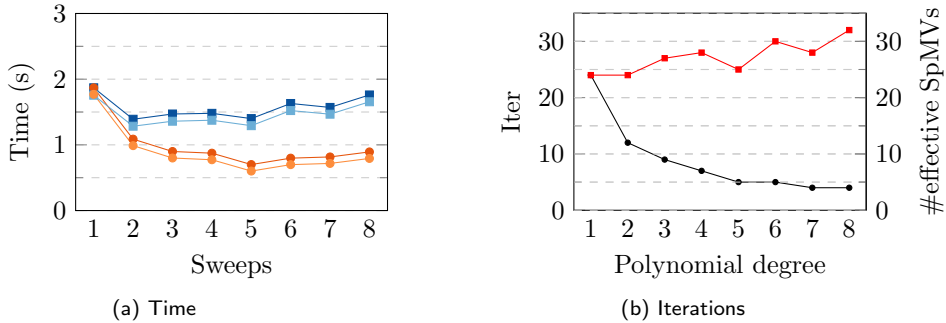
\begin{figure}[tbp]
	\centering
	\hspace{-3em}
	\subfloat[Time]{%
	\input{plots/tg097/NALU_wind/NALU_wind_gmres_Jacobi.tex}%
}
	\hfill
	\subfloat[Iterations]{%
	\input{plots/tg097/NALU_wind/NALU_wind_gmres_Jacobi_iter.tex}%
}
	\caption{\label{fig:NALU_wind_Jacobi}
		(a) Time required to solve the momentum equation using multiple sweeps of the Jacobi preconditioner with and without RACE on ROME.
		(b) Number of solver iterations and effective number of SpMVs performed.
		The same legends as in Fig.~\ref{fig:NALU_wind_poly} apply here.
	}
\end{figure}

To de\-mon\-strate a practical use case of RACE's cache-blocking technique 
we consider the dominant sparse linear system of equations (LSE) in the Nalu-Wind simulation code \cite{Nalu_Wind}.
The LSE arises when solving the unsteady compressible Navier-Stokes equations for the velocities in the simulation of wind turbines.
Specifically, we focus on the case of a large-eddy simulation of two aligned wind turbines under uniform flow, where the turbine blades are modeled using the actuator line model (Example 1.3.4 in~\cite{Nalu-Doc}).
We use a mesh with $256\times 256$ horizontal cells and 64 layers, which translates to a momentum matrix with 12 million rows and almost 300 million nonzeros.
In this scenario the numerical behavior of the linear systems is very similar between time steps after a short start-up phase.
The purpose of this case study is not to claim or find an optimal solver but to demonstrate that 
cache blocking techniques should be taken
into account when selecting and tuning linear solvers for best time to solution.

This application employs an established approach in computational fluid dynamics (CFD), where the model state is 
propagated forward in time using an ODE (ordinary differential equation) time-stepping method, and discretized PDEs (partial differential equations) 
in space are solved in each time step for momentum and conservation, which leads to the sparse linear systems.
Nalu-Wind uses the Trilinos library for solving these time-consuming sparse linear systems.
While the time-integration scheme is implicit and thus unconditionally stable, a small time step is used to cover the relevant physical scales as wind turbines rotate at high velocities.
Consequently, the LSE arising from the PDEs is diagonally dominant and GMRES converges quickly: For the matrix studied here, the default method (GMRES with GS preconditioning) achieves a residual norm of $10^{-12}$ in 13 iterations. 
 This results in a total runtime of 1.83\,s on ROME. 
 We will use this a baseline for comparison when investigating cache blocking in combination with different preconditioning strategies.

Due to the small number of iterations, $s$-step GMRES is unsuitable as it needs a few steps of standard GMRES to calculate the Newton shifts (see Sec.~\ref{sec:$s$-step_GMRES}).
However, using a polynomial preconditioner instead of GS may enable
accelerating the polynomial application using RACE (see Sec.~\ref{sec:gmres_poly_precon}).
Figure~\ref{fig:NALU_wind_poly}(a) compares the time required to solve
the LSE with and without RACE for different polynomial degrees.
The polynomial preconditioner was combined with Jacobi in this case as this proved to be most effective.
Clearly the time to solution can be reduced by increasing the polynomial degree (blue lines in Fig.~\ref{fig:NALU_wind_poly}(a)). This is correlated with a decrease in the number of iterations (see Fig.~\ref{fig:NALU_wind_poly}(b)), which has two positive effects:
First, the total number of SpMVs, the product of iterations and polynomial degree, goes down up to a certain point (red line in Fig.~\ref{fig:NALU_wind_poly}(b)). 
Second, the number of orthogonalization steps decreases.
The combination of theses effects leads to the decrease in solver time as the polynomial degree  increases up to seven.
However, the total time, which includes solver and setup time, is lowest at degree five because the setup cost, which primarily includes $d$ SpMVs where $d$ is the polynomial degree, increases linearly with the degree.
Nevertheless, the default Trilinos implementation of the polynomial preconditioner (without RACE) achieves a minimum runtime (at $d=5$) of  2.24\,s, which is approximately 20\% slower than the above specified baseline (GMRES and GS preconditioner).
This picture changes if the polynomial preconditioner uses RACE's cache blocking (orange lines in Fig.~\ref{fig:NALU_wind_poly}(a)). Time to solution reduces to 1.68\,s, which is a 9\% improvement over the baseline.
Again, the main reason for the limited speedup is the relatively high setup cost of the polynomial preconditioner, which incurs a huge overhead for the small number of iterations at hand. We may conclude that any sophisticated preconditioner with high setup cost would not be effective in this context because the setup has to be repeated in every time step due to the updated momentum matrix.
On the other hand, the RACE setup cost of 30--50 SpMVs for creating the graph traversal scheme may be neglected
in applications where the matrix sparsity pattern stays constant for all time steps, as then the RACE setup is done only once for the entire simulation.

Given these observations, a viable strategy may be increasing the sweeps of basic relaxation preconditioners, which have negligible setup cost.
Figure~\ref{fig:NALU_wind_Jacobi}(a) shows the benefit of increasing the number of Jacobi sweeps in the standard GMRES preconditioner.
In this case the number of iterations also decreases linearly while the SpMV counts remain almost constant up to a certain sweep count (see Fig.~\ref{fig:NALU_wind_Jacobi}(b)). 
Remember that the cost of every Jacobi sweep is similar to an SpMV (see \eqref{eq:Jacobi}).
Although the number of SpMVs remains the same, we reduce the orthogonalization cost in GMRES,
which effectively reduces the solver time (blue lines).
In this case, Jacobi with five sweeps converges in 1.4\,s, achieving a speedup of 1.3$\times$ over the baseline.
Multiple sweeps of Jacobi imply that the same SpMV-like operator is applied back-to-back, which can be accelerated via RACE.
This  further reduces the time to solution to 0.7\,s, which is an additional improvement of a factor of two.
Table~\ref{tab:NALU_solver_comparison} summarizes the results, showing that a total speedup of $2.6\times$ over the default solver is possible.

This study does not make any claims on the optimality of the chosen preconditioners but demonstrates the runtime impact of RACE's cache-blocking technique. 
The above findings indicate that RACE may be particularly useful in situations where the sparsity pattern of matrix is constant (over a long time), 
but the values of matrix elements change too frequently to afford the cost of computing a strong preconditioner repeatedly. 


\begin{table}[!tb]
	\centering
\resizebox{\linewidth}{!}{%
	\begin{tabular}{l | c c c c c c c}
		Preconditioner & None & Jacobi  &  GS & \multicolumn{1}{p{1.2cm}}{\centering Poly, \\ $d$=5} & \multicolumn{1}{p{1.3cm}}{\centering RACE+\\ Poly, \\
			$d$=5} & \multicolumn{1}{p{1.3cm}}{\centering Jacobi, \\ $d$=5} & \multicolumn{1}{p{1.3cm}}{\centering RACE+ \\ Jacobi, \\ $d$=5}\\
		\midrule
		Iter. &  43 &  24 & 13 & 6 & 6 & 5 & 5\\
		\#eff. SpMVs  &  43 &  24 &   26 & 30 & 30 & 25 & 25\\
		Solve time (s) & 3.72 & 1.73 &  1.73 & 1.76 & 1.16 & 1.29 & 0.60\\
		Time (s) & 3.83 & 1.83 &  1.83 & 2.24 & 1.64 & 1.40 & 0.70\\
	\end{tabular}
}
	\caption{Overview of the effectiveness of different preconditioners 
		on the momentum equation in Nalu-Wind. 
		The rows show the number of solver iterations, the effective number of SpMV-like operations, pure solve time, and the total time including the setup cost for the preconditioner. 
		\label{tab:NALU_solver_comparison}
	}
\end{table}

\section{Conclusion and outlook}
\label{sec:conclusion}
In this article, we demonstrated that the node-level performance of various sparse iterative 
solvers can be boosted by performing temporal cache blocking using the Recursive Algebraic Coloring Engine (RACE).
The key is to identify steps in the solver and/or preconditioner which can be (re)formulated into matrix polynomials and then replace the related back-to-back sparse matrix-vector operations with RACE's cache-blocked matrix powers kernel.

First we investigated $s$-step GMRES as a method representative of the broad class of $s$-step Krylov methods. Their basic structure allows to easily benefit from cache-blocked MPKs.
The raw performance improvement of the MPK can be utilized in parts of these solvers, leading to speedups up to 1.5$\times$ for the full solver.
For short-recurrence Krylov methods like conjugate gradients (CG), where the orthogonalization cost is low, the overall improvement may be substantially higher. 
Second, we showed how to apply cache blocking when combining the $s$-step method with preconditioners. Here we addressed the additional problem of calculating polynomials on parts of the sparse matrix, e.g., the triangular factors of the matrix. 
Using a two-stage Gauss-Seidel preconditioner, we further illustrated that cache blocking can be performed across multiple chained operators.  
Third, the challenges and benefits of using RACE in the context of polynomial and algebraic multigrid preconditioners were evaluated. These preconditioners are suitable even for standard Krylov methods, which broadens the scope of our work.
Fourth, we showed that a thorough performance analysis is required to perform fair comparisons with baseline implementations. Performance problems in the baseline libraries have been identified (e.g., a scheduling issue in Trilinos SpMV Kernels) and fixed, if possible. 
Furthermore we showed that the efficiency of cache blocking can be improved by combining it with inter-kernel optimizations (e.g., fusing the pre-smoother with the residual computation of AMG).
Finally, using a case study from wind turbine simulation we illustrated the potential impact of our approach on a real-world application.

RACE's optimizations can be applied to accelerate a variety of other applications like eigenvalue solvers, Chebyshev time propagation, and exponential time integration.
In the future we plan to extend this work to multi-node distributed systems and GPUs.


\appendix
\section{Performance of SpMV routine in Trilinos}
\label{appendix:spmv_perf_trilinos}
\begin{figure}[tbp]
	\centering
	\hspace{-2em}
	\subfloat[SpMV]{%
	\input{plots/GMRES_Transport_spmv.tex}%
}
	\subfloat[Ortho]{%
	\input{plots/GMRES_Transport_ortho.tex}%
}
	\caption{\label{fig:spmv_ortho_perf}
		(a) Performance of the SpMV kernel in \GFS\ with the \texttt{Transport} matrix with default and static scheduling.
		Bars with blue and green colors show the result on ICL and ROME, respectively.
		(b) Performance of orthogonalization routine (in number of routines executed per second) with MKL and BLIS libraries.
		BLIS* represents the BLIS routine called with default setting, while the other one 
		uses an optimized thread parallelization setting.
	}
\end{figure}
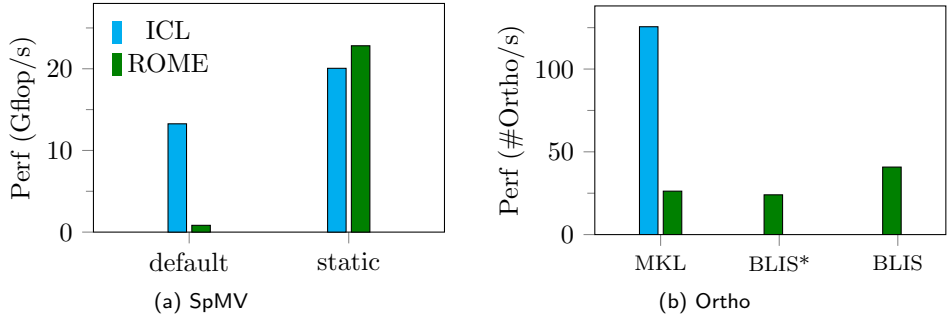

Trilinos uses the Kokkos Kernels package for SpMV, which has been shown to achieve good performance on a wide range of architectures; see~\cite{kokkos_kernels_spmv} and \cite{kokkos_kernels_spmv_2}.
However, initial tests with some of our matrices showed inferior SpMV performance.
For example, the performance on the \texttt{Transport} matrix was well below
the Roof{}line prediction of 24 and 21\,\GFS\ (see~\cite{RACE} for derivation of the performance model) on ICL and ROME, respectively (see Fig.~\ref{fig:spmv_ortho_perf}(a)).
A closer investigation revealed that Trilinos (tested until version 13.4.1) by default calls a dynamically scheduled version of SpMV for matrices with $\NNZ>10^7$.
Especially on ROME, the overhead associated with dynamic scheduling was too high, leading to inferior performance on our benchmark matrices.
Therefore we manually modified the routine to always call a statically scheduled version of SpMV from the Kokkos Kernels.
This led to a huge performance improvement and we ended up close to the memory-bound \rlm\ prediction (see Fig.~\ref{fig:spmv_ortho_perf}(a)).
Note that on ROME the performance slightly exceeds the limit; this is because of a residual caching effect from ROME's large L3 cache~\cite{RACE_MPK}.

\section{Performance of Ortho routines in Trilinos}
\label{appendix:ortho_perf_trilinos}
The Ortho routines in the $s$-step GMRES solver use
tall-skinny \texttt{DGEMM} and \texttt{TRSM} computations, for which
Trilinos employs BLAS libraries.
Figure~\ref{fig:spmv_ortho_perf}(b) reports the performance of the Ortho routine with the Intel MKL~\cite{MKL} and AMD AOCL-BLIS~\cite{AOCL_BLIS,BLIS1} BLAS libraries, respectively.
On ICL, MKL achieves near-optimal performance of 190 orthogonalization steps per second (equivalent to 170\,\GFS) as predicted by the \rlm.
One would expect ROME to match this level due to its practically identical bandwidth and floating-point performance.
However, the MKL version has 4.8$\times$ lower performance on ROME compared to ICL despite
the use of the \texttt{LD\_PRELOAD} trick mentioned in Sec.~\ref{sec:testbed}.
Performance did not improve with the AOCL-BLIS library\footnote{AOCL-BLIS was compiled with gcc v10.2.0 as the library did not support our de-facto Intel compiler.} with default configuration (BLIS* in Fig.~\ref{fig:spmv_ortho_perf}(b)).
However, changing the OpenMP loop used for parallelism via an environment variable (\texttt{BLIS\_JC\_NT}=64) 
yielded a $1.5\times$ speedup compared to MKL.
Thus, we used the AOCL-BLIS library for our runs on ROME with the $s$-step GMRES solver.
The attained performance of Ortho is still far from optimal but it is challenging to do further optimizations from the user level
as the solver requires BLAS computations with various matrix shapes, 
and tuning environment variables globally will not fix the issue in all cases.
We expect that the performance of the AOCL-BLIS library will improve in the future for tall-skinny matrices, leading to a performance boost for both the baseline and RACE-accelerated variants on AMD multi-core processors.

\section*{Acknowledgment}
This work was partially supported by NHR@FAU, which is funded by the State of Bavaria and  by
the Federal Ministry of Education and Research.
The authors would also like to thank NHR@KIT 
for providing access to the HoreKa supercomputer (ICL system), which is funded by the 
Ministry of Science, Research and the Arts Baden-Württemberg and by
the Federal Ministry of Education and Research.

\bibliographystyle{IEEEtran}
\bibliography{IEEEabrv,pub}

\end{document}

%% file: macros.tex
\newcommand{\bq}{\begin{equation}}
\newcommand{\eq}{\end{equation}}

\newcommand{\byte}{\mbox{byte}}
\newcommand{\second}{\mbox{s}}

\newcommand{\flop}{\mbox{flop}}

\newcommand{\iter}{\mbox{it}}

\newcommand{\bits}{\mbox{bits}}

\newcommand{\GBS}{\mbox{G\byte/\second}}

\newcommand{\GFS}{\mbox{G\flop/\second}}

\newcommand{\GHZ}{\mbox{GHz}}

\newcommand{\MB}{\mbox{MB}}

\newcommand{\MiB}{\mbox{MiB}}
\newcommand{\KiB}{\mbox{KiB}}

\newcommand{\eos}{\;.}
\newcommand{\cma}{\;,}
\newcommand{\rlm}{roof{}line model}

\newcommand{\NNZR}{\mbox{$N_\mathrm{nzr}$}}
\newcommand{\NR}{\mbox{$N_\mathrm{r}$}}
\newcommand{\NNZ}{\mbox{$N_\mathrm{nz}$}}

\definecolor{tumbleweed}{rgb}{0.87, 0.67, 0.53}

\newcommand{\spmv}{SpMV}

\newcommand*{\rom}[1]{\expandafter\@slowromancap\romannumeral #1@}

\newenvironment{customlegend}[1][]{%
	\begingroup
	\csname pgfplots@init@cleared@structures\endcsname
	\pgfplotsset{#1}%
}{%
	\csname pgfplots@createlegend\endcsname
	\endgroup
}%
\def\addlegendimage{\csname pgfplots@addlegendimage\endcsname}

\definecolor{applegreen}{rgb}{0.55, 0.71, 0.0}
\definecolor{amethyst}{rgb}{0.6, 0.4, 0.8}
\definecolor{amber}{rgb}{1.0, 0.75, 0.0}

\definecolor{col0}{rgb}{  0.5508,    0.8242,    0.7773 }
\definecolor{col1}{rgb}{  0.98,    0.93,    0.36 }
\definecolor{col2}{rgb}{  0.7422,    0.7266,    0.8516 }
\definecolor{col3}{rgb}{  0.9805,    0.5000,    0.4453 }
\definecolor{col4}{rgb}{  0.5000,    0.6914,    0.8242 }
\definecolor{col5}{rgb}{  0.9883,    0.7031,    0.3828 }
\definecolor{col6}{rgb}{  0.6992,    0.8672,    0.4102 }
\definecolor{col7}{rgb}{  0.9844,    0.8008,    0.8945 }
\definecolor{col8}{rgb}{  0.6484,    0.8047,    0.8867 }
\definecolor{col9}{rgb}{  0.1211,    0.4688,    0.7031 }
\definecolor{col10}{rgb}{  0.6953,    0.8711,    0.5391 }
\definecolor{col11}{rgb}{  0.1992,    0.6250,    0.1719 }
\definecolor{col12}{rgb}{  0.9805,    0.6016,    0.5977 }
\definecolor{col13}{rgb}{  0.8867,    0.1016,    0.1094 }
\definecolor{col14}{rgb}{  0.9883,    0.7461,    0.4336 }
\definecolor{armygreen}{rgb}{0.29, 0.33, 0.13}
\definecolor{applegreen}{rgb}{0.55, 0.71, 0.0}
\definecolor{ao(english)}{rgb}{0.0, 0.5, 0.0}

%% file: matrices.tex
\begin{tabular}{l l S[table-format=8.0, table-space-text-pre=(, table-space-text-post=)] S[table-format=9.0, table-space-text-pre=(, table-space-text-post=)] S[round-mode=places,round-precision=2]}
{ID} & {Matrix name} &  {\NR} & {\NNZ} & {\NNZR} \\
\midrule
{1} & {G3\_circuit} & 1585478 & 7660826 & 4.831871524 \\
{2} & {thermal2} & 1228045 & 8580313 & 6.986969533 \\
{3} & {Transport} & 1602111 & 23487281 & 14.6602 \\
{4} & {Fault\_639} & 638802 & 28614564 & 44.79410522 \\
{5} & {Emilia\_923} & 923136 & 41005206 & 44.41946365 \\
{6} & {af\_shell10} & 1508065 & 52672325 & 34.927092 \\
{7} & {ML\_Geer} & 1504002 & 110879972 & 73.7232876 \\
{8} & {Flan\_1565} & 1564794 & 117406044 & 75.02971254 \\
\end{tabular}

%% file: plots/perf_vs_power/plot_p_Flan_1565_ICL.tex
\begin{tikzpicture}
    \begin{axis}[
                ymin=0, ymax=65,
                xtick={1, 2, 3, 4, 5, 6, 7, 8},
                xticklabels={1, 2, 3, 4, 5, 6, 8, 10},
                width  = 0.48\linewidth,
				height = 0.25\textheight,
                major x tick style = transparent,
                grid = minor,	
                ymajorgrids = true,
                grid style={dashed, gray!40},
                ylabel = {\normalsize{Perf (\GFS)}},
                xlabel = {Power ($p_m$)},
                tick label style={font={\normalsize}},
                scaled y ticks = false,
                enlarge x limits=0.035,
                legend cell align=left,
                legend style={font=\normalsize},
                legend columns=1,
                legend style={
                    legend pos=south west,
                    draw=none
                },
            ]
\addplot[mark=*, mark size=1.5pt, mark options={cyan}, draw=cyan ] plot coordinates{(1,26.1573) (2,26.0637) (3,26.029) (4,25.9954) (5,26.0812) (6,25.903) (7,26.005) (8,26.0761)}; 
\addplot[mark=square*, mark size=1.5pt, mark options={orange}, draw=orange ] plot coordinates{(1,26.224) (2,42.7551) (3,52.316) (4,58.5135) (5,59.7939) (6,57.3075) (7,52.4784) (8,45.7434)};

\legend{Baseline, RACE}

    \end{axis}
\end{tikzpicture}

%% file: plots/perf_vs_power/plot_p_Flan_1565_ROME.tex
\begin{tikzpicture}
    \begin{axis}[
                ymin=0, ymax=100,
                xtick={1, 2, 3, 4, 5, 6, 7},
                xticklabels={1, 2, 3, 4, 6, 8, 10},
                width  = 0.48\linewidth,
                height = 0.25\textheight,
                major x tick style = transparent,
                grid = minor,	
                ymajorgrids = true,
                grid style={dashed, gray!40},
                ylabel = {\normalsize{Perf (\GFS)}},
                xlabel = {Power ($p_m$)},
                tick label style={font={\normalsize}},
                scaled y ticks = false,
                enlarge x limits=0.035,
                legend cell align=left,
                legend style={font=\normalsize},
                legend columns=-1,
                legend style={
                    legend pos=north east
                },
            ]
\addplot[mark=*, mark size=1.5pt, mark options={cyan}, draw=cyan ] plot coordinates{(1,22.5899) (2,22.5018) (3,22.5302) (4,22.5464) (5,22.5135) (6,22.5463) (7,22.5007)}; \addplot[mark=square*, mark size=1.5pt, mark options={orange}, draw=orange ] plot coordinates{(1,22.3019) (2,40.7389) (3,53.7471) (4,64.4385) (5,77.2183) (6,84.6459) (7,90.3267)};


    \end{axis}
\end{tikzpicture}

%% file: plots/gmres_s_step_legend.tex
\begin{tikzpicture}
	\definecolor{burlywood253190133}{RGB}{253,190,133}
	\definecolor{coral25314160}{RGB}{253,141,60}
	\definecolor{cornflowerblue107174214}{RGB}{107,174,214}
	\definecolor{darkgray176}{RGB}{176,176,176}
	\definecolor{lightgray204}{RGB}{204,204,204}
	\definecolor{orangered217711}{RGB}{217,71,1}
	\definecolor{powderblue189215231}{RGB}{189,215,231}
	\definecolor{steelblue33113181}{RGB}{33,113,181}
  \centering
  \begin{customlegend}[
      legend columns=3,
      legend style={
      anchor=north,
      /tikz/every even column/.append style={column sep=0.5cm}},
      legend entries={Baseline: MPK, Baseline: Ortho, Baseline: Misc,  RACE: MPK, RACE: Ortho, RACE: Misc},
      width  = \linewidth,
      legend image post style={xscale=0.5},
      legend image code/.code={\draw[#1, draw=none] (0cm,-0.15cm) rectangle (0.26cm,0.15cm);}, 
	  legend style={draw=none} 
      ]
    \addlegendimage{fill=steelblue33113181, color=steelblue33113181}
    \addlegendimage{fill=cornflowerblue107174214, color=cornflowerblue107174214}
    \addlegendimage{fill=powderblue189215231, color=powderblue189215231}  
    
    \addlegendimage{fill=orangered217711, color=orangered217711}
    \addlegendimage{fill=coral25314160, color=coral25314160}
    \addlegendimage{fill=burlywood253190133, color=burlywood253190133}  
  \end{customlegend}
\end{tikzpicture}

%% file: plots/horeka/gmres_s_step_none.tex
\begin{tikzpicture}

\definecolor{burlywood253190133}{RGB}{253,190,133}
\definecolor{coral25314160}{RGB}{253,141,60}
\definecolor{cornflowerblue107174214}{RGB}{107,174,214}
\definecolor{darkgray176}{RGB}{176,176,176}
\definecolor{lightgray204}{RGB}{204,204,204}
\definecolor{orangered217711}{RGB}{217,71,1}
\definecolor{powderblue189215231}{RGB}{189,215,231}
\definecolor{steelblue33113181}{RGB}{33,113,181}

\begin{axis}[
	legend cell align={left},
	legend style={fill opacity=0.8, draw opacity=1, text opacity=1, draw=lightgray204},
	tick align=outside,
	tick pos=left,
	x grid style={darkgray176},
	xtick={0, 1, 2, 3, 4, 5, 6, 7},
	xticklabels={1, 2, 3, 4, 5, 6, 7, 8},
	width  = 0.49\linewidth,
	height = 0.2\textheight,
	xmin=-0.68, xmax=7.68,
	xtick style={color=black},
	y grid style={darkgray176},
	ylabel={\small Normalized time},
	ymin=0, ymax=1.15,
	ytick style={color=black}
	]
\draw[draw=none,fill=powderblue189215231] (axis cs:0,0) rectangle (axis cs:-0.3,0.196674187937946);
\draw[draw=none,fill=powderblue189215231] (axis cs:1,0) rectangle (axis cs:0.7,0.180751401609995);
\draw[draw=none,fill=powderblue189215231] (axis cs:2,0) rectangle (axis cs:1.7,0.152093534241466);
\draw[draw=none,fill=powderblue189215231] (axis cs:3,0) rectangle (axis cs:2.7,0.134423862076933);
\draw[draw=none,fill=powderblue189215231] (axis cs:4,0) rectangle (axis cs:3.7,0.0807818153497478);
\draw[draw=none,fill=powderblue189215231] (axis cs:5,0) rectangle (axis cs:4.7,0.122622997637991);
\draw[draw=none,fill=powderblue189215231] (axis cs:6,0) rectangle (axis cs:5.7,0.0940149553772603);
\draw[draw=none,fill=powderblue189215231] (axis cs:7,0) rectangle (axis cs:6.7,0.0934263342092186);
\draw[draw=none,fill=cornflowerblue107174214] (axis cs:0,0.196674187937946) rectangle (axis cs:-0.3,0.70447192687524);
\draw[draw=none,fill=cornflowerblue107174214] (axis cs:1,0.180751401609995) rectangle (axis cs:0.7,0.684359417658018);
\draw[draw=none,fill=cornflowerblue107174214] (axis cs:2,0.152093534241466) rectangle (axis cs:1.7,0.558777278509191);
\draw[draw=none,fill=cornflowerblue107174214] (axis cs:3,0.134423862076933) rectangle (axis cs:2.7,0.363522823438357);
\draw[draw=none,fill=cornflowerblue107174214] (axis cs:4,0.0807818153497478) rectangle (axis cs:3.7,0.564906833430566);
\draw[draw=none,fill=cornflowerblue107174214] (axis cs:5,0.122622997637991) rectangle (axis cs:4.7,0.396250195702676);
\draw[draw=none,fill=cornflowerblue107174214] (axis cs:6,0.0940149553772603) rectangle (axis cs:5.7,0.264582114045204);
\draw[draw=none,fill=cornflowerblue107174214] (axis cs:7,0.0934263342092186) rectangle (axis cs:6.7,0.254066169062031);
\draw[draw=none,fill=steelblue33113181] (axis cs:0,0.70447192687524) rectangle (axis cs:-0.3,1);
\draw[draw=none,fill=steelblue33113181] (axis cs:1,0.684359417658018) rectangle (axis cs:0.7,1);
\draw[draw=none,fill=steelblue33113181] (axis cs:2,0.558777278509191) rectangle (axis cs:1.7,1);
\draw[draw=none,fill=steelblue33113181] (axis cs:3,0.363522823438357) rectangle (axis cs:2.7,1);
\draw[draw=none,fill=steelblue33113181] (axis cs:4,0.564906833430566) rectangle (axis cs:3.7,1);
\draw[draw=none,fill=steelblue33113181] (axis cs:5,0.396250195702676) rectangle (axis cs:4.7,1);
\draw[draw=none,fill=steelblue33113181] (axis cs:6,0.264582114045204) rectangle (axis cs:5.7,1);
\draw[draw=none,fill=steelblue33113181] (axis cs:7,0.254066169062031) rectangle (axis cs:6.7,1);
\draw[draw=none,fill=burlywood253190133] (axis cs:0,0) rectangle (axis cs:0.3,0.168673158282794);
\draw[draw=none,fill=burlywood253190133] (axis cs:1,0) rectangle (axis cs:1.3,0.208414684895012);
\draw[draw=none,fill=burlywood253190133] (axis cs:2,0) rectangle (axis cs:2.3,0.15938476513917);
\draw[draw=none,fill=burlywood253190133] (axis cs:3,0) rectangle (axis cs:3.3,0.145863208658906);
\draw[draw=none,fill=burlywood253190133] (axis cs:4,0) rectangle (axis cs:4.3,0.0928391902133725);
\draw[draw=none,fill=burlywood253190133] (axis cs:5,0) rectangle (axis cs:5.3,0.139231755751938);
\draw[draw=none,fill=burlywood253190133] (axis cs:6,0) rectangle (axis cs:6.3,0.120677630664924);
\draw[draw=none,fill=burlywood253190133] (axis cs:7,0) rectangle (axis cs:7.3,0.131039573417848);
\draw[draw=none,fill=coral25314160] (axis cs:0,0.168673158282794) rectangle (axis cs:0.3,0.676102714404733);
\draw[draw=none,fill=coral25314160] (axis cs:1,0.208414684895012) rectangle (axis cs:1.3,0.711449905234869);
\draw[draw=none,fill=coral25314160] (axis cs:2,0.15938476513917) rectangle (axis cs:2.3,0.574126300242631);
\draw[draw=none,fill=coral25314160] (axis cs:3,0.145863208658906) rectangle (axis cs:3.3,0.385492626100099);
\draw[draw=none,fill=coral25314160] (axis cs:4,0.0928391902133725) rectangle (axis cs:4.3,0.572964241460367);
\draw[draw=none,fill=coral25314160] (axis cs:5,0.139231755751938) rectangle (axis cs:5.3,0.417884449557859);
\draw[draw=none,fill=coral25314160] (axis cs:6,0.120677630664924) rectangle (axis cs:6.3,0.285463701327861);
\draw[draw=none,fill=coral25314160] (axis cs:7,0.131039573417848) rectangle (axis cs:7.3,0.290187228110322);
\draw[draw=none,fill=orangered217711] (axis cs:0,0.676102714404733) rectangle (axis cs:0.3,0.887422468553653);
\draw[draw=none,fill=orangered217711] (axis cs:1,0.711449905234869) rectangle (axis cs:1.3,0.923501891112163);
\draw[draw=none,fill=orangered217711] (axis cs:2,0.574126300242631) rectangle (axis cs:2.3,0.835079421077534);
\draw[draw=none,fill=orangered217711] (axis cs:3,0.385492626100099) rectangle (axis cs:3.3,0.769657339807071);
\draw[draw=none,fill=orangered217711] (axis cs:4,0.572964241460367) rectangle (axis cs:4.3,0.831182745573777);
\draw[draw=none,fill=orangered217711] (axis cs:5,0.417884449557859) rectangle (axis cs:5.3,0.711247763452478);
\draw[draw=none,fill=orangered217711] (axis cs:6,0.285463701327861) rectangle (axis cs:6.3,0.637483506626932);
\draw[draw=none,fill=orangered217711] (axis cs:7,0.290187228110322) rectangle (axis cs:7.3,0.689083703538283);
\draw (axis cs:0.15,0.937422468553653) node[
scale=0.72,
anchor=base,
text=black,
rotate=0.0
]{2};
\draw (axis cs:1.15,0.973501891112163) node[
scale=0.72,
anchor=base,
text=black,
rotate=0.0
]{4};
\draw (axis cs:2.15,0.885079421077534) node[
scale=0.72,
anchor=base,
text=black,
rotate=0.0
]{4};
\draw (axis cs:3.15,0.819657339807071) node[
scale=0.72,
anchor=base,
text=black,
rotate=0.0
]{4};
\draw (axis cs:4.15,0.881182745573777) node[
scale=0.72,
anchor=base,
text=black,
rotate=0.0
]{4};
\draw (axis cs:5.15,0.761247763452478) node[
scale=0.72,
anchor=base,
text=black,
rotate=0.0
]{4};
\draw (axis cs:6.15,0.687483506626932) node[
scale=0.72,
anchor=base,
text=black,
rotate=0.0
]{4};
\draw (axis cs:7.15,0.739083703538283) node[
scale=0.72,
anchor=base,
text=black,
rotate=0.0
]{4};
\end{axis}

\end{tikzpicture}

%% file: plots/tg097/gmres_s_step_none.tex
\begin{tikzpicture}

\definecolor{burlywood253190133}{RGB}{253,190,133}
\definecolor{coral25314160}{RGB}{253,141,60}
\definecolor{cornflowerblue107174214}{RGB}{107,174,214}
\definecolor{darkgray176}{RGB}{176,176,176}
\definecolor{lightgray204}{RGB}{204,204,204}
\definecolor{orangered217711}{RGB}{217,71,1}
\definecolor{powderblue189215231}{RGB}{189,215,231}
\definecolor{steelblue33113181}{RGB}{33,113,181}

	\begin{axis}[
		legend cell align={left},
		legend style={fill opacity=0.8, draw opacity=1, text opacity=1, draw=lightgray204},
		tick align=outside,
		tick pos=left,
		x grid style={darkgray176},
		xtick={0, 1, 2, 3, 4, 5, 6, 7},
		xticklabels={1, 2, 3, 4, 5, 6, 7, 8},
		width  = 0.49\linewidth,
		height = 0.2\textheight,
		xmin=-0.68, xmax=7.68,
		xtick style={color=black},
		y grid style={darkgray176},
		ylabel={\small Normalized time},
		ymin=0, ymax=1.15,
		ytick style={color=black},
		]
\draw[draw=none,fill=powderblue189215231] (axis cs:0,0) rectangle (axis cs:-0.3,0.231315499080866);
\draw[draw=none,fill=powderblue189215231] (axis cs:1,0) rectangle (axis cs:0.7,0.227247200783896);
\draw[draw=none,fill=powderblue189215231] (axis cs:2,0) rectangle (axis cs:1.7,0.210069959450192);
\draw[draw=none,fill=powderblue189215231] (axis cs:3,0) rectangle (axis cs:2.7,0.167245196929856);
\draw[draw=none,fill=powderblue189215231] (axis cs:4,0) rectangle (axis cs:3.7,0.167563668855043);
\draw[draw=none,fill=powderblue189215231] (axis cs:5,0) rectangle (axis cs:4.7,0.172692482604868);
\draw[draw=none,fill=powderblue189215231] (axis cs:6,0) rectangle (axis cs:5.7,0.141108604249127);
\draw[draw=none,fill=powderblue189215231] (axis cs:7,0) rectangle (axis cs:6.7,0.139610203812645);
\draw[draw=none,fill=cornflowerblue107174214] (axis cs:0,0.231315499080866) rectangle (axis cs:-0.3,0.917504269545065);
\draw[draw=none,fill=cornflowerblue107174214] (axis cs:1,0.227247200783896) rectangle (axis cs:0.7,0.885204113936916);
\draw[draw=none,fill=cornflowerblue107174214] (axis cs:2,0.210069959450192) rectangle (axis cs:1.7,0.802961714454209);
\draw[draw=none,fill=cornflowerblue107174214] (axis cs:3,0.167245196929856) rectangle (axis cs:2.7,0.608064492684909);
\draw[draw=none,fill=cornflowerblue107174214] (axis cs:4,0.167563668855043) rectangle (axis cs:3.7,0.572458504284165);
\draw[draw=none,fill=cornflowerblue107174214] (axis cs:5,0.172692482604868) rectangle (axis cs:4.7,0.614945588470496);
\draw[draw=none,fill=cornflowerblue107174214] (axis cs:6,0.141108604249127) rectangle (axis cs:5.7,0.457285679366861);
\draw[draw=none,fill=cornflowerblue107174214] (axis cs:7,0.139610203812645) rectangle (axis cs:6.7,0.447785130504407);
\draw[draw=none,fill=steelblue33113181] (axis cs:0,0.917504269545065) rectangle (axis cs:-0.3,1);
\draw[draw=none,fill=steelblue33113181] (axis cs:1,0.885204113936916) rectangle (axis cs:0.7,1);
\draw[draw=none,fill=steelblue33113181] (axis cs:2,0.802961714454209) rectangle (axis cs:1.7,1);
\draw[draw=none,fill=steelblue33113181] (axis cs:3,0.608064492684909) rectangle (axis cs:2.7,1);
\draw[draw=none,fill=steelblue33113181] (axis cs:4,0.572458504284165) rectangle (axis cs:3.7,1);
\draw[draw=none,fill=steelblue33113181] (axis cs:5,0.614945588470496) rectangle (axis cs:4.7,1);
\draw[draw=none,fill=steelblue33113181] (axis cs:6,0.457285679366861) rectangle (axis cs:5.7,1);
\draw[draw=none,fill=steelblue33113181] (axis cs:7,0.447785130504407) rectangle (axis cs:6.7,1);
\draw[draw=none,fill=burlywood253190133] (axis cs:0,0) rectangle (axis cs:0.3,0.23558948118982);
\draw[draw=none,fill=burlywood253190133] (axis cs:1,0) rectangle (axis cs:1.3,0.22724306855855);
\draw[draw=none,fill=burlywood253190133] (axis cs:2,0) rectangle (axis cs:2.3,0.215088710170206);
\draw[draw=none,fill=burlywood253190133] (axis cs:3,0) rectangle (axis cs:3.3,0.189520630862161);
\draw[draw=none,fill=burlywood253190133] (axis cs:4,0) rectangle (axis cs:4.3,0.184711193889611);
\draw[draw=none,fill=burlywood253190133] (axis cs:5,0) rectangle (axis cs:5.3,0.183764381441895);
\draw[draw=none,fill=burlywood253190133] (axis cs:6,0) rectangle (axis cs:6.3,0.15649917573462);
\draw[draw=none,fill=burlywood253190133] (axis cs:7,0) rectangle (axis cs:7.3,0.155506819062499);
\draw[draw=none,fill=coral25314160] (axis cs:0,0.23558948118982) rectangle (axis cs:0.3,0.925549738591585);
\draw[draw=none,fill=coral25314160] (axis cs:1,0.22724306855855) rectangle (axis cs:1.3,0.887732635955879);
\draw[draw=none,fill=coral25314160] (axis cs:2,0.215088710170206) rectangle (axis cs:2.3,0.811534261907514);
\draw[draw=none,fill=coral25314160] (axis cs:3,0.189520630862161) rectangle (axis cs:3.3,0.646346122832069);
\draw[draw=none,fill=coral25314160] (axis cs:4,0.184711193889611) rectangle (axis cs:4.3,0.594066735132782);
\draw[draw=none,fill=coral25314160] (axis cs:5,0.183764381441895) rectangle (axis cs:5.3,0.625860992931775);
\draw[draw=none,fill=coral25314160] (axis cs:6,0.15649917573462) rectangle (axis cs:6.3,0.474227496358451);
\draw[draw=none,fill=coral25314160] (axis cs:7,0.155506819062499) rectangle (axis cs:7.3,0.464314656779471);
\draw[draw=none,fill=orangered217711] (axis cs:0,0.925549738591585) rectangle (axis cs:0.3,0.982595697043997);
\draw[draw=none,fill=orangered217711] (axis cs:1,0.887732635955879) rectangle (axis cs:1.3,0.949786663869195);
\draw[draw=none,fill=orangered217711] (axis cs:2,0.811534261907514) rectangle (axis cs:2.3,0.922389028007966);
\draw[draw=none,fill=orangered217711] (axis cs:3,0.646346122832069) rectangle (axis cs:3.3,0.848728467772692);
\draw[draw=none,fill=orangered217711] (axis cs:4,0.594066735132782) rectangle (axis cs:4.3,0.772572744103833);
\draw[draw=none,fill=orangered217711] (axis cs:5,0.625860992931775) rectangle (axis cs:5.3,0.773153499386386);
\draw[draw=none,fill=orangered217711] (axis cs:6,0.474227496358451) rectangle (axis cs:6.3,0.678299382848056);
\draw[draw=none,fill=orangered217711] (axis cs:7,0.464314656779471) rectangle (axis cs:7.3,0.686550814764823);
\draw (axis cs:0.15,1.032595697044) node[
scale=0.72,
anchor=base,
text=black,
rotate=0.0
]{3};
\draw (axis cs:1.15,0.999786663869195) node[
scale=0.72,
anchor=base,
text=black,
rotate=0.0
]{3};
\draw (axis cs:2.15,0.972389028007966) node[
scale=0.72,
anchor=base,
text=black,
rotate=0.0
]{4};
\draw (axis cs:3.15,0.898728467772692) node[
scale=0.72,
anchor=base,
text=black,
rotate=0.0
]{4};
\draw (axis cs:4.15,0.822572744103833) node[
scale=0.72,
anchor=base,
text=black,
rotate=0.0
]{4};
\draw (axis cs:5.15,0.823153499386386) node[
scale=0.72,
anchor=base,
text=black,
rotate=0.0
]{4};
\draw (axis cs:6.15,0.728299382848056) node[
scale=0.72,
anchor=base,
text=black,
rotate=0.0
]{4};
\draw (axis cs:7.15,0.736550814764823) node[
scale=0.72,
anchor=base,
text=black,
rotate=0.0
]{4};
\end{axis}

\end{tikzpicture}

%% file: plots/gmres_s_step_legend_w_precon.tex
\begin{tikzpicture}
	\definecolor{burlywood253190133}{RGB}{253,190,133}
	\definecolor{coral25314160}{RGB}{253,141,60}
	\definecolor{cornflowerblue107174214}{RGB}{107,174,214}
	\definecolor{darkgray176}{RGB}{176,176,176}
	\definecolor{lightgray204}{RGB}{204,204,204}
	\definecolor{orangered217711}{RGB}{217,71,1}
	\definecolor{powderblue189215231}{RGB}{189,215,231}
	\definecolor{steelblue33113181}{RGB}{33,113,181}
  \centering
  \begin{customlegend}[
      legend columns=3,
      legend style={
      anchor=north,
      /tikz/every even column/.append style={column sep=0.5cm}},
      legend entries={Baseline: MPK+Precon, Baseline: Ortho, Baseline: Misc, RACE: MPK+Precon, RACE: Ortho, RACE: Misc},
      width  = \linewidth,
      legend image post style={xscale=0.5},
      legend image code/.code={\draw[#1, draw=none] (0cm,-0.15cm) rectangle (0.26cm,0.15cm);}, 
	  legend style={draw=none} 
      ]
    \addlegendimage{fill=steelblue33113181, color=steelblue33113181}
    \addlegendimage{fill=cornflowerblue107174214, color=cornflowerblue107174214}
    \addlegendimage{fill=powderblue189215231, color=powderblue189215231}   
    
    \addlegendimage{fill=orangered217711, color=orangered217711}
    \addlegendimage{fill=coral25314160, color=coral25314160}
    \addlegendimage{fill=burlywood253190133, color=burlywood253190133}  
  \end{customlegend}
\end{tikzpicture}

%% file: plots/horeka/gmres_s_step_Jacobi.tex
\begin{tikzpicture}

\definecolor{burlywood253190133}{RGB}{253,190,133}
\definecolor{coral25314160}{RGB}{253,141,60}
\definecolor{cornflowerblue107174214}{RGB}{107,174,214}
\definecolor{darkgray176}{RGB}{176,176,176}
\definecolor{lightgray204}{RGB}{204,204,204}
\definecolor{orangered217711}{RGB}{217,71,1}
\definecolor{powderblue189215231}{RGB}{189,215,231}
\definecolor{steelblue33113181}{RGB}{33,113,181}

\begin{axis}[
	legend cell align={left},
legend style={fill opacity=0.8, draw opacity=1, text opacity=1, draw=lightgray204},
tick align=outside,
tick pos=left,
x grid style={darkgray176},
xtick={0, 1, 2, 3, 4, 5, 6, 7},
xticklabels={1, 2, 3, 4, 5, 6, 7, 8},
width  = 0.49\linewidth,
height = 0.2\textheight,
xmin=-0.68, xmax=7.68,
xtick style={color=black},
y grid style={darkgray176},
ylabel={\small Normalized time},
ymin=0, ymax=1.15,
ytick style={color=black}
]
\draw[draw=none,fill=powderblue189215231] (axis cs:0,0) rectangle (axis cs:-0.3,0.167680203029951);
\draw[draw=none,fill=powderblue189215231] (axis cs:1,0) rectangle (axis cs:0.7,0.179686309360189);
\draw[draw=none,fill=powderblue189215231] (axis cs:2,0) rectangle (axis cs:1.7,0.149188967145305);
\draw[draw=none,fill=powderblue189215231] (axis cs:3,0) rectangle (axis cs:2.7,0.132474077154232);
\draw[draw=none,fill=powderblue189215231] (axis cs:4,0) rectangle (axis cs:3.7,0.0832110183885939);
\draw[draw=none,fill=powderblue189215231] (axis cs:5,0) rectangle (axis cs:4.7,0.120578073079956);
\draw[draw=none,fill=powderblue189215231] (axis cs:6,0) rectangle (axis cs:5.7,0.0996581981807221);
\draw[draw=none,fill=powderblue189215231] (axis cs:7,0) rectangle (axis cs:6.7,0.0934691124341076);
\draw[draw=none,fill=cornflowerblue107174214] (axis cs:0,0.167680203029951) rectangle (axis cs:-0.3,0.6690796819598);
\draw[draw=none,fill=cornflowerblue107174214] (axis cs:1,0.179686309360189) rectangle (axis cs:0.7,0.65842400780074);
\draw[draw=none,fill=cornflowerblue107174214] (axis cs:2,0.149188967145305) rectangle (axis cs:1.7,0.545466531918552);
\draw[draw=none,fill=cornflowerblue107174214] (axis cs:3,0.132474077154232) rectangle (axis cs:2.7,0.360529407233889);
\draw[draw=none,fill=cornflowerblue107174214] (axis cs:4,0.0832110183885939) rectangle (axis cs:3.7,0.558412336185681);
\draw[draw=none,fill=cornflowerblue107174214] (axis cs:5,0.120578073079956) rectangle (axis cs:4.7,0.389758117969075);
\draw[draw=none,fill=cornflowerblue107174214] (axis cs:6,0.0996581981807221) rectangle (axis cs:5.7,0.261319446598451);
\draw[draw=none,fill=cornflowerblue107174214] (axis cs:7,0.0934691124341076) rectangle (axis cs:6.7,0.252547097955951);
\draw[draw=none,fill=steelblue33113181] (axis cs:0,0.6690796819598) rectangle (axis cs:-0.3,1);
\draw[draw=none,fill=steelblue33113181] (axis cs:1,0.65842400780074) rectangle (axis cs:0.7,1);
\draw[draw=none,fill=steelblue33113181] (axis cs:2,0.545466531918552) rectangle (axis cs:1.7,1);
\draw[draw=none,fill=steelblue33113181] (axis cs:3,0.360529407233889) rectangle (axis cs:2.7,1);
\draw[draw=none,fill=steelblue33113181] (axis cs:4,0.558412336185681) rectangle (axis cs:3.7,1);
\draw[draw=none,fill=steelblue33113181] (axis cs:5,0.389758117969075) rectangle (axis cs:4.7,1);
\draw[draw=none,fill=steelblue33113181] (axis cs:6,0.261319446598451) rectangle (axis cs:5.7,1);
\draw[draw=none,fill=steelblue33113181] (axis cs:7,0.252547097955951) rectangle (axis cs:6.7,1);
\draw[draw=none,fill=burlywood253190133] (axis cs:0,0) rectangle (axis cs:0.3,0.177196039828701);
\draw[draw=none,fill=burlywood253190133] (axis cs:1,0) rectangle (axis cs:1.3,0.187597785390726);
\draw[draw=none,fill=burlywood253190133] (axis cs:2,0) rectangle (axis cs:2.3,0.153329277519089);
\draw[draw=none,fill=burlywood253190133] (axis cs:3,0) rectangle (axis cs:3.3,0.147155744985233);
\draw[draw=none,fill=burlywood253190133] (axis cs:4,0) rectangle (axis cs:4.3,0.0924646521525088);
\draw[draw=none,fill=burlywood253190133] (axis cs:5,0) rectangle (axis cs:5.3,0.147512477977784);
\draw[draw=none,fill=burlywood253190133] (axis cs:6,0) rectangle (axis cs:6.3,0.12146363454793);
\draw[draw=none,fill=burlywood253190133] (axis cs:7,0) rectangle (axis cs:7.3,0.12649675408428);
\draw[draw=none,fill=coral25314160] (axis cs:0,0.177196039828701) rectangle (axis cs:0.3,0.668629044737308);
\draw[draw=none,fill=coral25314160] (axis cs:1,0.187597785390726) rectangle (axis cs:1.3,0.655709296599536);
\draw[draw=none,fill=coral25314160] (axis cs:2,0.153329277519089) rectangle (axis cs:2.3,0.548548319344117);
\draw[draw=none,fill=coral25314160] (axis cs:3,0.147155744985233) rectangle (axis cs:3.3,0.373056560845905);
\draw[draw=none,fill=coral25314160] (axis cs:4,0.0924646521525088) rectangle (axis cs:4.3,0.56829923396585);
\draw[draw=none,fill=coral25314160] (axis cs:5,0.147512477977784) rectangle (axis cs:5.3,0.413879119386043);
\draw[draw=none,fill=coral25314160] (axis cs:6,0.12146363454793) rectangle (axis cs:6.3,0.283650050210912);
\draw[draw=none,fill=coral25314160] (axis cs:7,0.12649675408428) rectangle (axis cs:7.3,0.285933332667012);
\draw[draw=none,fill=orangered217711] (axis cs:0,0.668629044737308) rectangle (axis cs:0.3,0.893585715610717);
\draw[draw=none,fill=orangered217711] (axis cs:1,0.655709296599536) rectangle (axis cs:1.3,0.917737651181982);
\draw[draw=none,fill=orangered217711] (axis cs:2,0.548548319344117) rectangle (axis cs:2.3,0.837980806193362);
\draw[draw=none,fill=orangered217711] (axis cs:3,0.373056560845905) rectangle (axis cs:3.3,0.806433594909314);
\draw[draw=none,fill=orangered217711] (axis cs:4,0.56829923396585) rectangle (axis cs:4.3,0.863594930898864);
\draw[draw=none,fill=orangered217711] (axis cs:5,0.413879119386043) rectangle (axis cs:5.3,0.747919448124906);
\draw[draw=none,fill=orangered217711] (axis cs:6,0.283650050210912) rectangle (axis cs:6.3,0.686325317804138);
\draw[draw=none,fill=orangered217711] (axis cs:7,0.285933332667012) rectangle (axis cs:7.3,0.728544747721771);
\draw (axis cs:0.15,0.943585715610717) node[
scale=0.72,
anchor=base,
text=black,
rotate=0.0
]{3};
\draw (axis cs:1.15,0.967737651181982) node[
scale=0.72,
anchor=base,
text=black,
rotate=0.0
]{4};
\draw (axis cs:2.15,0.887980806193362) node[
scale=0.72,
anchor=base,
text=black,
rotate=0.0
]{4};
\draw (axis cs:3.15,0.856433594909314) node[
scale=0.72,
anchor=base,
text=black,
rotate=0.0
]{4};
\draw (axis cs:4.15,0.913594930898864) node[
scale=0.72,
anchor=base,
text=black,
rotate=0.0
]{4};
\draw (axis cs:5.15,0.797919448124907) node[
scale=0.72,
anchor=base,
text=black,
rotate=0.0
]{4};
\draw (axis cs:6.15,0.736325317804138) node[
scale=0.72,
anchor=base,
text=black,
rotate=0.0
]{4};
\draw (axis cs:7.15,0.778544747721771) node[
scale=0.72,
anchor=base,
text=black,
rotate=0.0
]{4};
\end{axis}

\end{tikzpicture}

%% file: plots/tg097/gmres_s_step_Jacobi.tex
\begin{tikzpicture}

\definecolor{burlywood253190133}{RGB}{253,190,133}
\definecolor{coral25314160}{RGB}{253,141,60}
\definecolor{cornflowerblue107174214}{RGB}{107,174,214}
\definecolor{darkgray176}{RGB}{176,176,176}
\definecolor{lightgray204}{RGB}{204,204,204}
\definecolor{orangered217711}{RGB}{217,71,1}
\definecolor{powderblue189215231}{RGB}{189,215,231}
\definecolor{steelblue33113181}{RGB}{33,113,181}

\begin{axis}[
		legend cell align={left},
legend style={fill opacity=0.8, draw opacity=1, text opacity=1, draw=lightgray204},
tick align=outside,
tick pos=left,
x grid style={darkgray176},
xtick={0, 1, 2, 3, 4, 5, 6, 7},
xticklabels={1, 2, 3, 4, 5, 6, 7, 8},
width  = 0.49\linewidth,
height = 0.2\textheight,
xmin=-0.68, xmax=7.68,
xtick style={color=black},
y grid style={darkgray176},
ylabel={\small Normalized time},
ymin=0, ymax=1.15,
ytick style={color=black},
]
\draw[draw=none,fill=powderblue189215231] (axis cs:0,0) rectangle (axis cs:-0.3,0.229454291367302);
\draw[draw=none,fill=powderblue189215231] (axis cs:1,0) rectangle (axis cs:0.7,0.21620656478596);
\draw[draw=none,fill=powderblue189215231] (axis cs:2,0) rectangle (axis cs:1.7,0.207178634657311);
\draw[draw=none,fill=powderblue189215231] (axis cs:3,0) rectangle (axis cs:2.7,0.199026191526559);
\draw[draw=none,fill=powderblue189215231] (axis cs:4,0) rectangle (axis cs:3.7,0.16721942204528);
\draw[draw=none,fill=powderblue189215231] (axis cs:5,0) rectangle (axis cs:4.7,0.178744802303775);
\draw[draw=none,fill=powderblue189215231] (axis cs:6,0) rectangle (axis cs:5.7,0.142509966775606);
\draw[draw=none,fill=powderblue189215231] (axis cs:7,0) rectangle (axis cs:6.7,0.13997491198335);
\draw[draw=none,fill=cornflowerblue107174214] (axis cs:0,0.229454291367302) rectangle (axis cs:-0.3,0.909644564635418);
\draw[draw=none,fill=cornflowerblue107174214] (axis cs:1,0.21620656478596) rectangle (axis cs:0.7,0.855906070855127);
\draw[draw=none,fill=cornflowerblue107174214] (axis cs:2,0.207178634657311) rectangle (axis cs:1.7,0.781896001576411);
\draw[draw=none,fill=cornflowerblue107174214] (axis cs:3,0.199026191526559) rectangle (axis cs:2.7,0.613433981994347);
\draw[draw=none,fill=cornflowerblue107174214] (axis cs:4,0.16721942204528) rectangle (axis cs:3.7,0.565518511563733);
\draw[draw=none,fill=cornflowerblue107174214] (axis cs:5,0.178744802303775) rectangle (axis cs:4.7,0.608311813807368);
\draw[draw=none,fill=cornflowerblue107174214] (axis cs:6,0.142509966775606) rectangle (axis cs:5.7,0.451797492920415);
\draw[draw=none,fill=cornflowerblue107174214] (axis cs:7,0.13997491198335) rectangle (axis cs:6.7,0.441383515240024);
\draw[draw=none,fill=steelblue33113181] (axis cs:0,0.909644564635418) rectangle (axis cs:-0.3,1);
\draw[draw=none,fill=steelblue33113181] (axis cs:1,0.855906070855127) rectangle (axis cs:0.7,1);
\draw[draw=none,fill=steelblue33113181] (axis cs:2,0.781896001576411) rectangle (axis cs:1.7,1);
\draw[draw=none,fill=steelblue33113181] (axis cs:3,0.613433981994347) rectangle (axis cs:2.7,1);
\draw[draw=none,fill=steelblue33113181] (axis cs:4,0.565518511563733) rectangle (axis cs:3.7,1);
\draw[draw=none,fill=steelblue33113181] (axis cs:5,0.608311813807368) rectangle (axis cs:4.7,1);
\draw[draw=none,fill=steelblue33113181] (axis cs:6,0.451797492920414) rectangle (axis cs:5.7,1);
\draw[draw=none,fill=steelblue33113181] (axis cs:7,0.441383515240024) rectangle (axis cs:6.7,1);
\draw[draw=none,fill=burlywood253190133] (axis cs:0,0) rectangle (axis cs:0.3,0.237720298998341);
\draw[draw=none,fill=burlywood253190133] (axis cs:1,0) rectangle (axis cs:1.3,0.231776805928014);
\draw[draw=none,fill=burlywood253190133] (axis cs:2,0) rectangle (axis cs:2.3,0.215751709100889);
\draw[draw=none,fill=burlywood253190133] (axis cs:3,0) rectangle (axis cs:3.3,0.180808457101776);
\draw[draw=none,fill=burlywood253190133] (axis cs:4,0) rectangle (axis cs:4.3,0.183997517891937);
\draw[draw=none,fill=burlywood253190133] (axis cs:5,0) rectangle (axis cs:5.3,0.183534121007149);
\draw[draw=none,fill=burlywood253190133] (axis cs:6,0) rectangle (axis cs:6.3,0.156964043629133);
\draw[draw=none,fill=burlywood253190133] (axis cs:7,0) rectangle (axis cs:7.3,0.157136597539447);
\draw[draw=none,fill=coral25314160] (axis cs:0,0.237720298998341) rectangle (axis cs:0.3,0.918349617141893);
\draw[draw=none,fill=coral25314160] (axis cs:1,0.231776805928014) rectangle (axis cs:1.3,0.882025881258268);
\draw[draw=none,fill=coral25314160] (axis cs:2,0.215751709100889) rectangle (axis cs:2.3,0.790851329007631);
\draw[draw=none,fill=coral25314160] (axis cs:3,0.180808457101776) rectangle (axis cs:3.3,0.606245118703332);
\draw[draw=none,fill=coral25314160] (axis cs:4,0.183997517891937) rectangle (axis cs:4.3,0.585223621635228);
\draw[draw=none,fill=coral25314160] (axis cs:5,0.183534121007149) rectangle (axis cs:5.3,0.614089341880685);
\draw[draw=none,fill=coral25314160] (axis cs:6,0.156964043629133) rectangle (axis cs:6.3,0.468205307436902);
\draw[draw=none,fill=coral25314160] (axis cs:7,0.157136597539447) rectangle (axis cs:7.3,0.46181011517115);
\draw[draw=none,fill=orangered217711] (axis cs:0,0.918349617141893) rectangle (axis cs:0.3,0.980342753553355);
\draw[draw=none,fill=orangered217711] (axis cs:1,0.882025881258268) rectangle (axis cs:1.3,0.951027976402723);
\draw[draw=none,fill=orangered217711] (axis cs:2,0.790851329007631) rectangle (axis cs:2.3,0.907740470098547);
\draw[draw=none,fill=orangered217711] (axis cs:3,0.606245118703332) rectangle (axis cs:3.3,0.814224302730258);
\draw[draw=none,fill=orangered217711] (axis cs:4,0.585223621635228) rectangle (axis cs:4.3,0.777428599991084);
\draw[draw=none,fill=orangered217711] (axis cs:5,0.614089341880685) rectangle (axis cs:5.3,0.768250003591761);
\draw[draw=none,fill=orangered217711] (axis cs:6,0.468205307436902) rectangle (axis cs:6.3,0.683034501921718);
\draw[draw=none,fill=orangered217711] (axis cs:7,0.46181011517115) rectangle (axis cs:7.3,0.691832142978374);
\draw (axis cs:0.15,1.03034275355335) node[
scale=0.72,
anchor=base,
text=black,
rotate=0.0
]{2};
\draw (axis cs:1.15,1.00102797640272) node[
scale=0.72,
anchor=base,
text=black,
rotate=0.0
]{2};
\draw (axis cs:2.15,0.957740470098547) node[
scale=0.72,
anchor=base,
text=black,
rotate=0.0
]{4};
\draw (axis cs:3.15,0.864224302730258) node[
scale=0.72,
anchor=base,
text=black,
rotate=0.0
]{4};
\draw (axis cs:4.15,0.827428599991084) node[
scale=0.72,
anchor=base,
text=black,
rotate=0.0
]{4};
\draw (axis cs:5.15,0.818250003591761) node[
scale=0.72,
anchor=base,
text=black,
rotate=0.0
]{4};
\draw (axis cs:6.15,0.733034501921718) node[
scale=0.72,
anchor=base,
text=black,
rotate=0.0
]{4};
\draw (axis cs:7.15,0.741832142978374) node[
scale=0.72,
anchor=base,
text=black,
rotate=0.0
]{4};
\end{axis}

\end{tikzpicture}

%% file: plots/horeka/gmres_s_step_GS2_1.tex
\begin{tikzpicture}

\definecolor{burlywood253190133}{RGB}{253,190,133}
\definecolor{coral25314160}{RGB}{253,141,60}
\definecolor{cornflowerblue107174214}{RGB}{107,174,214}
\definecolor{darkgray176}{RGB}{176,176,176}
\definecolor{lightgray204}{RGB}{204,204,204}
\definecolor{orangered217711}{RGB}{217,71,1}
\definecolor{powderblue189215231}{RGB}{189,215,231}
\definecolor{steelblue33113181}{RGB}{33,113,181}

\begin{axis}[
	legend cell align={left},
	legend style={fill opacity=0.8, draw opacity=1, text opacity=1, draw=lightgray204},
	tick align=outside,
	tick pos=left,
	x grid style={darkgray176},
	xtick={0, 1, 2, 3, 4, 5, 6, 7},
	xticklabels={1, 2, 3, 4, 5, 6, 7, 8},
	width  = 0.49\linewidth,
	height = 0.2\textheight,
	xmin=-0.68, xmax=7.68,
	xtick style={color=black},
	y grid style={darkgray176},
	ylabel={\small Normalized time},
	ymin=0, ymax=1.1,
	ytick style={color=black}
	]
\draw[draw=none,fill=powderblue189215231] (axis cs:0,0) rectangle (axis cs:-0.3,0.149930296762672);
\draw[draw=none,fill=powderblue189215231] (axis cs:1,0) rectangle (axis cs:0.7,0.149865952810609);
\draw[draw=none,fill=powderblue189215231] (axis cs:2,0) rectangle (axis cs:1.7,0.129246734418198);
\draw[draw=none,fill=powderblue189215231] (axis cs:3,0) rectangle (axis cs:2.7,0.11306078950874);
\draw[draw=none,fill=powderblue189215231] (axis cs:4,0) rectangle (axis cs:3.7,0.0783948309245592);
\draw[draw=none,fill=powderblue189215231] (axis cs:5,0) rectangle (axis cs:4.7,0.108532620239308);
\draw[draw=none,fill=powderblue189215231] (axis cs:6,0) rectangle (axis cs:5.7,0.0864679310231214);
\draw[draw=none,fill=powderblue189215231] (axis cs:7,0) rectangle (axis cs:6.7,0.0850534336479496);
\draw[draw=none,fill=cornflowerblue107174214] (axis cs:0,0.149930296762672) rectangle (axis cs:-0.3,0.536769596632377);
\draw[draw=none,fill=cornflowerblue107174214] (axis cs:1,0.149865952810609) rectangle (axis cs:0.7,0.49787078089957);
\draw[draw=none,fill=cornflowerblue107174214] (axis cs:2,0.129246734418198) rectangle (axis cs:1.7,0.414553769155531);
\draw[draw=none,fill=cornflowerblue107174214] (axis cs:3,0.11306078950874) rectangle (axis cs:2.7,0.273619395805272);
\draw[draw=none,fill=cornflowerblue107174214] (axis cs:4,0.0783948309245592) rectangle (axis cs:3.7,0.446703232536755);
\draw[draw=none,fill=cornflowerblue107174214] (axis cs:5,0.108532620239308) rectangle (axis cs:4.7,0.295619954287914);
\draw[draw=none,fill=cornflowerblue107174214] (axis cs:6,0.0864679310231214) rectangle (axis cs:5.7,0.20125627761859);
\draw[draw=none,fill=cornflowerblue107174214] (axis cs:7,0.0850534336479496) rectangle (axis cs:6.7,0.19367919234591);
\draw[draw=none,fill=steelblue33113181] (axis cs:0,0.536769596632377) rectangle (axis cs:-0.3,1);
\draw[draw=none,fill=steelblue33113181] (axis cs:1,0.49787078089957) rectangle (axis cs:0.7,1);
\draw[draw=none,fill=steelblue33113181] (axis cs:2,0.414553769155531) rectangle (axis cs:1.7,1);
\draw[draw=none,fill=steelblue33113181] (axis cs:3,0.273619395805272) rectangle (axis cs:2.7,1);
\draw[draw=none,fill=steelblue33113181] (axis cs:4,0.446703232536755) rectangle (axis cs:3.7,1);
\draw[draw=none,fill=steelblue33113181] (axis cs:5,0.295619954287914) rectangle (axis cs:4.7,1);
\draw[draw=none,fill=steelblue33113181] (axis cs:6,0.20125627761859) rectangle (axis cs:5.7,1);
\draw[draw=none,fill=steelblue33113181] (axis cs:7,0.19367919234591) rectangle (axis cs:6.7,1);
\draw[draw=none,fill=burlywood253190133] (axis cs:0,0) rectangle (axis cs:0.3,0.151581306776248);
\draw[draw=none,fill=burlywood253190133] (axis cs:1,0) rectangle (axis cs:1.3,0.1504734475025);
\draw[draw=none,fill=burlywood253190133] (axis cs:2,0) rectangle (axis cs:2.3,0.128769894568267);
\draw[draw=none,fill=burlywood253190133] (axis cs:3,0) rectangle (axis cs:3.3,0.128644854279836);
\draw[draw=none,fill=burlywood253190133] (axis cs:4,0) rectangle (axis cs:4.3,0.0896540962901245);
\draw[draw=none,fill=burlywood253190133] (axis cs:5,0) rectangle (axis cs:5.3,0.122946087695316);
\draw[draw=none,fill=burlywood253190133] (axis cs:6,0) rectangle (axis cs:6.3,0.106894169667843);
\draw[draw=none,fill=burlywood253190133] (axis cs:7,0) rectangle (axis cs:7.3,0.108621034469545);
\draw[draw=none,fill=coral25314160] (axis cs:0,0.151581306776248) rectangle (axis cs:0.3,0.536598299787701);
\draw[draw=none,fill=coral25314160] (axis cs:1,0.1504734475025) rectangle (axis cs:1.3,0.492722014412563);
\draw[draw=none,fill=coral25314160] (axis cs:2,0.128769894568267) rectangle (axis cs:2.3,0.403267702518768);
\draw[draw=none,fill=coral25314160] (axis cs:3,0.128644854279836) rectangle (axis cs:3.3,0.286687215997395);
\draw[draw=none,fill=coral25314160] (axis cs:4,0.0896540962901245) rectangle (axis cs:4.3,0.455116523043073);
\draw[draw=none,fill=coral25314160] (axis cs:5,0.122946087695316) rectangle (axis cs:5.3,0.31077774086439);
\draw[draw=none,fill=coral25314160] (axis cs:6,0.106894169667843) rectangle (axis cs:6.3,0.221057057989476);
\draw[draw=none,fill=coral25314160] (axis cs:7,0.108621034469545) rectangle (axis cs:7.3,0.222228229574405);
\draw[draw=none,fill=orangered217711] (axis cs:0,0.536598299787701) rectangle (axis cs:0.3,0.81510510154805);
\draw[draw=none,fill=orangered217711] (axis cs:1,0.492722014412563) rectangle (axis cs:1.3,0.896542658078086);
\draw[draw=none,fill=orangered217711] (axis cs:2,0.403267702518768) rectangle (axis cs:2.3,0.752615228799585);
\draw[draw=none,fill=orangered217711] (axis cs:3,0.286687215997395) rectangle (axis cs:3.3,0.836581015674723);
\draw[draw=none,fill=orangered217711] (axis cs:4,0.455116523043073) rectangle (axis cs:4.3,0.872218672931737);
\draw[draw=none,fill=orangered217711] (axis cs:5,0.31077774086439) rectangle (axis cs:5.3,0.726469077171369);
\draw[draw=none,fill=orangered217711] (axis cs:6,0.221057057989476) rectangle (axis cs:6.3,0.797919457727676);
\draw[draw=none,fill=orangered217711] (axis cs:7,0.222228229574405) rectangle (axis cs:7.3,0.838855518922372);
\draw (axis cs:0.15,0.86510510154805) node[
scale=0.72,
anchor=base,
text=black,
rotate=0.0
]{2};
\draw (axis cs:1.15,0.946542658078086) node[
scale=0.72,
anchor=base,
text=black,
rotate=0.0
]{2};
\draw (axis cs:2.15,0.802615228799585) node[
scale=0.72,
anchor=base,
text=black,
rotate=0.0
]{2};
\draw (axis cs:3.15,0.886581015674723) node[
scale=0.72,
anchor=base,
text=black,
rotate=0.0
]{2};
\draw (axis cs:4.15,0.922218672931737) node[
scale=0.72,
anchor=base,
text=black,
rotate=0.0
]{2};
\draw (axis cs:5.15,0.776469077171369) node[
scale=0.72,
anchor=base,
text=black,
rotate=0.0
]{2};
\draw (axis cs:6.15,0.847919457727677) node[
scale=0.72,
anchor=base,
text=black,
rotate=0.0
]{1};
\draw (axis cs:7.15,0.888855518922372) node[
scale=0.72,
anchor=base,
text=black,
rotate=0.0
]{2};
\end{axis}

\end{tikzpicture}

%% file: plots/tg097/gmres_s_step_GS2_1.tex
\begin{tikzpicture}

\definecolor{burlywood253190133}{RGB}{253,190,133}
\definecolor{coral25314160}{RGB}{253,141,60}
\definecolor{cornflowerblue107174214}{RGB}{107,174,214}
\definecolor{darkgray176}{RGB}{176,176,176}
\definecolor{lightgray204}{RGB}{204,204,204}
\definecolor{orangered217711}{RGB}{217,71,1}
\definecolor{powderblue189215231}{RGB}{189,215,231}
\definecolor{steelblue33113181}{RGB}{33,113,181}

\begin{axis}[
		legend cell align={left},
legend style={fill opacity=0.8, draw opacity=1, text opacity=1, draw=lightgray204},
tick align=outside,
tick pos=left,
x grid style={darkgray176},
xtick={0, 1, 2, 3, 4, 5, 6, 7},
xticklabels={1, 2, 3, 4, 5, 6, 7, 8},
width  = 0.49\linewidth,
height = 0.2\textheight,
xmin=-0.68, xmax=7.68,
xtick style={color=black},
y grid style={darkgray176},
ylabel={\small Normalized time},
ymin=0, ymax=1.1,
ytick style={color=black},
]
\draw[draw=none,fill=powderblue189215231] (axis cs:0,0) rectangle (axis cs:-0.3,0.212235428532621);
\draw[draw=none,fill=powderblue189215231] (axis cs:1,0) rectangle (axis cs:0.7,0.192783508595114);
\draw[draw=none,fill=powderblue189215231] (axis cs:2,0) rectangle (axis cs:1.7,0.183365966780589);
\draw[draw=none,fill=powderblue189215231] (axis cs:3,0) rectangle (axis cs:2.7,0.14031496561995);
\draw[draw=none,fill=powderblue189215231] (axis cs:4,0) rectangle (axis cs:3.7,0.144038116778659);
\draw[draw=none,fill=powderblue189215231] (axis cs:5,0) rectangle (axis cs:4.7,0.150411221616356);
\draw[draw=none,fill=powderblue189215231] (axis cs:6,0) rectangle (axis cs:5.7,0.119947238359667);
\draw[draw=none,fill=powderblue189215231] (axis cs:7,0) rectangle (axis cs:6.7,0.11862669345085);
\draw[draw=none,fill=cornflowerblue107174214] (axis cs:0,0.212235428532621) rectangle (axis cs:-0.3,0.822848348478811);
\draw[draw=none,fill=cornflowerblue107174214] (axis cs:1,0.192783508595114) rectangle (axis cs:0.7,0.717734534302029);
\draw[draw=none,fill=cornflowerblue107174214] (axis cs:2,0.183365966780589) rectangle (axis cs:1.7,0.651404153799671);
\draw[draw=none,fill=cornflowerblue107174214] (axis cs:3,0.14031496561995) rectangle (axis cs:2.7,0.461706644697234);
\draw[draw=none,fill=cornflowerblue107174214] (axis cs:4,0.144038116778659) rectangle (axis cs:3.7,0.452005557966345);
\draw[draw=none,fill=cornflowerblue107174214] (axis cs:5,0.150411221616356) rectangle (axis cs:4.7,0.491886396802291);
\draw[draw=none,fill=cornflowerblue107174214] (axis cs:6,0.119947238359667) rectangle (axis cs:5.7,0.351064392813669);
\draw[draw=none,fill=cornflowerblue107174214] (axis cs:7,0.11862669345085) rectangle (axis cs:6.7,0.34134720008782);
\draw[draw=none,fill=steelblue33113181] (axis cs:0,0.822848348478811) rectangle (axis cs:-0.3,1);
\draw[draw=none,fill=steelblue33113181] (axis cs:1,0.717734534302029) rectangle (axis cs:0.7,1);
\draw[draw=none,fill=steelblue33113181] (axis cs:2,0.651404153799671) rectangle (axis cs:1.7,1);
\draw[draw=none,fill=steelblue33113181] (axis cs:3,0.461706644697234) rectangle (axis cs:2.7,1);
\draw[draw=none,fill=steelblue33113181] (axis cs:4,0.452005557966345) rectangle (axis cs:3.7,1);
\draw[draw=none,fill=steelblue33113181] (axis cs:5,0.491886396802291) rectangle (axis cs:4.7,1);
\draw[draw=none,fill=steelblue33113181] (axis cs:6,0.351064392813669) rectangle (axis cs:5.7,1);
\draw[draw=none,fill=steelblue33113181] (axis cs:7,0.34134720008782) rectangle (axis cs:6.7,1);
\draw[draw=none,fill=burlywood253190133] (axis cs:0,0) rectangle (axis cs:0.3,0.220622617446112);
\draw[draw=none,fill=burlywood253190133] (axis cs:1,0) rectangle (axis cs:1.3,0.199658264557963);
\draw[draw=none,fill=burlywood253190133] (axis cs:2,0) rectangle (axis cs:2.3,0.192136581767385);
\draw[draw=none,fill=burlywood253190133] (axis cs:3,0) rectangle (axis cs:3.3,0.158935327481859);
\draw[draw=none,fill=burlywood253190133] (axis cs:4,0) rectangle (axis cs:4.3,0.160690099543703);
\draw[draw=none,fill=burlywood253190133] (axis cs:5,0) rectangle (axis cs:5.3,0.160827337932894);
\draw[draw=none,fill=burlywood253190133] (axis cs:6,0) rectangle (axis cs:6.3,0.137001542952291);
\draw[draw=none,fill=burlywood253190133] (axis cs:7,0) rectangle (axis cs:7.3,0.136313685715158);
\draw[draw=none,fill=coral25314160] (axis cs:0,0.220622617446112) rectangle (axis cs:0.3,0.844818029559897);
\draw[draw=none,fill=coral25314160] (axis cs:1,0.199658264557963) rectangle (axis cs:1.3,0.747302659462634);
\draw[draw=none,fill=coral25314160] (axis cs:2,0.192136581767385) rectangle (axis cs:2.3,0.668885741080082);
\draw[draw=none,fill=coral25314160] (axis cs:3,0.158935327481859) rectangle (axis cs:3.3,0.484940587471086);
\draw[draw=none,fill=coral25314160] (axis cs:4,0.160690099543703) rectangle (axis cs:4.3,0.475057091886158);
\draw[draw=none,fill=coral25314160] (axis cs:5,0.160827337932894) rectangle (axis cs:5.3,0.507707330936995);
\draw[draw=none,fill=coral25314160] (axis cs:6,0.137001542952291) rectangle (axis cs:6.3,0.371579144496057);
\draw[draw=none,fill=coral25314160] (axis cs:7,0.136313685715158) rectangle (axis cs:7.3,0.365141419017981);
\draw[draw=none,fill=orangered217711] (axis cs:0,0.844818029559897) rectangle (axis cs:0.3,0.941378563032251);
\draw[draw=none,fill=orangered217711] (axis cs:1,0.747302659462634) rectangle (axis cs:1.3,0.862101072261695);
\draw[draw=none,fill=orangered217711] (axis cs:2,0.668885741080082) rectangle (axis cs:2.3,0.839197805776709);
\draw[draw=none,fill=orangered217711] (axis cs:3,0.484940587471086) rectangle (axis cs:3.3,0.771729142241516);
\draw[draw=none,fill=orangered217711] (axis cs:4,0.475057091886158) rectangle (axis cs:4.3,0.758334538728242);
\draw[draw=none,fill=orangered217711] (axis cs:5,0.507707330936995) rectangle (axis cs:5.3,0.697458946157691);
\draw[draw=none,fill=orangered217711] (axis cs:6,0.371579144496057) rectangle (axis cs:6.3,0.631772196732524);
\draw[draw=none,fill=orangered217711] (axis cs:7,0.365141419017981) rectangle (axis cs:7.3,0.660924456311877);
\draw (axis cs:0.15,0.991378563032251) node[
scale=0.72,
anchor=base,
text=black,
rotate=0.0
]{1};
\draw (axis cs:1.15,0.912101072261695) node[
scale=0.72,
anchor=base,
text=black,
rotate=0.0
]{1};
\draw (axis cs:2.15,0.889197805776709) node[
scale=0.72,
anchor=base,
text=black,
rotate=0.0
]{4};
\draw (axis cs:3.15,0.821729142241516) node[
scale=0.72,
anchor=base,
text=black,
rotate=0.0
]{4};
\draw (axis cs:4.15,0.808334538728243) node[
scale=0.72,
anchor=base,
text=black,
rotate=0.0
]{4};
\draw (axis cs:5.15,0.747458946157691) node[
scale=0.72,
anchor=base,
text=black,
rotate=0.0
]{4};
\draw (axis cs:6.15,0.681772196732524) node[
scale=0.72,
anchor=base,
text=black,
rotate=0.0
]{4};
\draw (axis cs:7.15,0.710924456311877) node[
scale=0.72,
anchor=base,
text=black,
rotate=0.0
]{4};
\end{axis}

\end{tikzpicture}

%% file: plots/horeka/gmres_poly_none.tex
\begin{tikzpicture}

\definecolor{burlywood253190133}{RGB}{253,190,133}
\definecolor{coral25314160}{RGB}{253,141,60}
\definecolor{cornflowerblue107174214}{RGB}{107,174,214}
\definecolor{darkgray176}{RGB}{176,176,176}
\definecolor{lightgray204}{RGB}{204,204,204}
\definecolor{orangered217711}{RGB}{217,71,1}
\definecolor{powderblue189215231}{RGB}{189,215,231}
\definecolor{steelblue33113181}{RGB}{33,113,181}

\begin{axis}[
	legend cell align={left},
legend style={fill opacity=0.8, draw opacity=1, text opacity=1, draw=lightgray204},
tick align=outside,
tick pos=left,
x grid style={darkgray176},
xtick={0, 1, 2, 3, 4, 5, 6, 7},
xticklabels={1, 2, 3, 4, 5, 6, 7, 8},
width  = 0.49\linewidth,
height = 0.2\textheight,
xmin=-0.68, xmax=7.68,
xtick style={color=black},
y grid style={darkgray176},
ylabel={\small Normalized time},
ymin=0, ymax=1.05,
ytick style={color=black}
]
\draw[draw=none,fill=powderblue189215231] (axis cs:0,0) rectangle (axis cs:-0.3,0.0489399586179578);
\draw[draw=none,fill=powderblue189215231] (axis cs:1,0) rectangle (axis cs:0.7,0.0429111104120785);
\draw[draw=none,fill=powderblue189215231] (axis cs:2,0) rectangle (axis cs:1.7,0.0827880196733246);
\draw[draw=none,fill=powderblue189215231] (axis cs:3,0) rectangle (axis cs:2.7,0.0079932137321337);
\draw[draw=none,fill=powderblue189215231] (axis cs:4,0) rectangle (axis cs:3.7,0.0101288565871879);
\draw[draw=none,fill=powderblue189215231] (axis cs:5,0) rectangle (axis cs:4.7,0.0152214637621114);
\draw[draw=none,fill=powderblue189215231] (axis cs:6,0) rectangle (axis cs:5.7,0.00956820469454421);
\draw[draw=none,fill=powderblue189215231] (axis cs:7,0) rectangle (axis cs:6.7,0.00860919969129599);
\draw[draw=none,fill=cornflowerblue107174214] (axis cs:0,0.0489399586179578) rectangle (axis cs:-0.3,0.0993999096972312);
\draw[draw=none,fill=cornflowerblue107174214] (axis cs:1,0.0429111104120785) rectangle (axis cs:0.7,0.0836596074186584);
\draw[draw=none,fill=cornflowerblue107174214] (axis cs:2,0.0827880196733246) rectangle (axis cs:1.7,0.103242524442689);
\draw[draw=none,fill=cornflowerblue107174214] (axis cs:3,0.0079932137321337) rectangle (axis cs:2.7,0.0129982426216811);
\draw[draw=none,fill=cornflowerblue107174214] (axis cs:4,0.0101288565871879) rectangle (axis cs:3.7,0.0167781329984498);
\draw[draw=none,fill=cornflowerblue107174214] (axis cs:5,0.0152214637621114) rectangle (axis cs:4.7,0.0272214595989445);
\draw[draw=none,fill=cornflowerblue107174214] (axis cs:6,0.00956820469454421) rectangle (axis cs:5.7,0.0152229068918686);
\draw[draw=none,fill=cornflowerblue107174214] (axis cs:7,0.00860919969129599) rectangle (axis cs:6.7,0.0138164161001058);
\draw[draw=none,fill=steelblue33113181] (axis cs:0,0.0993999096972312) rectangle (axis cs:-0.3,1);
\draw[draw=none,fill=steelblue33113181] (axis cs:1,0.0836596074186584) rectangle (axis cs:0.7,1);
\draw[draw=none,fill=steelblue33113181] (axis cs:2,0.103242524442689) rectangle (axis cs:1.7,1);
\draw[draw=none,fill=steelblue33113181] (axis cs:3,0.0129982426216811) rectangle (axis cs:2.7,1);
\draw[draw=none,fill=steelblue33113181] (axis cs:4,0.0167781329984498) rectangle (axis cs:3.7,1);
\draw[draw=none,fill=steelblue33113181] (axis cs:5,0.0272214595989445) rectangle (axis cs:4.7,1);
\draw[draw=none,fill=steelblue33113181] (axis cs:6,0.0152229068918686) rectangle (axis cs:5.7,1);
\draw[draw=none,fill=steelblue33113181] (axis cs:7,0.0138164161001058) rectangle (axis cs:6.7,1);
\draw[draw=none,fill=burlywood253190133] (axis cs:0,0) rectangle (axis cs:0.3,0.0565924592840429);
\draw[draw=none,fill=burlywood253190133] (axis cs:1,0) rectangle (axis cs:1.3,0.0521884113949259);
\draw[draw=none,fill=burlywood253190133] (axis cs:2,0) rectangle (axis cs:2.3,0.0985869495278101);
\draw[draw=none,fill=burlywood253190133] (axis cs:3,0) rectangle (axis cs:3.3,0.00954245452953778);
\draw[draw=none,fill=burlywood253190133] (axis cs:4,0) rectangle (axis cs:4.3,0.0120758823959777);
\draw[draw=none,fill=burlywood253190133] (axis cs:5,0) rectangle (axis cs:5.3,0.019506459899277);
\draw[draw=none,fill=burlywood253190133] (axis cs:6,0) rectangle (axis cs:6.3,0.0124499841892317);
\draw[draw=none,fill=burlywood253190133] (axis cs:7,0) rectangle (axis cs:7.3,0.011899480418009);
\draw[draw=none,fill=coral25314160] (axis cs:0,0.0565924592840429) rectangle (axis cs:0.3,0.106534899210633);
\draw[draw=none,fill=coral25314160] (axis cs:1,0.0521884113949259) rectangle (axis cs:1.3,0.092024943332665);
\draw[draw=none,fill=coral25314160] (axis cs:2,0.0985869495278101) rectangle (axis cs:2.3,0.119698018647796);
\draw[draw=none,fill=coral25314160] (axis cs:3,0.00954245452953778) rectangle (axis cs:3.3,0.0144493949132096);
\draw[draw=none,fill=coral25314160] (axis cs:4,0.0120758823959777) rectangle (axis cs:4.3,0.0187773242671347);
\draw[draw=none,fill=coral25314160] (axis cs:5,0.019506459899277) rectangle (axis cs:5.3,0.0316701520489559);
\draw[draw=none,fill=coral25314160] (axis cs:6,0.0124499841892317) rectangle (axis cs:6.3,0.0182074091447917);
\draw[draw=none,fill=coral25314160] (axis cs:7,0.011899480418009) rectangle (axis cs:7.3,0.0170941216506951);
\draw[draw=none,fill=orangered217711] (axis cs:0,0.106534899210633) rectangle (axis cs:0.3,0.559503524598832);
\draw[draw=none,fill=orangered217711] (axis cs:1,0.092024943332665) rectangle (axis cs:1.3,0.599631026418641);
\draw[draw=none,fill=orangered217711] (axis cs:2,0.119698018647796) rectangle (axis cs:2.3,0.548981676562541);
\draw[draw=none,fill=orangered217711] (axis cs:3,0.0144493949132096) rectangle (axis cs:3.3,0.595134607241422);
\draw[draw=none,fill=orangered217711] (axis cs:4,0.0187773242671347) rectangle (axis cs:4.3,0.586622618228355);
\draw[draw=none,fill=orangered217711] (axis cs:5,0.0316701520489559) rectangle (axis cs:5.3,0.398904348076485);
\draw[draw=none,fill=orangered217711] (axis cs:6,0.0182074091447917) rectangle (axis cs:6.3,0.436060074798471);
\draw[draw=none,fill=orangered217711] (axis cs:7,0.0170941216506951) rectangle (axis cs:7.3,0.523392582264523);
\draw (axis cs:0.15,0.609503524598832) node[
scale=0.72,
anchor=base,
text=black,
rotate=0.0
]{7};
\draw (axis cs:1.15,0.649631026418641) node[
scale=0.72,
anchor=base,
text=black,
rotate=0.0
]{8};
\draw (axis cs:2.15,0.598981676562541) node[
scale=0.72,
anchor=base,
text=black,
rotate=0.0
]{7};
\draw (axis cs:3.15,0.645134607241422) node[
scale=0.72,
anchor=base,
text=black,
rotate=0.0
]{5};
\draw (axis cs:4.15,0.636622618228355) node[
scale=0.72,
anchor=base,
text=black,
rotate=0.0
]{5};
\draw (axis cs:5.15,0.448904348076485) node[
scale=0.72,
anchor=base,
text=black,
rotate=0.0
]{8};
\draw (axis cs:6.15,0.486060074798471) node[
scale=0.72,
anchor=base,
text=black,
rotate=0.0
]{5};
\draw (axis cs:7.15,0.573392582264523) node[
scale=0.72,
anchor=base,
text=black,
rotate=0.0
]{4};
\end{axis}

\end{tikzpicture}

%% file: plots/tg097/gmres_poly_none.tex
\begin{tikzpicture}

\definecolor{burlywood253190133}{RGB}{253,190,133}
\definecolor{coral25314160}{RGB}{253,141,60}
\definecolor{cornflowerblue107174214}{RGB}{107,174,214}
\definecolor{darkgray176}{RGB}{176,176,176}
\definecolor{lightgray204}{RGB}{204,204,204}
\definecolor{orangered217711}{RGB}{217,71,1}
\definecolor{powderblue189215231}{RGB}{189,215,231}
\definecolor{steelblue33113181}{RGB}{33,113,181}

\begin{axis}[
		legend cell align={left},
legend style={fill opacity=0.8, draw opacity=1, text opacity=1, draw=lightgray204},
tick align=outside,
tick pos=left,
x grid style={darkgray176},
xtick={0, 1, 2, 3, 4, 5, 6, 7},
xticklabels={1, 2, 3, 4, 5, 6, 7, 8},
width  = 0.49\linewidth,
height = 0.2\textheight,
xmin=-0.68, xmax=7.68,
xtick style={color=black},
y grid style={darkgray176},
ylabel={\small Normalized time},
ymin=0, ymax=1.05,
ytick style={color=black},
]
\draw[draw=none,fill=powderblue189215231] (axis cs:0,0) rectangle (axis cs:-0.3,0.0588943357331463);
\draw[draw=none,fill=powderblue189215231] (axis cs:1,0) rectangle (axis cs:0.7,0.0385641691824735);
\draw[draw=none,fill=powderblue189215231] (axis cs:2,0) rectangle (axis cs:1.7,0.106219874541256);
\draw[draw=none,fill=powderblue189215231] (axis cs:3,0) rectangle (axis cs:2.7,0.00821941408088296);
\draw[draw=none,fill=powderblue189215231] (axis cs:4,0) rectangle (axis cs:3.7,0.00911502869987132);
\draw[draw=none,fill=powderblue189215231] (axis cs:5,0) rectangle (axis cs:4.7,0.0140275983573406);
\draw[draw=none,fill=powderblue189215231] (axis cs:6,0) rectangle (axis cs:5.7,0.00915563452788127);
\draw[draw=none,fill=powderblue189215231] (axis cs:7,0) rectangle (axis cs:6.7,0.00855974024384967);
\draw[draw=none,fill=cornflowerblue107174214] (axis cs:0,0.0588943357331463) rectangle (axis cs:-0.3,0.119753379648506);
\draw[draw=none,fill=cornflowerblue107174214] (axis cs:1,0.0385641691824735) rectangle (axis cs:0.7,0.0727345319425532);
\draw[draw=none,fill=cornflowerblue107174214] (axis cs:2,0.106219874541256) rectangle (axis cs:1.7,0.13115755370296);
\draw[draw=none,fill=cornflowerblue107174214] (axis cs:3,0.00821941408088296) rectangle (axis cs:2.7,0.0117396793470236);
\draw[draw=none,fill=cornflowerblue107174214] (axis cs:4,0.00911502869987132) rectangle (axis cs:3.7,0.0134287741601448);
\draw[draw=none,fill=cornflowerblue107174214] (axis cs:5,0.0140275983573406) rectangle (axis cs:4.7,0.0235254616324589);
\draw[draw=none,fill=cornflowerblue107174214] (axis cs:6,0.00915563452788127) rectangle (axis cs:5.7,0.0136878110562186);
\draw[draw=none,fill=cornflowerblue107174214] (axis cs:7,0.00855974024384967) rectangle (axis cs:6.7,0.0127299944225542);
\draw[draw=none,fill=steelblue33113181] (axis cs:0,0.119753379648506) rectangle (axis cs:-0.3,1);
\draw[draw=none,fill=steelblue33113181] (axis cs:1,0.0727345319425532) rectangle (axis cs:0.7,1);
\draw[draw=none,fill=steelblue33113181] (axis cs:2,0.13115755370296) rectangle (axis cs:1.7,1);
\draw[draw=none,fill=steelblue33113181] (axis cs:3,0.0117396793470236) rectangle (axis cs:2.7,1);
\draw[draw=none,fill=steelblue33113181] (axis cs:4,0.0134287741601448) rectangle (axis cs:3.7,1);
\draw[draw=none,fill=steelblue33113181] (axis cs:5,0.0235254616324589) rectangle (axis cs:4.7,1);
\draw[draw=none,fill=steelblue33113181] (axis cs:6,0.0136878110562186) rectangle (axis cs:5.7,1);
\draw[draw=none,fill=steelblue33113181] (axis cs:7,0.0127299944225542) rectangle (axis cs:6.7,1);
\draw[draw=none,fill=burlywood253190133] (axis cs:0,0) rectangle (axis cs:0.3,0.0623739996221263);
\draw[draw=none,fill=burlywood253190133] (axis cs:1,0) rectangle (axis cs:1.3,0.0437647607027059);
\draw[draw=none,fill=burlywood253190133] (axis cs:2,0) rectangle (axis cs:2.3,0.122239270698507);
\draw[draw=none,fill=burlywood253190133] (axis cs:3,0) rectangle (axis cs:3.3,0.0106101246948655);
\draw[draw=none,fill=burlywood253190133] (axis cs:4,0) rectangle (axis cs:4.3,0.0113392325985642);
\draw[draw=none,fill=burlywood253190133] (axis cs:5,0) rectangle (axis cs:5.3,0.0170727439939954);
\draw[draw=none,fill=burlywood253190133] (axis cs:6,0) rectangle (axis cs:6.3,0.0118914933256972);
\draw[draw=none,fill=burlywood253190133] (axis cs:7,0) rectangle (axis cs:7.3,0.0108268585171456);
\draw[draw=none,fill=coral25314160] (axis cs:0,0.0623739996221263) rectangle (axis cs:0.3,0.121953687176109);
\draw[draw=none,fill=coral25314160] (axis cs:1,0.0437647607027059) rectangle (axis cs:1.3,0.0798221952429843);
\draw[draw=none,fill=coral25314160] (axis cs:2,0.122239270698507) rectangle (axis cs:2.3,0.147041983236938);
\draw[draw=none,fill=coral25314160] (axis cs:3,0.0106101246948655) rectangle (axis cs:3.3,0.0140838783398982);
\draw[draw=none,fill=coral25314160] (axis cs:4,0.0113392325985642) rectangle (axis cs:4.3,0.0156769839009829);
\draw[draw=none,fill=coral25314160] (axis cs:5,0.0170727439939954) rectangle (axis cs:5.3,0.0264886495281642);
\draw[draw=none,fill=coral25314160] (axis cs:6,0.0118914933256972) rectangle (axis cs:6.3,0.0164757123799186);
\draw[draw=none,fill=coral25314160] (axis cs:7,0.0108268585171456) rectangle (axis cs:7.3,0.0150298915590567);
\draw[draw=none,fill=orangered217711] (axis cs:0,0.121953687176109) rectangle (axis cs:0.3,0.510445388448747);
\draw[draw=none,fill=orangered217711] (axis cs:1,0.0798221952429843) rectangle (axis cs:1.3,0.358633336059104);
\draw[draw=none,fill=orangered217711] (axis cs:2,0.147041983236938) rectangle (axis cs:2.3,0.617982060028794);
\draw[draw=none,fill=orangered217711] (axis cs:3,0.0140838783398982) rectangle (axis cs:3.3,0.439849929058482);
\draw[draw=none,fill=orangered217711] (axis cs:4,0.0156769839009829) rectangle (axis cs:4.3,0.34809840487337);
\draw[draw=none,fill=orangered217711] (axis cs:5,0.0264886495281642) rectangle (axis cs:5.3,0.290946655888108);
\draw[draw=none,fill=orangered217711] (axis cs:6,0.0164757123799186) rectangle (axis cs:6.3,0.279195336595725);
\draw[draw=none,fill=orangered217711] (axis cs:7,0.0150298915590567) rectangle (axis cs:7.3,0.304799934673193);
\draw (axis cs:0.15,0.560445388448747) node[
scale=0.72,
anchor=base,
text=black,
rotate=0.0
]{3};
\draw (axis cs:1.15,0.408633336059104) node[
scale=0.72,
anchor=base,
text=black,
rotate=0.0
]{3};
\draw (axis cs:2.15,0.667982060028794) node[
scale=0.72,
anchor=base,
text=black,
rotate=0.0
]{8};
\draw (axis cs:3.15,0.489849929058482) node[
scale=0.72,
anchor=base,
text=black,
rotate=0.0
]{8};
\draw (axis cs:4.15,0.39809840487337) node[
scale=0.72,
anchor=base,
text=black,
rotate=0.0
]{8};
\draw (axis cs:5.15,0.340946655888108) node[
scale=0.72,
anchor=base,
text=black,
rotate=0.0
]{8};
\draw (axis cs:6.15,0.329195336595725) node[
scale=0.72,
anchor=base,
text=black,
rotate=0.0
]{8};
\draw (axis cs:7.15,0.354799934673193) node[
scale=0.72,
anchor=base,
text=black,
rotate=0.0
]{8};
\end{axis}

\end{tikzpicture}

%% file: plots/horeka/gmres_poly_GS2_1.tex
\begin{tikzpicture}

\definecolor{burlywood253190133}{RGB}{253,190,133}
\definecolor{coral25314160}{RGB}{253,141,60}
\definecolor{cornflowerblue107174214}{RGB}{107,174,214}
\definecolor{darkgray176}{RGB}{176,176,176}
\definecolor{lightgray204}{RGB}{204,204,204}
\definecolor{orangered217711}{RGB}{217,71,1}
\definecolor{powderblue189215231}{RGB}{189,215,231}
\definecolor{steelblue33113181}{RGB}{33,113,181}

\begin{axis}[
	legend cell align={left},
legend style={fill opacity=0.8, draw opacity=1, text opacity=1, draw=lightgray204},
tick align=outside,
tick pos=left,
x grid style={darkgray176},
xtick={0, 1, 2, 3, 4, 5, 6, 7},
xticklabels={1, 2, 3, 4, 5, 6, 7, 8},
width  = 0.49\linewidth,
height = 0.2\textheight,
xmin=-0.68, xmax=7.68,
xtick style={color=black},
y grid style={darkgray176},
ylabel={\small Normalized time},
ymin=0, ymax=1.05,
ytick style={color=black}
]
\draw[draw=none,fill=powderblue189215231] (axis cs:0,0) rectangle (axis cs:-0.3,0.088054288092441);
\draw[draw=none,fill=powderblue189215231] (axis cs:1,0) rectangle (axis cs:0.7,0.0199479439807887);
\draw[draw=none,fill=powderblue189215231] (axis cs:2,0) rectangle (axis cs:1.7,0.0835137872263346);
\draw[draw=none,fill=powderblue189215231] (axis cs:3,0) rectangle (axis cs:2.7,0.0117122428216393);
\draw[draw=none,fill=powderblue189215231] (axis cs:4,0) rectangle (axis cs:3.7,0.0106974116650482);
\draw[draw=none,fill=powderblue189215231] (axis cs:5,0) rectangle (axis cs:4.7,0.0110815524479049);
\draw[draw=none,fill=powderblue189215231] (axis cs:6,0) rectangle (axis cs:5.7,0.00853351485017232);
\draw[draw=none,fill=powderblue189215231] (axis cs:7,0) rectangle (axis cs:6.7,0.00736106552790031);
\draw[draw=none,fill=cornflowerblue107174214] (axis cs:0,0.088054288092441) rectangle (axis cs:-0.3,0.105932379563377);
\draw[draw=none,fill=cornflowerblue107174214] (axis cs:1,0.0199479439807887) rectangle (axis cs:0.7,0.0363459400649538);
\draw[draw=none,fill=cornflowerblue107174214] (axis cs:2,0.0835137872263346) rectangle (axis cs:1.7,0.0925776500649086);
\draw[draw=none,fill=cornflowerblue107174214] (axis cs:3,0.0117122428216393) rectangle (axis cs:2.7,0.0165799790319716);
\draw[draw=none,fill=cornflowerblue107174214] (axis cs:4,0.0106974116650482) rectangle (axis cs:3.7,0.0159432372581223);
\draw[draw=none,fill=cornflowerblue107174214] (axis cs:5,0.0110815524479049) rectangle (axis cs:4.7,0.0181525630720984);
\draw[draw=none,fill=cornflowerblue107174214] (axis cs:6,0.00853351485017232) rectangle (axis cs:5.7,0.0121583600478652);
\draw[draw=none,fill=cornflowerblue107174214] (axis cs:7,0.00736106552790031) rectangle (axis cs:6.7,0.0104564639906858);
\draw[draw=none,fill=steelblue33113181] (axis cs:0,0.105932379563377) rectangle (axis cs:-0.3,1);
\draw[draw=none,fill=steelblue33113181] (axis cs:1,0.0363459400649538) rectangle (axis cs:0.7,1);
\draw[draw=none,fill=steelblue33113181] (axis cs:2,0.0925776500649086) rectangle (axis cs:1.7,1);
\draw[draw=none,fill=steelblue33113181] (axis cs:3,0.0165799790319716) rectangle (axis cs:2.7,1);
\draw[draw=none,fill=steelblue33113181] (axis cs:4,0.0159432372581223) rectangle (axis cs:3.7,1);
\draw[draw=none,fill=steelblue33113181] (axis cs:5,0.0181525630720984) rectangle (axis cs:4.7,1);
\draw[draw=none,fill=steelblue33113181] (axis cs:6,0.0121583600478652) rectangle (axis cs:5.7,1);
\draw[draw=none,fill=steelblue33113181] (axis cs:7,0.0104564639906858) rectangle (axis cs:6.7,1);
\draw[draw=none,fill=burlywood253190133] (axis cs:0,0) rectangle (axis cs:0.3,0.0963419615029989);
\draw[draw=none,fill=burlywood253190133] (axis cs:1,0) rectangle (axis cs:1.3,0.0242813471147331);
\draw[draw=none,fill=burlywood253190133] (axis cs:2,0) rectangle (axis cs:2.3,0.100138187886123);
\draw[draw=none,fill=burlywood253190133] (axis cs:3,0) rectangle (axis cs:3.3,0.014441393681451);
\draw[draw=none,fill=burlywood253190133] (axis cs:4,0) rectangle (axis cs:4.3,0.0125524558317971);
\draw[draw=none,fill=burlywood253190133] (axis cs:5,0) rectangle (axis cs:5.3,0.0143425671454519);
\draw[draw=none,fill=burlywood253190133] (axis cs:6,0) rectangle (axis cs:6.3,0.0103874055816428);
\draw[draw=none,fill=burlywood253190133] (axis cs:7,0) rectangle (axis cs:7.3,0.00944409164419791);
\draw[draw=none,fill=coral25314160] (axis cs:0,0.0963419615029989) rectangle (axis cs:0.3,0.114586908476074);
\draw[draw=none,fill=coral25314160] (axis cs:1,0.0242813471147331) rectangle (axis cs:1.3,0.0404096122236313);
\draw[draw=none,fill=coral25314160] (axis cs:2,0.100138187886123) rectangle (axis cs:2.3,0.108810952933782);
\draw[draw=none,fill=coral25314160] (axis cs:3,0.014441393681451) rectangle (axis cs:3.3,0.0193150355468761);
\draw[draw=none,fill=coral25314160] (axis cs:4,0.0125524558317971) rectangle (axis cs:4.3,0.0178137958813387);
\draw[draw=none,fill=coral25314160] (axis cs:5,0.0143425671454519) rectangle (axis cs:5.3,0.0212292801523127);
\draw[draw=none,fill=coral25314160] (axis cs:6,0.0103874055816428) rectangle (axis cs:6.3,0.0141227521017352);
\draw[draw=none,fill=coral25314160] (axis cs:7,0.00944409164419791) rectangle (axis cs:7.3,0.012633280088048);
\draw[draw=none,fill=orangered217711] (axis cs:0,0.114586908476074) rectangle (axis cs:0.3,0.559142405968821);
\draw[draw=none,fill=orangered217711] (axis cs:1,0.0404096122236313) rectangle (axis cs:1.3,0.820337756715426);
\draw[draw=none,fill=orangered217711] (axis cs:2,0.108810952933782) rectangle (axis cs:2.3,0.64234381475413);
\draw[draw=none,fill=orangered217711] (axis cs:3,0.0193150355468762) rectangle (axis cs:3.3,0.776173457357996);
\draw[draw=none,fill=orangered217711] (axis cs:4,0.0178137958813387) rectangle (axis cs:4.3,0.740170767622338);
\draw[draw=none,fill=orangered217711] (axis cs:5,0.0212292801523127) rectangle (axis cs:5.3,0.554676852933886);
\draw[draw=none,fill=orangered217711] (axis cs:6,0.0141227521017352) rectangle (axis cs:6.3,0.721624546241388);
\draw[draw=none,fill=orangered217711] (axis cs:7,0.012633280088048) rectangle (axis cs:7.3,0.755992924747098);
\draw (axis cs:0.15,0.609142405968821) node[
scale=0.72,
anchor=base,
text=black,
rotate=0.0
]{2};
\draw (axis cs:1.15,0.870337756715426) node[
scale=0.72,
anchor=base,
text=black,
rotate=0.0
]{2};
\draw (axis cs:2.15,0.69234381475413) node[
scale=0.72,
anchor=base,
text=black,
rotate=0.0
]{2};
\draw (axis cs:3.15,0.826173457357996) node[
scale=0.72,
anchor=base,
text=black,
rotate=0.0
]{2};
\draw (axis cs:4.15,0.790170767622338) node[
scale=0.72,
anchor=base,
text=black,
rotate=0.0
]{3};
\draw (axis cs:5.15,0.604676852933887) node[
scale=0.72,
anchor=base,
text=black,
rotate=0.0
]{3};
\draw (axis cs:6.15,0.771624546241388) node[
scale=0.72,
anchor=base,
text=black,
rotate=0.0
]{1};
\draw (axis cs:7.15,0.805992924747098) node[
scale=0.72,
anchor=base,
text=black,
rotate=0.0
]{4};
\end{axis}

\end{tikzpicture}

%% file: plots/tg097/gmres_poly_GS2_1.tex
\begin{tikzpicture}

\definecolor{burlywood253190133}{RGB}{253,190,133}
\definecolor{coral25314160}{RGB}{253,141,60}
\definecolor{cornflowerblue107174214}{RGB}{107,174,214}
\definecolor{darkgray176}{RGB}{176,176,176}
\definecolor{lightgray204}{RGB}{204,204,204}
\definecolor{orangered217711}{RGB}{217,71,1}
\definecolor{powderblue189215231}{RGB}{189,215,231}
\definecolor{steelblue33113181}{RGB}{33,113,181}

\begin{axis}[
		legend cell align={left},
legend style={fill opacity=0.8, draw opacity=1, text opacity=1, draw=lightgray204},
tick align=outside,
tick pos=left,
x grid style={darkgray176},
xtick={0, 1, 2, 3, 4, 5, 6, 7},
xticklabels={1, 2, 3, 4, 5, 6, 7, 8},
width  = 0.49\linewidth,
height = 0.2\textheight,
xmin=-0.68, xmax=7.68,
xtick style={color=black},
y grid style={darkgray176},
ylabel={\small Normalized time},
ymin=0, ymax=1.05,
ytick style={color=black},
]
\draw[draw=none,fill=powderblue189215231] (axis cs:0,0) rectangle (axis cs:-0.3,0.126609776734236);
\draw[draw=none,fill=powderblue189215231] (axis cs:1,0) rectangle (axis cs:0.7,0.0223196756431066);
\draw[draw=none,fill=powderblue189215231] (axis cs:2,0) rectangle (axis cs:1.7,0.0766312425801975);
\draw[draw=none,fill=powderblue189215231] (axis cs:3,0) rectangle (axis cs:2.7,0.0102913557754016);
\draw[draw=none,fill=powderblue189215231] (axis cs:4,0) rectangle (axis cs:3.7,0.00958645935109688);
\draw[draw=none,fill=powderblue189215231] (axis cs:5,0) rectangle (axis cs:4.7,0.0106340347420368);
\draw[draw=none,fill=powderblue189215231] (axis cs:6,0) rectangle (axis cs:5.7,0.00808821466598975);
\draw[draw=none,fill=powderblue189215231] (axis cs:7,0) rectangle (axis cs:6.7,0.00718505175327792);
\draw[draw=none,fill=cornflowerblue107174214] (axis cs:0,0.126609776734236) rectangle (axis cs:-0.3,0.155836580414849);
\draw[draw=none,fill=cornflowerblue107174214] (axis cs:1,0.0223196756431066) rectangle (axis cs:0.7,0.039286086953859);
\draw[draw=none,fill=cornflowerblue107174214] (axis cs:2,0.0766312425801975) rectangle (axis cs:1.7,0.0844417284917232);
\draw[draw=none,fill=cornflowerblue107174214] (axis cs:3,0.0102913557754016) rectangle (axis cs:2.7,0.013100068899399);
\draw[draw=none,fill=cornflowerblue107174214] (axis cs:4,0.00958645935109688) rectangle (axis cs:3.7,0.0132230869587621);
\draw[draw=none,fill=cornflowerblue107174214] (axis cs:5,0.0106340347420368) rectangle (axis cs:4.7,0.0162483023114969);
\draw[draw=none,fill=cornflowerblue107174214] (axis cs:6,0.00808821466598975) rectangle (axis cs:5.7,0.011070960636269);
\draw[draw=none,fill=cornflowerblue107174214] (axis cs:7,0.00718505175327792) rectangle (axis cs:6.7,0.00974646467014213);
\draw[draw=none,fill=steelblue33113181] (axis cs:0,0.155836580414849) rectangle (axis cs:-0.3,1);
\draw[draw=none,fill=steelblue33113181] (axis cs:1,0.039286086953859) rectangle (axis cs:0.7,1);
\draw[draw=none,fill=steelblue33113181] (axis cs:2,0.0844417284917232) rectangle (axis cs:1.7,1);
\draw[draw=none,fill=steelblue33113181] (axis cs:3,0.013100068899399) rectangle (axis cs:2.7,1);
\draw[draw=none,fill=steelblue33113181] (axis cs:4,0.0132230869587621) rectangle (axis cs:3.7,1);
\draw[draw=none,fill=steelblue33113181] (axis cs:5,0.0162483023114969) rectangle (axis cs:4.7,1);
\draw[draw=none,fill=steelblue33113181] (axis cs:6,0.011070960636269) rectangle (axis cs:5.7,1);
\draw[draw=none,fill=steelblue33113181] (axis cs:7,0.00974646467014211) rectangle (axis cs:6.7,1);
\draw[draw=none,fill=burlywood253190133] (axis cs:0,0) rectangle (axis cs:0.3,0.150773466483119);
\draw[draw=none,fill=burlywood253190133] (axis cs:1,0) rectangle (axis cs:1.3,0.0252060619294033);
\draw[draw=none,fill=burlywood253190133] (axis cs:2,0) rectangle (axis cs:2.3,0.0941262217808799);
\draw[draw=none,fill=burlywood253190133] (axis cs:3,0) rectangle (axis cs:3.3,0.0134256749701026);
\draw[draw=none,fill=burlywood253190133] (axis cs:4,0) rectangle (axis cs:4.3,0.0124705692546704);
\draw[draw=none,fill=burlywood253190133] (axis cs:5,0) rectangle (axis cs:5.3,0.0130147867362189);
\draw[draw=none,fill=burlywood253190133] (axis cs:6,0) rectangle (axis cs:6.3,0.0106316406050266);
\draw[draw=none,fill=burlywood253190133] (axis cs:7,0) rectangle (axis cs:7.3,0.00874375537875779);
\draw[draw=none,fill=coral25314160] (axis cs:0,0.150773466483119) rectangle (axis cs:0.3,0.179275740156523);
\draw[draw=none,fill=coral25314160] (axis cs:1,0.0252060619294033) rectangle (axis cs:1.3,0.0436281503455348);
\draw[draw=none,fill=coral25314160] (axis cs:2,0.0941262217808799) rectangle (axis cs:2.3,0.101662990582222);
\draw[draw=none,fill=coral25314160] (axis cs:3,0.0134256749701026) rectangle (axis cs:3.3,0.016118943556542);
\draw[draw=none,fill=coral25314160] (axis cs:4,0.0124705692546704) rectangle (axis cs:4.3,0.0161115857443145);
\draw[draw=none,fill=coral25314160] (axis cs:5,0.0130147867362189) rectangle (axis cs:5.3,0.0186180526790886);
\draw[draw=none,fill=coral25314160] (axis cs:6,0.0106316406050266) rectangle (axis cs:6.3,0.0136302385190362);
\draw[draw=none,fill=coral25314160] (axis cs:7,0.00874375537875779) rectangle (axis cs:7.3,0.0113081581554526);
\draw[draw=none,fill=orangered217711] (axis cs:0,0.179275740156523) rectangle (axis cs:0.3,0.820939162811717);
\draw[draw=none,fill=orangered217711] (axis cs:1,0.0436281503455348) rectangle (axis cs:1.3,0.535932023176829);
\draw[draw=none,fill=orangered217711] (axis cs:2,0.101662990582222) rectangle (axis cs:2.3,0.519202528382236);
\draw[draw=none,fill=orangered217711] (axis cs:3,0.016118943556542) rectangle (axis cs:3.3,0.49630918565948);
\draw[draw=none,fill=orangered217711] (axis cs:4,0.0161115857443145) rectangle (axis cs:4.3,0.483681127359825);
\draw[draw=none,fill=orangered217711] (axis cs:5,0.0186180526790886) rectangle (axis cs:5.3,0.362171347033659);
\draw[draw=none,fill=orangered217711] (axis cs:6,0.0136302385190362) rectangle (axis cs:6.3,0.374645354232143);
\draw[draw=none,fill=orangered217711] (axis cs:7,0.0113081581554526) rectangle (axis cs:7.3,0.4152083226695);
\draw (axis cs:0.15,0.870939162811717) node[
scale=0.72,
anchor=base,
text=black,
rotate=0.0
]{3};
\draw (axis cs:1.15,0.585932023176829) node[
scale=0.72,
anchor=base,
text=black,
rotate=0.0
]{2};
\draw (axis cs:2.15,0.569202528382236) node[
scale=0.72,
anchor=base,
text=black,
rotate=0.0
]{4};
\draw (axis cs:3.15,0.54630918565948) node[
scale=0.72,
anchor=base,
text=black,
rotate=0.0
]{4};
\draw (axis cs:4.15,0.533681127359825) node[
scale=0.72,
anchor=base,
text=black,
rotate=0.0
]{4};
\draw (axis cs:5.15,0.412171347033659) node[
scale=0.72,
anchor=base,
text=black,
rotate=0.0
]{4};
\draw (axis cs:6.15,0.424645354232143) node[
scale=0.72,
anchor=base,
text=black,
rotate=0.0
]{4};
\draw (axis cs:7.15,0.4652083226695) node[
scale=0.72,
anchor=base,
text=black,
rotate=0.0
]{4};
\end{axis}

\end{tikzpicture}

%% file: plots/gmres_mg_legend.tex
\begin{tikzpicture}
\definecolor{burlywood253190133}{RGB}{253,190,133}
\definecolor{coral25314160}{RGB}{253,141,60}
\definecolor{cornflowerblue107174214}{RGB}{107,174,214}
\definecolor{darkgray176}{RGB}{176,176,176}
\definecolor{lightgray204}{RGB}{204,204,204}
\definecolor{orangered2308513}{RGB}{230,85,13}
\definecolor{powderblue189215231}{RGB}{189,215,231}
\definecolor{saddlebrown166543}{RGB}{166,54,3}
\definecolor{steelblue49130189}{RGB}{49,130,189}
\definecolor{teal881156}{RGB}{8,81,156}

  \centering
  \begin{customlegend}[
      legend columns=4,
      legend style={
      anchor=north,
      /tikz/every even column/.append style={column sep=0.5cm}},
      legend entries={Baseline: AMG, Baseline: SpMV, Baseline: Ortho, Baseline: Misc, RACE: AMG, RACE: SpMV, RACE: Ortho, RACE: Misc},
      width  = \linewidth,
      legend image post style={xscale=0.5},
      legend image code/.code={\draw[#1, draw=none] (0cm,-0.15cm) rectangle (0.26cm,0.15cm);}, 
	  legend style={draw=none} 
      ]
    \addlegendimage{fill=teal881156, color=teal881156} 
    \addlegendimage{fill=steelblue49130189, color=steelblue49130189}
    \addlegendimage{fill=cornflowerblue107174214, color=cornflowerblue107174214}
    \addlegendimage{fill=powderblue189215231, color=powderblue189215231}
	
	\addlegendimage{fill=saddlebrown166543, color=saddlebrown166543}
	\addlegendimage{fill=orangered2308513, color=orangered2308513}
	\addlegendimage{fill=coral25314160, color=coral25314160}
	\addlegendimage{fill=burlywood253190133, color=burlywood253190133}
  \end{customlegend}
\end{tikzpicture}

%% file: plots/horeka/gmres_amg_gs2_1_2.tex
\begin{tikzpicture}

\definecolor{burlywood253190133}{RGB}{253,190,133}
\definecolor{coral25314160}{RGB}{253,141,60}
\definecolor{cornflowerblue107174214}{RGB}{107,174,214}
\definecolor{darkgray176}{RGB}{176,176,176}
\definecolor{lightgray204}{RGB}{204,204,204}
\definecolor{orangered2308513}{RGB}{230,85,13}
\definecolor{powderblue189215231}{RGB}{189,215,231}
\definecolor{saddlebrown166543}{RGB}{166,54,3}
\definecolor{steelblue49130189}{RGB}{49,130,189}
\definecolor{teal881156}{RGB}{8,81,156}

\begin{axis}[
	legend cell align={left},
legend style={fill opacity=0.8, draw opacity=1, text opacity=1, draw=lightgray204},
tick align=outside,
tick pos=left,
x grid style={darkgray176},
xtick={0, 1, 2, 3, 4, 5, 6, 7},
xticklabels={1, 2, 3, 4, 5, 6, 7, 8},
width  = 0.49\linewidth,
height = 0.2\textheight,
xmin=-0.68, xmax=7.68,
xtick style={color=black},
y grid style={darkgray176},
ylabel={\small Normalized time},
ymin=0, ymax=1.25,
ytick style={color=black},
]
\draw[draw=none,fill=powderblue189215231] (axis cs:0,0) rectangle (axis cs:-0.3,0.25680893785718);
\draw[draw=none,fill=powderblue189215231] (axis cs:1,0) rectangle (axis cs:0.7,0.360542005420054);
\draw[draw=none,fill=powderblue189215231] (axis cs:2,0) rectangle (axis cs:1.7,0.108980456950216);
\draw[draw=none,fill=powderblue189215231] (axis cs:3,0) rectangle (axis cs:2.7,0.0248289840385104);
\draw[draw=none,fill=powderblue189215231] (axis cs:4,0) rectangle (axis cs:3.7,0.0387769434091639);
\draw[draw=none,fill=powderblue189215231] (axis cs:5,0) rectangle (axis cs:4.7,0.0230804887221336);
\draw[draw=none,fill=powderblue189215231] (axis cs:6,0) rectangle (axis cs:5.7,0.0160257364248793);
\draw[draw=none,fill=powderblue189215231] (axis cs:7,0) rectangle (axis cs:6.7,0.0216599469452987);
\draw[draw=none,fill=cornflowerblue107174214] (axis cs:0,0.25680893785718) rectangle (axis cs:-0.3,0.373680687844308);
\draw[draw=none,fill=cornflowerblue107174214] (axis cs:1,0.360542005420054) rectangle (axis cs:0.7,0.430135501355014);
\draw[draw=none,fill=cornflowerblue107174214] (axis cs:2,0.108980456950216) rectangle (axis cs:1.7,0.231195447187036);
\draw[draw=none,fill=cornflowerblue107174214] (axis cs:3,0.0248289840385104) rectangle (axis cs:2.7,0.0964812754504432);
\draw[draw=none,fill=cornflowerblue107174214] (axis cs:4,0.0387769434091639) rectangle (axis cs:3.7,0.114016655422966);
\draw[draw=none,fill=cornflowerblue107174214] (axis cs:5,0.0230804887221336) rectangle (axis cs:4.7,0.126621275330377);
\draw[draw=none,fill=cornflowerblue107174214] (axis cs:6,0.0160257364248793) rectangle (axis cs:5.7,0.0741695687114922);
\draw[draw=none,fill=cornflowerblue107174214] (axis cs:7,0.0216599469452987) rectangle (axis cs:6.7,0.0756429318922729);
\draw[draw=none,fill=steelblue49130189] (axis cs:0,0.373680687844308) rectangle (axis cs:-0.3,0.404829326056737);
\draw[draw=none,fill=steelblue49130189] (axis cs:1,0.430135501355014) rectangle (axis cs:0.7,0.463739837398374);
\draw[draw=none,fill=steelblue49130189] (axis cs:2,0.231195447187036) rectangle (axis cs:1.7,0.293045611575627);
\draw[draw=none,fill=steelblue49130189] (axis cs:3,0.0964812754504432) rectangle (axis cs:2.7,0.182596945763247);
\draw[draw=none,fill=steelblue49130189] (axis cs:4,0.114016655422966) rectangle (axis cs:3.7,0.208585394985678);
\draw[draw=none,fill=steelblue49130189] (axis cs:5,0.126621275330377) rectangle (axis cs:4.7,0.226805814289088);
\draw[draw=none,fill=steelblue49130189] (axis cs:6,0.0741695687114922) rectangle (axis cs:5.7,0.188848706729768);
\draw[draw=none,fill=steelblue49130189] (axis cs:7,0.0756429318922729) rectangle (axis cs:6.7,0.19143962738865);
\draw[draw=none,fill=teal881156] (axis cs:0,0.404829326056737) rectangle (axis cs:-0.3,1);
\draw[draw=none,fill=teal881156] (axis cs:1,0.463739837398374) rectangle (axis cs:0.7,1);
\draw[draw=none,fill=teal881156] (axis cs:2,0.293045611575627) rectangle (axis cs:1.7,1);
\draw[draw=none,fill=teal881156] (axis cs:3,0.182596945763247) rectangle (axis cs:2.7,1);
\draw[draw=none,fill=teal881156] (axis cs:4,0.208585394985678) rectangle (axis cs:3.7,1);
\draw[draw=none,fill=teal881156] (axis cs:5,0.226805814289088) rectangle (axis cs:4.7,1);
\draw[draw=none,fill=teal881156] (axis cs:6,0.188848706729768) rectangle (axis cs:5.7,1);
\draw[draw=none,fill=teal881156] (axis cs:7,0.19143962738865) rectangle (axis cs:6.7,1);
\draw[draw=none,fill=burlywood253190133] (axis cs:0,0) rectangle (axis cs:0.3,0.356536065489368);
\draw[draw=none,fill=burlywood253190133] (axis cs:1,0) rectangle (axis cs:1.3,0.584932249322493);
\draw[draw=none,fill=burlywood253190133] (axis cs:2,0) rectangle (axis cs:2.3,0.137335402793772);
\draw[draw=none,fill=burlywood253190133] (axis cs:3,0) rectangle (axis cs:3.3,0.0316182167340448);
\draw[draw=none,fill=burlywood253190133] (axis cs:4,0) rectangle (axis cs:4.3,0.0279664016102453);
\draw[draw=none,fill=burlywood253190133] (axis cs:5,0) rectangle (axis cs:5.3,0.0326507141045488);
\draw[draw=none,fill=burlywood253190133] (axis cs:6,0) rectangle (axis cs:6.3,0.0236132194342749);
\draw[draw=none,fill=burlywood253190133] (axis cs:7,0) rectangle (axis cs:7.3,0.025302387461568);
\draw[draw=none,fill=coral25314160] (axis cs:0,0.356536065489368) rectangle (axis cs:0.3,0.459506770323843);
\draw[draw=none,fill=coral25314160] (axis cs:1,0.584932249322493) rectangle (axis cs:1.3,0.648021680216802);
\draw[draw=none,fill=coral25314160] (axis cs:2,0.137335402793772) rectangle (axis cs:2.3,0.243929304561158);
\draw[draw=none,fill=coral25314160] (axis cs:3,0.0316182167340448) rectangle (axis cs:3.3,0.104474871760978);
\draw[draw=none,fill=coral25314160] (axis cs:4,0.0279664016102453) rectangle (axis cs:4.3,0.102528726733834);
\draw[draw=none,fill=coral25314160] (axis cs:5,0.0326507141045488) rectangle (axis cs:5.3,0.13096143922896);
\draw[draw=none,fill=coral25314160] (axis cs:6,0.0236132194342749) rectangle (axis cs:6.3,0.0790323540792333);
\draw[draw=none,fill=coral25314160] (axis cs:7,0.025302387461568) rectangle (axis cs:7.3,0.0789411774249814);
\draw[draw=none,fill=orangered2308513] (axis cs:0,0.459506770323843) rectangle (axis cs:0.3,0.483859342017196);
\draw[draw=none,fill=orangered2308513] (axis cs:1,0.648021680216802) rectangle (axis cs:1.3,0.683252032520325);
\draw[draw=none,fill=orangered2308513] (axis cs:2,0.243929304561158) rectangle (axis cs:2.3,0.298820073766251);
\draw[draw=none,fill=orangered2308513] (axis cs:3,0.104474871760978) rectangle (axis cs:3.3,0.195415340214362);
\draw[draw=none,fill=orangered2308513] (axis cs:4,0.102528726733834) rectangle (axis cs:4.3,0.200951659459639);
\draw[draw=none,fill=orangered2308513] (axis cs:5,0.13096143922896) rectangle (axis cs:5.3,0.232505969091911);
\draw[draw=none,fill=orangered2308513] (axis cs:6,0.0790323540792333) rectangle (axis cs:6.3,0.192507267513965);
\draw[draw=none,fill=orangered2308513] (axis cs:7,0.0789411774249813) rectangle (axis cs:7.3,0.193030804871576);
\draw[draw=none,fill=saddlebrown166543] (axis cs:0,0.483859342017196) rectangle (axis cs:0.3,0.975235545487309);
\draw[draw=none,fill=saddlebrown166543] (axis cs:1,0.683252032520325) rectangle (axis cs:1.3,1.19208672086721);
\draw[draw=none,fill=saddlebrown166543] (axis cs:2,0.298820073766251) rectangle (axis cs:2.3,0.76156144127906);
\draw[draw=none,fill=saddlebrown166543] (axis cs:3,0.195415340214362) rectangle (axis cs:3.3,0.93033050234447);
\draw[draw=none,fill=saddlebrown166543] (axis cs:4,0.200951659459639) rectangle (axis cs:4.3,0.82918790988819);
\draw[draw=none,fill=saddlebrown166543] (axis cs:5,0.232505969091911) rectangle (axis cs:5.3,0.709178408067239);
\draw[draw=none,fill=saddlebrown166543] (axis cs:6,0.192507267513965) rectangle (axis cs:6.3,0.7946356817573);
\draw[draw=none,fill=saddlebrown166543] (axis cs:7,0.193030804871576) rectangle (axis cs:7.3,0.825814160204802);
\end{axis}

\end{tikzpicture}

%% file: plots/tg097/gmres_amg_gs2_1_2.tex
\begin{tikzpicture}

\definecolor{burlywood253190133}{RGB}{253,190,133}
\definecolor{coral25314160}{RGB}{253,141,60}
\definecolor{cornflowerblue107174214}{RGB}{107,174,214}
\definecolor{darkgray176}{RGB}{176,176,176}
\definecolor{lightgray204}{RGB}{204,204,204}
\definecolor{orangered2308513}{RGB}{230,85,13}
\definecolor{powderblue189215231}{RGB}{189,215,231}
\definecolor{saddlebrown166543}{RGB}{166,54,3}
\definecolor{steelblue49130189}{RGB}{49,130,189}
\definecolor{teal881156}{RGB}{8,81,156}

\begin{axis}[
	legend cell align={left},
legend style={fill opacity=0.8, draw opacity=1, text opacity=1, draw=lightgray204},
tick align=outside,
tick pos=left,
x grid style={darkgray176},
xtick={0, 1, 2, 3, 4, 5, 6, 7},
xticklabels={1, 2, 3, 4, 5, 6, 7, 8},
width  = 0.49\linewidth,
height = 0.2\textheight,
xmin=-0.68, xmax=7.68,
xtick style={color=black},
y grid style={darkgray176},
ylabel={\small Normalized time},
ymin=0, ymax=1.25,
ytick style={color=black},
]
\draw[draw=none,fill=powderblue189215231] (axis cs:0,0) rectangle (axis cs:-0.3,0.263030022560305);
\draw[draw=none,fill=powderblue189215231] (axis cs:1,0) rectangle (axis cs:0.7,0.384076990376203);
\draw[draw=none,fill=powderblue189215231] (axis cs:2,0) rectangle (axis cs:1.7,0.0831852233277411);
\draw[draw=none,fill=powderblue189215231] (axis cs:3,0) rectangle (axis cs:2.7,0.0188838845152147);
\draw[draw=none,fill=powderblue189215231] (axis cs:4,0) rectangle (axis cs:3.7,0.0172514508003316);
\draw[draw=none,fill=powderblue189215231] (axis cs:5,0) rectangle (axis cs:4.7,0.0182273320355001);
\draw[draw=none,fill=powderblue189215231] (axis cs:6,0) rectangle (axis cs:5.7,0.0141864981018969);
\draw[draw=none,fill=powderblue189215231] (axis cs:7,0) rectangle (axis cs:6.7,0.0133002315023929);
\draw[draw=none,fill=cornflowerblue107174214] (axis cs:0,0.263030022560305) rectangle (axis cs:-0.3,0.379244519002719);
\draw[draw=none,fill=cornflowerblue107174214] (axis cs:1,0.384076990376203) rectangle (axis cs:0.7,0.433070866141732);
\draw[draw=none,fill=cornflowerblue107174214] (axis cs:2,0.0831852233277411) rectangle (axis cs:1.7,0.19972426623258);
\draw[draw=none,fill=cornflowerblue107174214] (axis cs:3,0.0188838845152147) rectangle (axis cs:2.7,0.0694993730105142);
\draw[draw=none,fill=cornflowerblue107174214] (axis cs:4,0.0172514508003316) rectangle (axis cs:3.7,0.0731892098718193);
\draw[draw=none,fill=cornflowerblue107174214] (axis cs:5,0.0182273320355001) rectangle (axis cs:4.7,0.104668031198022);
\draw[draw=none,fill=cornflowerblue107174214] (axis cs:6,0.0141864981018969) rectangle (axis cs:5.7,0.0616093348183826);
\draw[draw=none,fill=cornflowerblue107174214] (axis cs:7,0.0133002315023929) rectangle (axis cs:6.7,0.0597802958485211);
\draw[draw=none,fill=steelblue49130189] (axis cs:0,0.379244519002719) rectangle (axis cs:-0.3,0.399375253080349);
\draw[draw=none,fill=steelblue49130189] (axis cs:1,0.433070866141732) rectangle (axis cs:0.7,0.459645669291339);
\draw[draw=none,fill=steelblue49130189] (axis cs:2,0.19972426623258) rectangle (axis cs:1.7,0.269571283823065);
\draw[draw=none,fill=steelblue49130189] (axis cs:3,0.0694993730105141) rectangle (axis cs:2.7,0.171527999347041);
\draw[draw=none,fill=steelblue49130189] (axis cs:4,0.0731892098718193) rectangle (axis cs:3.7,0.176807601556023);
\draw[draw=none,fill=steelblue49130189] (axis cs:5,0.104668031198022) rectangle (axis cs:4.7,0.21673537957455);
\draw[draw=none,fill=steelblue49130189] (axis cs:6,0.0616093348183827) rectangle (axis cs:5.7,0.18256077720177);
\draw[draw=none,fill=steelblue49130189] (axis cs:7,0.0597802958485211) rectangle (axis cs:6.7,0.180551186844308);
\draw[draw=none,fill=teal881156] (axis cs:0,0.399375253080349) rectangle (axis cs:-0.3,1);
\draw[draw=none,fill=teal881156] (axis cs:1,0.459645669291339) rectangle (axis cs:0.7,1);
\draw[draw=none,fill=teal881156] (axis cs:2,0.269571283823065) rectangle (axis cs:1.7,1);
\draw[draw=none,fill=teal881156] (axis cs:3,0.171527999347041) rectangle (axis cs:2.7,1);
\draw[draw=none,fill=teal881156] (axis cs:4,0.176807601556023) rectangle (axis cs:3.7,1);
\draw[draw=none,fill=teal881156] (axis cs:5,0.21673537957455) rectangle (axis cs:4.7,1);
\draw[draw=none,fill=teal881156] (axis cs:6,0.18256077720177) rectangle (axis cs:5.7,1);
\draw[draw=none,fill=teal881156] (axis cs:7,0.180551186844308) rectangle (axis cs:6.7,1);
\draw[draw=none,fill=burlywood253190133] (axis cs:0,0) rectangle (axis cs:0.3,0.327934285879563);
\draw[draw=none,fill=burlywood253190133] (axis cs:1,0) rectangle (axis cs:1.3,0.567366579177603);
\draw[draw=none,fill=burlywood253190133] (axis cs:2,0) rectangle (axis cs:2.3,0.112984402770626);
\draw[draw=none,fill=burlywood253190133] (axis cs:3,0) rectangle (axis cs:3.3,0.0259662687076595);
\draw[draw=none,fill=burlywood253190133] (axis cs:4,0) rectangle (axis cs:4.3,0.02327147503348);
\draw[draw=none,fill=burlywood253190133] (axis cs:5,0) rectangle (axis cs:5.3,0.0247116350615081);
\draw[draw=none,fill=burlywood253190133] (axis cs:6,0) rectangle (axis cs:6.3,0.0210173365911119);
\draw[draw=none,fill=burlywood253190133] (axis cs:7,0) rectangle (axis cs:7.3,0.0199427284630397);
\draw[draw=none,fill=coral25314160] (axis cs:0,0.327934285879563) rectangle (axis cs:0.3,0.442413374211836);
\draw[draw=none,fill=coral25314160] (axis cs:1,0.567366579177603) rectangle (axis cs:1.3,0.618219597550306);
\draw[draw=none,fill=coral25314160] (axis cs:2,0.112984402770626) rectangle (axis cs:2.3,0.219341228842417);
\draw[draw=none,fill=coral25314160] (axis cs:3,0.0259662687076595) rectangle (axis cs:3.3,0.0740404092868642);
\draw[draw=none,fill=coral25314160] (axis cs:4,0.02327147503348) rectangle (axis cs:4.3,0.0777884063516357);
\draw[draw=none,fill=coral25314160] (axis cs:5,0.0247116350615081) rectangle (axis cs:5.3,0.110103923930851);
\draw[draw=none,fill=coral25314160] (axis cs:6,0.0210173365911119) rectangle (axis cs:6.3,0.0684436721150875);
\draw[draw=none,fill=coral25314160] (axis cs:7,0.0199427284630397) rectangle (axis cs:7.3,0.0666437377351078);
\draw[draw=none,fill=orangered2308513] (axis cs:0,0.442413374211836) rectangle (axis cs:0.3,0.473130097761324);
\draw[draw=none,fill=orangered2308513] (axis cs:1,0.618219597550306) rectangle (axis cs:1.3,0.654855643044619);
\draw[draw=none,fill=orangered2308513] (axis cs:2,0.219341228842417) rectangle (axis cs:2.3,0.285234955068684);
\draw[draw=none,fill=orangered2308513] (axis cs:3,0.0740404092868642) rectangle (axis cs:3.3,0.178209703868043);
\draw[draw=none,fill=orangered2308513] (axis cs:4,0.0777884063516357) rectangle (axis cs:4.3,0.180024233148396);
\draw[draw=none,fill=orangered2308513] (axis cs:5,0.110103923930851) rectangle (axis cs:5.3,0.218388241130199);
\draw[draw=none,fill=orangered2308513] (axis cs:6,0.0684436721150875) rectangle (axis cs:6.3,0.186790835456915);
\draw[draw=none,fill=orangered2308513] (axis cs:7,0.0666437377351077) rectangle (axis cs:7.3,0.184949405788757);
\draw[draw=none,fill=saddlebrown166543] (axis cs:0,0.473130097761324) rectangle (axis cs:0.3,0.982067449528548);
\draw[draw=none,fill=saddlebrown166543] (axis cs:1,0.654855643044619) rectangle (axis cs:1.3,1.05675853018373);
\draw[draw=none,fill=saddlebrown166543] (axis cs:2,0.285234955068684) rectangle (axis cs:2.3,0.718834609570952);
\draw[draw=none,fill=saddlebrown166543] (axis cs:3,0.178209703868043) rectangle (axis cs:3.3,0.745924568341854);
\draw[draw=none,fill=saddlebrown166543] (axis cs:4,0.180024233148396) rectangle (axis cs:4.3,0.729578470760793);
\draw[draw=none,fill=saddlebrown166543] (axis cs:5,0.218388241130199) rectangle (axis cs:5.3,0.624840028846915);
\draw[draw=none,fill=saddlebrown166543] (axis cs:6,0.186790835456915) rectangle (axis cs:6.3,0.585628065101145);
\draw[draw=none,fill=saddlebrown166543] (axis cs:7,0.184949405788757) rectangle (axis cs:7.3,0.626214516843513);
\end{axis}

\end{tikzpicture}

%% file: plots/horeka/gmres_amg_gs2_avg_speedup.tex
\begin{tikzpicture}
    \begin{axis}[
                ymin=1, ymax=1.6,
                xmin=0.8, xmax=3.2,
                xtick={1, 2, 3},
                xticklabels={1, 2, 3},
                width  = 0.44\linewidth,
                height = 0.2\textheight,
                major x tick style = transparent,
                grid = minor,	
                ymajorgrids = true,
                grid style={dashed, gray!40},
                minor y tick num=1,
                ylabel = {\normalsize{Average speedup}},
                xlabel = {outer sweeps},
                tick label style={font={\normalsize}},
                scaled y ticks = false,
                enlarge x limits=0.035,
                legend cell align=left,
                legend style={font=\normalsize},
                legend columns=1,
                legend style={
                    legend pos=north west,
                    draw=none
                },
            ]


\addplot[mark=*, mark size=1.5pt, mark options={red}, draw=red ] plot coordinates{(1,1.09504957) (2,1.16721039) (3,1.16131315)}; 
\addplot[mark=square*, mark size=1.5pt, mark options={black}, draw=black ] plot coordinates{(1,1.15536367) (2,1.18935935) (3,1.17449426)};  

\legend{$\gamma=1$, $\gamma=2$}
    \end{axis}
\end{tikzpicture}

%% file: plots/tg097/gmres_amg_gs2_avg_speedup.tex
\begin{tikzpicture}
    \begin{axis}[
                ymin=1, ymax=1.6,
                xmin=0.8, xmax=3.2,
                xtick={1, 2, 3},
                xticklabels={1, 2, 3},
                width  = 0.44\linewidth,
                height = 0.2\textheight,
                major x tick style = transparent,
                grid = minor,	
                ymajorgrids = true,
                grid style={dashed, gray!40},
                minor y tick num=1,
				ylabel = {\normalsize{Average speedup}},
                xlabel = {outer sweeps},
                tick label style={font={\normalsize}},
                scaled y ticks = false,
                enlarge x limits=0.035,
                legend cell align=left,
                legend style={font=\normalsize},
                legend columns=1,
                legend style={
                    legend pos=north west,
                    draw=none
                },
            ]
            

\addplot[mark=*, mark size=1.5pt, mark options={red}, draw=red ] plot coordinates{(1,1.12414905) (2,1.37147976) (3,1.56555855)}; 
\addplot[mark=square*, mark size=1.5pt, mark options={black}, draw=black ] plot coordinates{(1,1.16749168) (2,1.41731313) (3,1.55868256)};  

    \end{axis}
\end{tikzpicture}

%% file: plots/horeka/gmres_amg_cheb3.tex
\begin{tikzpicture}

\definecolor{burlywood253190133}{RGB}{253,190,133}
\definecolor{coral25314160}{RGB}{253,141,60}
\definecolor{cornflowerblue107174214}{RGB}{107,174,214}
\definecolor{darkgray176}{RGB}{176,176,176}
\definecolor{lightgray204}{RGB}{204,204,204}
\definecolor{orangered2308513}{RGB}{230,85,13}
\definecolor{powderblue189215231}{RGB}{189,215,231}
\definecolor{saddlebrown166543}{RGB}{166,54,3}
\definecolor{steelblue49130189}{RGB}{49,130,189}
\definecolor{teal881156}{RGB}{8,81,156}

\begin{axis}[
	legend cell align={left},
legend style={fill opacity=0.8, draw opacity=1, text opacity=1, draw=lightgray204},
tick align=outside,
tick pos=left,
x grid style={darkgray176},
xtick={0, 1, 2, 3, 4, 5, 6, 7},
xticklabels={1, 2, 3, 4, 5, 6, 7, 8},
width  = 0.49\linewidth,
height = 0.2\textheight,
xmin=-0.68, xmax=7.68,
xtick style={color=black},
y grid style={darkgray176},
ylabel={\small Normalized time},
ymin=0, ymax=1.15,
ytick style={color=black}
]
\draw[draw=none,fill=powderblue189215231] (axis cs:0,0) rectangle (axis cs:-0.3,0.322218128224024);
\draw[draw=none,fill=powderblue189215231] (axis cs:1,0) rectangle (axis cs:0.7,0.24815513986614);
\draw[draw=none,fill=powderblue189215231] (axis cs:2,0) rectangle (axis cs:1.7,0.132462743893071);
\draw[draw=none,fill=powderblue189215231] (axis cs:3,0) rectangle (axis cs:2.7,0.0237602450533985);
\draw[draw=none,fill=powderblue189215231] (axis cs:4,0) rectangle (axis cs:3.7,0.0223133225316277);
\draw[draw=none,fill=powderblue189215231] (axis cs:5,0) rectangle (axis cs:4.7,0.0188760919393411);
\draw[draw=none,fill=powderblue189215231] (axis cs:6,0) rectangle (axis cs:5.7,0.0168734987730373);
\draw[draw=none,fill=powderblue189215231] (axis cs:7,0) rectangle (axis cs:6.7,0.0135855515128303);
\draw[draw=none,fill=cornflowerblue107174214] (axis cs:0,0.322218128224024) rectangle (axis cs:-0.3,0.39867354458364);
\draw[draw=none,fill=cornflowerblue107174214] (axis cs:1,0.24815513986614) rectangle (axis cs:0.7,0.39917624849837);
\draw[draw=none,fill=cornflowerblue107174214] (axis cs:2,0.132462743893071) rectangle (axis cs:1.7,0.26593946842833);
\draw[draw=none,fill=cornflowerblue107174214] (axis cs:3,0.0237602450533985) rectangle (axis cs:2.7,0.0992119341320435);
\draw[draw=none,fill=cornflowerblue107174214] (axis cs:4,0.0223133225316277) rectangle (axis cs:3.7,0.101190458559068);
\draw[draw=none,fill=cornflowerblue107174214] (axis cs:5,0.0188760919393411) rectangle (axis cs:4.7,0.126604645221723);
\draw[draw=none,fill=cornflowerblue107174214] (axis cs:6,0.0168734987730373) rectangle (axis cs:5.7,0.0755494701781142);
\draw[draw=none,fill=cornflowerblue107174214] (axis cs:7,0.0135855515128303) rectangle (axis cs:6.7,0.0700713328226733);
\draw[draw=none,fill=steelblue49130189] (axis cs:0,0.39867354458364) rectangle (axis cs:-0.3,0.428610906411201);
\draw[draw=none,fill=steelblue49130189] (axis cs:1,0.39917624849837) rectangle (axis cs:0.7,0.437188948000686);
\draw[draw=none,fill=steelblue49130189] (axis cs:2,0.26593946842833) rectangle (axis cs:1.7,0.328283914579812);
\draw[draw=none,fill=steelblue49130189] (axis cs:3,0.0992119341320435) rectangle (axis cs:2.7,0.188745520561068);
\draw[draw=none,fill=steelblue49130189] (axis cs:4,0.101190458559068) rectangle (axis cs:3.7,0.200390401687199);
\draw[draw=none,fill=steelblue49130189] (axis cs:5,0.126604645221723) rectangle (axis cs:4.7,0.229997374378189);
\draw[draw=none,fill=steelblue49130189] (axis cs:6,0.0755494701781143) rectangle (axis cs:5.7,0.192907908291514);
\draw[draw=none,fill=steelblue49130189] (axis cs:7,0.0700713328226732) rectangle (axis cs:6.7,0.188817742244351);
\draw[draw=none,fill=teal881156] (axis cs:0,0.428610906411201) rectangle (axis cs:-0.3,1);
\draw[draw=none,fill=teal881156] (axis cs:1,0.437188948000686) rectangle (axis cs:0.7,1);
\draw[draw=none,fill=teal881156] (axis cs:2,0.328283914579812) rectangle (axis cs:1.7,1);
\draw[draw=none,fill=teal881156] (axis cs:3,0.188745520561068) rectangle (axis cs:2.7,1);
\draw[draw=none,fill=teal881156] (axis cs:4,0.200390401687199) rectangle (axis cs:3.7,1);
\draw[draw=none,fill=teal881156] (axis cs:5,0.229997374378189) rectangle (axis cs:4.7,1);
\draw[draw=none,fill=teal881156] (axis cs:6,0.192907908291514) rectangle (axis cs:5.7,1);
\draw[draw=none,fill=teal881156] (axis cs:7,0.188817742244351) rectangle (axis cs:6.7,1);
\draw[draw=none,fill=burlywood253190133] (axis cs:0,0) rectangle (axis cs:0.3,0.452192336035372);
\draw[draw=none,fill=burlywood253190133] (axis cs:1,0) rectangle (axis cs:1.3,0.381757336536811);
\draw[draw=none,fill=burlywood253190133] (axis cs:2,0) rectangle (axis cs:2.3,0.229528345367952);
\draw[draw=none,fill=burlywood253190133] (axis cs:3,0) rectangle (axis cs:3.3,0.0302847523643947);
\draw[draw=none,fill=burlywood253190133] (axis cs:4,0) rectangle (axis cs:4.3,0.0256367724029254);
\draw[draw=none,fill=burlywood253190133] (axis cs:5,0) rectangle (axis cs:5.3,0.0273987184127052);
\draw[draw=none,fill=burlywood253190133] (axis cs:6,0) rectangle (axis cs:6.3,0.0206836436572715);
\draw[draw=none,fill=burlywood253190133] (axis cs:7,0) rectangle (axis cs:7.3,0.0208205668326311);
\draw[draw=none,fill=coral25314160] (axis cs:0,0.452192336035372) rectangle (axis cs:0.3,0.525239498894621);
\draw[draw=none,fill=coral25314160] (axis cs:1,0.381757336536811) rectangle (axis cs:1.3,0.532521022824781);
\draw[draw=none,fill=coral25314160] (axis cs:2,0.229528345367952) rectangle (axis cs:2.3,0.374650483945306);
\draw[draw=none,fill=coral25314160] (axis cs:3,0.0302847523643947) rectangle (axis cs:3.3,0.110648547853614);
\draw[draw=none,fill=coral25314160] (axis cs:4,0.0256367724029254) rectangle (axis cs:4.3,0.10535513816606);
\draw[draw=none,fill=coral25314160] (axis cs:5,0.0273987184127052) rectangle (axis cs:5.3,0.135699704558411);
\draw[draw=none,fill=coral25314160] (axis cs:6,0.0206836436572715) rectangle (axis cs:6.3,0.0788321964387552);
\draw[draw=none,fill=coral25314160] (axis cs:7,0.0208205668326311) rectangle (axis cs:7.3,0.0763548448870165);
\draw[draw=none,fill=orangered2308513] (axis cs:0,0.525239498894621) rectangle (axis cs:0.3,0.550110537951363);
\draw[draw=none,fill=orangered2308513] (axis cs:1,0.532521022824781) rectangle (axis cs:1.3,0.56950403295006);
\draw[draw=none,fill=orangered2308513] (axis cs:2,0.374650483945306) rectangle (axis cs:2.3,0.441911199877093);
\draw[draw=none,fill=orangered2308513] (axis cs:3,0.110648547853614) rectangle (axis cs:3.3,0.208188946577886);
\draw[draw=none,fill=orangered2308513] (axis cs:4,0.10535513816606) rectangle (axis cs:4.3,0.208861940574821);
\draw[draw=none,fill=orangered2308513] (axis cs:5,0.135699704558411) rectangle (axis cs:5.3,0.241159815212995);
\draw[draw=none,fill=orangered2308513] (axis cs:6,0.0788321964387552) rectangle (axis cs:6.3,0.198812138942331);
\draw[draw=none,fill=orangered2308513] (axis cs:7,0.0763548448870164) rectangle (axis cs:7.3,0.194200019149751);
\draw[draw=none,fill=saddlebrown166543] (axis cs:0,0.550110537951363) rectangle (axis cs:0.3,1.04753131908622);
\draw[draw=none,fill=saddlebrown166543] (axis cs:1,0.56950403295006) rectangle (axis cs:1.3,1.07216406384074);
\draw[draw=none,fill=saddlebrown166543] (axis cs:2,0.441911199877093) rectangle (axis cs:2.3,0.964357044092795);
\draw[draw=none,fill=saddlebrown166543] (axis cs:3,0.208188946577886) rectangle (axis cs:3.3,0.816805238528891);
\draw[draw=none,fill=saddlebrown166543] (axis cs:4,0.208861940574821) rectangle (axis cs:4.3,0.786949683205915);
\draw[draw=none,fill=saddlebrown166543] (axis cs:5,0.241159815212995) rectangle (axis cs:5.3,0.675309929817839);
\draw[draw=none,fill=saddlebrown166543] (axis cs:6,0.198812138942331) rectangle (axis cs:6.3,0.647602518618787);
\draw[draw=none,fill=saddlebrown166543] (axis cs:7,0.194200019149751) rectangle (axis cs:7.3,0.693363414400613);
\end{axis}

\end{tikzpicture}

%% file: plots/tg097/gmres_amg_cheb3.tex
\begin{tikzpicture}

\definecolor{burlywood253190133}{RGB}{253,190,133}
\definecolor{coral25314160}{RGB}{253,141,60}
\definecolor{cornflowerblue107174214}{RGB}{107,174,214}
\definecolor{darkgray176}{RGB}{176,176,176}
\definecolor{lightgray204}{RGB}{204,204,204}
\definecolor{orangered2308513}{RGB}{230,85,13}
\definecolor{powderblue189215231}{RGB}{189,215,231}
\definecolor{saddlebrown166543}{RGB}{166,54,3}
\definecolor{steelblue49130189}{RGB}{49,130,189}
\definecolor{teal881156}{RGB}{8,81,156}

\begin{axis}[
	legend cell align={left},
legend style={fill opacity=0.8, draw opacity=1, text opacity=1, draw=lightgray204},
tick align=outside,
tick pos=left,
x grid style={darkgray176},
xtick={0, 1, 2, 3, 4, 5, 6, 7},
xticklabels={1, 2, 3, 4, 5, 6, 7, 8},
width  = 0.49\linewidth,
height = 0.2\textheight,
xmin=-0.68, xmax=7.68,
xtick style={color=black},
y grid style={darkgray176},
ylabel={\small Normalized time},
ymin=0, ymax=1.15,
ytick style={color=black},
]
\draw[draw=none,fill=powderblue189215231] (axis cs:0,0) rectangle (axis cs:-0.3,0.361250754071989);
\draw[draw=none,fill=powderblue189215231] (axis cs:1,0) rectangle (axis cs:0.7,0.270247516416905);
\draw[draw=none,fill=powderblue189215231] (axis cs:2,0) rectangle (axis cs:1.7,0.115607449009755);
\draw[draw=none,fill=powderblue189215231] (axis cs:3,0) rectangle (axis cs:2.7,0.020058952727337);
\draw[draw=none,fill=powderblue189215231] (axis cs:4,0) rectangle (axis cs:3.7,0.0158817221026092);
\draw[draw=none,fill=powderblue189215231] (axis cs:5,0) rectangle (axis cs:4.7,0.0165826673517806);
\draw[draw=none,fill=powderblue189215231] (axis cs:6,0) rectangle (axis cs:5.7,0.0126290488781739);
\draw[draw=none,fill=powderblue189215231] (axis cs:7,0) rectangle (axis cs:6.7,0.0119097760024637);
\draw[draw=none,fill=cornflowerblue107174214] (axis cs:0,0.361250754071989) rectangle (axis cs:-0.3,0.435049266036598);
\draw[draw=none,fill=cornflowerblue107174214] (axis cs:1,0.270247516416905) rectangle (axis cs:0.7,0.411096144132009);
\draw[draw=none,fill=cornflowerblue107174214] (axis cs:2,0.115607449009755) rectangle (axis cs:1.7,0.251108483594443);
\draw[draw=none,fill=cornflowerblue107174214] (axis cs:3,0.020058952727337) rectangle (axis cs:2.7,0.0805760927687896);
\draw[draw=none,fill=cornflowerblue107174214] (axis cs:4,0.0158817221026092) rectangle (axis cs:3.7,0.0758802485407646);
\draw[draw=none,fill=cornflowerblue107174214] (axis cs:5,0.0165826673517806) rectangle (axis cs:4.7,0.108054094829943);
\draw[draw=none,fill=cornflowerblue107174214] (axis cs:6,0.0126290488781739) rectangle (axis cs:5.7,0.0624258758797682);
\draw[draw=none,fill=cornflowerblue107174214] (axis cs:7,0.0119097760024637) rectangle (axis cs:6.7,0.0606256170575075);
\draw[draw=none,fill=steelblue49130189] (axis cs:0,0.435049266036598) rectangle (axis cs:-0.3,0.461492057108385);
\draw[draw=none,fill=steelblue49130189] (axis cs:1,0.411096144132009) rectangle (axis cs:0.7,0.447381714093282);
\draw[draw=none,fill=steelblue49130189] (axis cs:2,0.251108483594443) rectangle (axis cs:1.7,0.307123854566952);
\draw[draw=none,fill=steelblue49130189] (axis cs:3,0.0805760927687896) rectangle (axis cs:2.7,0.16934988487867);
\draw[draw=none,fill=steelblue49130189] (axis cs:4,0.0758802485407647) rectangle (axis cs:3.7,0.171686327802019);
\draw[draw=none,fill=steelblue49130189] (axis cs:5,0.108054094829943) rectangle (axis cs:4.7,0.214750352578469);
\draw[draw=none,fill=steelblue49130189] (axis cs:6,0.0624258758797682) rectangle (axis cs:5.7,0.186498716802115);
\draw[draw=none,fill=steelblue49130189] (axis cs:7,0.0606256170575075) rectangle (axis cs:6.7,0.184462831625952);
\draw[draw=none,fill=teal881156] (axis cs:0,0.461492057108385) rectangle (axis cs:-0.3,1);
\draw[draw=none,fill=teal881156] (axis cs:1,0.447381714093282) rectangle (axis cs:0.7,1);
\draw[draw=none,fill=teal881156] (axis cs:2,0.307123854566952) rectangle (axis cs:1.7,1);
\draw[draw=none,fill=teal881156] (axis cs:3,0.16934988487867) rectangle (axis cs:2.7,1);
\draw[draw=none,fill=teal881156] (axis cs:4,0.171686327802019) rectangle (axis cs:3.7,1);
\draw[draw=none,fill=teal881156] (axis cs:5,0.214750352578469) rectangle (axis cs:4.7,1);
\draw[draw=none,fill=teal881156] (axis cs:6,0.186498716802115) rectangle (axis cs:5.7,1);
\draw[draw=none,fill=teal881156] (axis cs:7,0.184462831625952) rectangle (axis cs:6.7,1);
\draw[draw=none,fill=burlywood253190133] (axis cs:0,0) rectangle (axis cs:0.3,0.455962195857631);
\draw[draw=none,fill=burlywood253190133] (axis cs:1,0) rectangle (axis cs:1.3,0.380703822192288);
\draw[draw=none,fill=burlywood253190133] (axis cs:2,0) rectangle (axis cs:2.3,0.170322199231451);
\draw[draw=none,fill=burlywood253190133] (axis cs:3,0) rectangle (axis cs:3.3,0.0281417517003045);
\draw[draw=none,fill=burlywood253190133] (axis cs:4,0) rectangle (axis cs:4.3,0.0224145129468782);
\draw[draw=none,fill=burlywood253190133] (axis cs:5,0) rectangle (axis cs:5.3,0.0227508773464233);
\draw[draw=none,fill=burlywood253190133] (axis cs:6,0) rectangle (axis cs:6.3,0.0190183050122487);
\draw[draw=none,fill=burlywood253190133] (axis cs:7,0) rectangle (axis cs:7.3,0.01758557285581);
\draw[draw=none,fill=coral25314160] (axis cs:0,0.455962195857631) rectangle (axis cs:0.3,0.52976070782224);
\draw[draw=none,fill=coral25314160] (axis cs:1,0.380703822192288) rectangle (axis cs:1.3,0.521131503620138);
\draw[draw=none,fill=coral25314160] (axis cs:2,0.170322199231451) rectangle (axis cs:2.3,0.302749039314218);
\draw[draw=none,fill=coral25314160] (axis cs:3,0.0281417517003045) rectangle (axis cs:3.3,0.0857112553749067);
\draw[draw=none,fill=coral25314160] (axis cs:4,0.0224145129468782) rectangle (axis cs:4.3,0.0810650031790733);
\draw[draw=none,fill=coral25314160] (axis cs:5,0.0227508773464233) rectangle (axis cs:5.3,0.112728902688342);
\draw[draw=none,fill=coral25314160] (axis cs:6,0.0190183050122487) rectangle (axis cs:6.3,0.0687507290897073);
\draw[draw=none,fill=coral25314160] (axis cs:7,0.01758557285581) rectangle (axis cs:7.3,0.0661417715098322);
\draw[draw=none,fill=orangered2308513] (axis cs:0,0.52976070782224) rectangle (axis cs:0.3,0.572591996782626);
\draw[draw=none,fill=orangered2308513] (axis cs:1,0.521131503620138) rectangle (axis cs:1.3,0.580232362350564);
\draw[draw=none,fill=orangered2308513] (axis cs:2,0.302749039314218) rectangle (axis cs:2.3,0.392314513745197);
\draw[draw=none,fill=orangered2308513] (axis cs:3,0.0857112553749067) rectangle (axis cs:3.3,0.224665573640088);
\draw[draw=none,fill=orangered2308513] (axis cs:4,0.0810650031790733) rectangle (axis cs:4.3,0.194480146482963);
\draw[draw=none,fill=orangered2308513] (axis cs:5,0.112728902688342) rectangle (axis cs:5.3,0.227819261170451);
\draw[draw=none,fill=orangered2308513] (axis cs:6,0.0687507290897072) rectangle (axis cs:6.3,0.192325358712136);
\draw[draw=none,fill=orangered2308513] (axis cs:7,0.0661417715098321) rectangle (axis cs:7.3,0.189623470376711);
\draw[draw=none,fill=saddlebrown166543] (axis cs:0,0.572591996782626) rectangle (axis cs:0.3,1.09843153026342);
\draw[draw=none,fill=saddlebrown166543] (axis cs:1,0.580232362350564) rectangle (axis cs:1.3,1.02037380030308);
\draw[draw=none,fill=saddlebrown166543] (axis cs:2,0.392314513745197) rectangle (axis cs:2.3,0.907981081879988);
\draw[draw=none,fill=saddlebrown166543] (axis cs:3,0.224665573640088) rectangle (axis cs:3.3,0.853325260845932);
\draw[draw=none,fill=saddlebrown166543] (axis cs:4,0.194480146482963) rectangle (axis cs:4.3,0.705470461907062);
\draw[draw=none,fill=saddlebrown166543] (axis cs:5,0.227819261170451) rectangle (axis cs:5.3,0.625584624298943);
\draw[draw=none,fill=saddlebrown166543] (axis cs:6,0.192325358712136) rectangle (axis cs:6.3,0.572512831978847);
\draw[draw=none,fill=saddlebrown166543] (axis cs:7,0.189623470376711) rectangle (axis cs:7.3,0.592702417506771);
\end{axis}

\end{tikzpicture}

%% file: plots/gmres_amg_cheb_sweep3_level0_contrib.tex
\begin{tikzpicture}
    \begin{axis}[
                ymin=0, ymax=1.1,
                xmin=0.8, xmax=8.2,
                width  = 0.44\linewidth,
                height = 0.2\textheight,
                major x tick style = transparent,
                grid = minor,	
                ymajorgrids = true,
                grid style={dashed, gray!40},
                minor y tick num=1,
                ylabel = {\normalsize{Runtime fraction}},
                xlabel = {Matrix},
                xtick = {1,2,3,4,5,6,7,8},
                xticklabels = {1,2,3,4,5,6,7,8},
                tick label style={font={\normalsize}},
                scaled y ticks = false,
                enlarge x limits=0.035,
                legend cell align=left,
                legend style={font=\normalsize},
                legend columns=1,
                legend style={
                    legend pos=south east,
                    draw=none
                },
            ]
            

\addplot[mark=*, mark size=1.5pt, mark options={red}, draw=red ] plot coordinates{(1, .40512655166854747702)
	(2, .56090867510291202927)
	(3, .62728146013448607108)
	(4, .71855653069753428385)
	(5, .78863044500998336722)
	(6, .87450387682628623388)
	(7, .90962498933668550384)
	(8, .91750164881617212017)}; 
\addplot[mark=square*, mark size=1.5pt, mark options={black}, draw=black ] plot coordinates{(1, .33793876026885735623)
	(2, .56870810481413772090)
	(3, .63472696245733788395)
	(4, .76047286153584235111)
	(5, .79321875319147534286)
	(6, .91422604627850639044)
	(7, .95092954683284886760)
	(8, .94497030412243753314)}; 

\legend{ICL, ROME}
    \end{axis}
\end{tikzpicture}

%% file: plots/gmres_amg_cheb_avg_speedup.tex
\begin{tikzpicture}
    \begin{axis}[
                ymin=1, ymax=1.6,
                xmin=0.8, xmax=3.2,
                xtick={1, 2, 3},
                xticklabels={2, 3, 4},
                width  = 0.44\linewidth,
                height = 0.2\textheight,
                major x tick style = transparent,
                grid = minor,	
                ymajorgrids = true,
                grid style={dashed, gray!40},
                minor y tick num=1,
                ylabel = {\normalsize{Average speedup}},
                xlabel = {polynomial degree},
                tick label style={font={\normalsize}},
                scaled y ticks = false,
                enlarge x limits=0.035,
                legend cell align=left,
                legend style={font=\normalsize},
                legend columns=1,
                legend style={
                    legend pos=north west,
                    draw=none
                },
            ]
            

\addplot[mark=*, mark size=1.5pt, mark options={red}, draw=red ] plot coordinates{(1,1.13830537) (2,1.23581173) (3,1.34084653)}; 
\addplot[mark=square*, mark size=1.5pt, mark options={black}, draw=black ] plot coordinates{(1,1.15460943) (2,1.32669053) (3,1.48285605)}; 

    \end{axis}
\end{tikzpicture}

%% file: plots/tg097/NALU_wind/NALU_wind_legend.tex
\begin{tikzpicture}
\definecolor{burlywood253190133}{RGB}{253,190,133}
\definecolor{coral25314160}{RGB}{253,141,60}
\definecolor{cornflowerblue107174214}{RGB}{107,174,214}
\definecolor{darkgray176}{RGB}{176,176,176}
\definecolor{lightgray204}{RGB}{204,204,204}
\definecolor{orangered2308513}{RGB}{230,85,13}
\definecolor{powderblue189215231}{RGB}{189,215,231}
\definecolor{saddlebrown166543}{RGB}{166,54,3}
\definecolor{steelblue49130189}{RGB}{49,130,189}
\definecolor{teal881156}{RGB}{8,81,156}

  \centering
  \begin{customlegend}[
      legend columns=4,
      legend style={
      anchor=north,
      /tikz/every even column/.append style={column sep=0.25cm}},
      legend entries={Total time, Solve time, Total time w RACE, Solve time w RACE},
      width  = \linewidth,
      legend image post style={xscale=0.5},
      legend image code/.code={\draw[#1, draw=none] (0cm,-0.15cm) rectangle (0.26cm,0.15cm);}, 
	  legend style={draw=none} 
      ]
    \addlegendimage{fill=teal881156, color=teal881156} 
    \addlegendimage{fill=cornflowerblue107174214, color=cornflowerblue107174214}
	
	\addlegendimage{fill=orangered2308513, color=orangered2308513}
	\addlegendimage{fill=coral25314160, color=coral25314160}
  \end{customlegend}
\end{tikzpicture}

%% file: plots/tg097/NALU_wind/NALU_wind_gmres_poly.tex
\definecolor{burlywood253190133}{RGB}{253,190,133}
\definecolor{coral25314160}{RGB}{253,141,60}
\definecolor{cornflowerblue107174214}{RGB}{107,174,214}
\definecolor{darkgray176}{RGB}{176,176,176}
\definecolor{lightgray204}{RGB}{204,204,204}
\definecolor{orangered2308513}{RGB}{230,85,13}
\definecolor{powderblue189215231}{RGB}{189,215,231}
\definecolor{saddlebrown166543}{RGB}{166,54,3}
\definecolor{steelblue49130189}{RGB}{49,130,189}
\definecolor{teal881156}{RGB}{8,81,156}

\begin{tikzpicture}
    \begin{axis}[
                ymin=0, ymax=3,
                xmin=2.8, xmax=8.2,
                width  = 0.42\linewidth,
                height = 0.2\textheight,
                major x tick style = transparent,
                grid = minor,	
                ymajorgrids = true,
                grid style={dashed, gray!40},
                minor y tick num=1,
                ylabel = {\normalsize{Time (s)}},
                xlabel = {Polynomial degree},
                xtick = {1,2,3,4,5,6,7,8},
                xticklabels = {1,2,3,4,5,6,7,8},
                tick label style={font={\normalsize}},
                scaled y ticks = false,
                enlarge x limits=0.035,
                legend cell align=left,
                legend style={font=\footnotesize, fill=none},
                legend columns=2,
                legend style={
                    legend pos=south west,
                    draw=none
                },
            	legend image post style={xscale=0.6, yscale=0.6},
            ]

\addlegendimage{no markers, blue, thick}
\addlegendimage{no markers, orange, thick}
\addlegendimage{only marks, mark=*}
\addlegendimage{only marks, mark=square*}

\addplot[mark=square*, mark size=1.5pt, mark options={teal881156}, draw=teal881156 ] plot coordinates{
	(3, 2.838)
	(4, 2.505)
	(5, 2.24)
	(6, 2.302)
	(7, 2.302)
	(8, 2.564)};

\addplot[mark=square*, mark size=1.5pt, mark options={orangered2308513}, draw=orangered2308513 ] plot coordinates{
	(3, 2.356)
	(4, 1.916)
	(5, 1.68)
	(6, 1.69)
	(7, 1.67)
	(8, 1.82)
};

\addplot[mark=*, mark size=1.5pt, mark options={cornflowerblue107174214}, draw=cornflowerblue107174214 ] plot coordinates{
	(3, 2.473)
	(4, 2.084)
	(5, 1.758)
	(6, 1.769)
	(7, 1.698)
	(8, 1.908)};

\addplot[mark=*, mark size=1.5pt, mark options={coral25314160}, draw=coral25314160 ] plot coordinates{
	(3, 2.002)
	(4, 1.51)
	(5, 1.211)
	(6, 1.174)
	(7, 1.086)
	(8, 1.185)
};




    \end{axis}
\end{tikzpicture}

%% file: plots/tg097/NALU_wind/NALU_wind_gmres_poly_iter.tex
\definecolor{burlywood253190133}{RGB}{253,190,133}
\definecolor{coral25314160}{RGB}{253,141,60}
\definecolor{cornflowerblue107174214}{RGB}{107,174,214}
\definecolor{darkgray176}{RGB}{176,176,176}
\definecolor{lightgray204}{RGB}{204,204,204}
\definecolor{orangered2308513}{RGB}{230,85,13}
\definecolor{powderblue189215231}{RGB}{189,215,231}
\definecolor{saddlebrown166543}{RGB}{166,54,3}
\definecolor{steelblue49130189}{RGB}{49,130,189}
\definecolor{teal881156}{RGB}{8,81,156}

\begin{tikzpicture}
    \begin{axis}[
	ymin=0, ymax=15,
	xmin=2.8, xmax=8.2,
	width  = 0.42\linewidth,
	height = 0.2\textheight,
	major x tick style = transparent,
	grid = minor,	
	ymajorgrids = true,
	grid style={dashed, gray!40},
	minor y tick num=1,
	ylabel = {\normalsize{Iter}},
	xlabel = {Polynomial degree},
	xtick = {1,2,3,4,5,6,7,8},
	xticklabels = {1,2,3,4,5,6,7,8},
	tick label style={font={\normalsize}},
	scaled y ticks = false,
	enlarge x limits=0.035,
	legend cell align=left,
	legend style={font=\normalsize},
	legend columns=1,
	legend style={
		legend pos=south east,
		draw=none
	},
	]
	

	\addplot[mark=*, mark size=1pt, mark options={black}, draw=black ] plot coordinates{
		(3, 14)
		(4, 9)
		(5, 6)
		(6, 5)
		(7, 4)
		(8, 4)};\label{plot_one}
	
\end{axis}
    \begin{axis}[
	axis y line*=right,
	axis x line=none,
	ymin=0, ymax=45,
	xmin=2.8, xmax=8.2,
	width  = 0.42\linewidth,
	height = 0.2\textheight,
	major x tick style = transparent,
	grid = minor,	
	ymajorgrids = true,
	grid style={dashed, gray!40},
	minor y tick num=1,
	ylabel = {\normalsize{\#effective SpMVs}},
	xlabel = {Polynomial degree},
	xtick = {1,2,3,4,5,6,7,8},
	xticklabels = {1,2,3,4,5,6,7,8},
	ytick = {0,15,30,45},
	yticklabels = {0,15,30,45},
	tick label style={font={\normalsize}},
	scaled y ticks = false,
	enlarge x limits=0.035,
	legend cell align=left,
	legend style={font=\footnotesize, fill=none},
	legend columns=1,
	legend style={
		legend pos=south west,
		draw=none
	},
	]
	

	\addplot[mark=square*, mark size=1pt, mark options={red}, draw=red ] plot coordinates{
		(3, 42)
		(4, 36)
		(5, 30)
		(6, 30)
		(7, 28)
		(8, 32)};\label{plot_two}
	
\addlegendimage{/pgfplots/refstyle=plot_one}\addlegendentry{\#SpMVs}
\addlegendimage{/pgfplots/refstyle=plot_two}\addlegendentry{Iter}
\end{axis}
\end{tikzpicture}

%% file: plots/tg097/NALU_wind/NALU_wind_gmres_Jacobi.tex
\definecolor{burlywood253190133}{RGB}{253,190,133}
\definecolor{coral25314160}{RGB}{253,141,60}
\definecolor{cornflowerblue107174214}{RGB}{107,174,214}
\definecolor{darkgray176}{RGB}{176,176,176}
\definecolor{lightgray204}{RGB}{204,204,204}
\definecolor{orangered2308513}{RGB}{230,85,13}
\definecolor{powderblue189215231}{RGB}{189,215,231}
\definecolor{saddlebrown166543}{RGB}{166,54,3}
\definecolor{steelblue49130189}{RGB}{49,130,189}
\definecolor{teal881156}{RGB}{8,81,156}

\begin{tikzpicture}
    \begin{axis}[
                ymin=0, ymax=3,
                xmin=0.8, xmax=8.2,
                width  = 0.44\linewidth,
                height = 0.2\textheight,
                major x tick style = transparent,
                grid = minor,	
                ymajorgrids = true,
                grid style={dashed, gray!40},
                minor y tick num=1,
                ylabel = {\normalsize{Time (s)}},
                xlabel = {Sweeps},
                xtick = {1,2,3,4,5,6,7,8},
                xticklabels = {1,2,3,4,5,6,7,8},
                tick label style={font={\normalsize}},
                scaled y ticks = false,
                enlarge x limits=0.035,
                legend cell align=left,
                legend style={font=\normalsize},
                legend columns=1,
                legend style={
                    legend pos=south east,
                    draw=none
                },
            ]

\addplot[mark=square*, mark size=1.5pt, mark options={teal881156}, draw=teal881156 ] plot coordinates{
	(1, 1.869)
	(2, 1.39)
	(3, 1.47)
	(4, 1.48)
	(5, 1.40)
	(6, 1.63)
	(7, 1.57)
	(8, 1.76)};

\addplot[mark=*, mark size=1.5pt, mark options={orangered2308513}, draw=orangered2308513 ] plot coordinates{
	(1, 1.868)
	(2, 1.087)
	(3, 0.8991)
	(4, 0.8735)
	(5, 0.7002)
	(6, 0.7977)
	(7, 0.8161)
	(8, 0.8927)
};

\addplot[mark=square*, mark size=1.5pt, mark options={cornflowerblue107174214}, draw=cornflowerblue107174214 ] plot coordinates{
	(1, 1.756)
	(2, 1.284)
	(3, 1.36)
	(4, 1.375)
	(5, 1.292)
	(6, 1.522)
	(7, 1.467)
	(8, 1.656)};

\addplot[mark=*, mark size=1.5pt, mark options={coral25314160}, draw=coral25314160 ] plot coordinates{
	(1, 1.769)
	(2, 0.9882)
	(3, 0.7999)
	(4, 0.7734)
	(5, 0.6011)
	(6, 0.6989)
	(7, 0.7173)
	(8, 0.7933)
}; 

    \end{axis}
\end{tikzpicture}

%% file: plots/tg097/NALU_wind/NALU_wind_gmres_Jacobi_iter.tex
\definecolor{burlywood253190133}{RGB}{253,190,133}
\definecolor{coral25314160}{RGB}{253,141,60}
\definecolor{cornflowerblue107174214}{RGB}{107,174,214}
\definecolor{darkgray176}{RGB}{176,176,176}
\definecolor{lightgray204}{RGB}{204,204,204}
\definecolor{orangered2308513}{RGB}{230,85,13}
\definecolor{powderblue189215231}{RGB}{189,215,231}
\definecolor{saddlebrown166543}{RGB}{166,54,3}
\definecolor{steelblue49130189}{RGB}{49,130,189}
\definecolor{teal881156}{RGB}{8,81,156}

\begin{tikzpicture}
    \begin{axis}[
	ymin=0, ymax=35,
	xmin=0.8, xmax=8.2,
	width  = 0.42\linewidth,
	height = 0.2\textheight,
	major x tick style = transparent,
	grid = minor,	
	ymajorgrids = true,
	grid style={dashed, gray!40},
	minor y tick num=1,
	ylabel = {\normalsize{Iter}},
	xlabel = {Polynomial degree},
	xtick = {1,2,3,4,5,6,7,8},
	xticklabels = {1,2,3,4,5,6,7,8},
	tick label style={font={\normalsize}},
	scaled y ticks = false,
	enlarge x limits=0.035,
	legend cell align=left,
	legend style={font=\normalsize},
	legend columns=1,
	legend style={
		legend pos=south east,
		draw=none
	},
	]
	

	\addplot[mark=*, mark size=1pt, mark options={black}, draw=black ] plot coordinates{
		(1, 24)
		(2, 12)
		(3, 9)
		(4, 7)
		(5, 5)
		(6, 5)
		(7, 4)
		(8, 4)
};\label{NALU_Jacobi_plot_one}
	
\end{axis}
    \begin{axis}[
	axis y line*=right,
	axis x line=none,
	ymin=0, ymax=35,
	xmin=0.8, xmax=8.2,
	width  = 0.42\linewidth,
	height = 0.2\textheight,
	major x tick style = transparent,
	grid = minor,	
	ymajorgrids = true,
	grid style={dashed, gray!40},
	minor y tick num=1,
	ylabel = {\normalsize{\#effective SpMVs}},
	xlabel = {Sweeps},
	xtick = {1,2,3,4,5,6,7,8},
	xticklabels = {1,2,3,4,5,6,7,8},,
	tick label style={font={\normalsize}},
	scaled y ticks = false,
	enlarge x limits=0.035,
	legend cell align=left,
	legend style={font=\footnotesize, fill=none},
	legend columns=1,
	legend style={
		legend pos=south west,
		draw=none
	},
	]
	

	\addplot[mark=square*, mark size=1pt, mark options={red}, draw=red ] plot coordinates{
		(1, 24)
		(2, 24)
		(3, 27)
		(4, 28)
		(5, 25)
		(6, 30)
		(7, 28)
		(8, 32)};\label{NALU_Jacobi_plot_two}
	
\end{axis}
\end{tikzpicture}

%% file: plots/GMRES_Transport_spmv.tex
\pgfplotstableread[col sep=comma] {plots/gmres_s_step_kernels_perf_data/tg097/dynamic_and_blis_optimal_config/table.txt} \fileone
\pgfplotstableread[col sep=comma] {plots/gmres_s_step_kernels_perf_data/tg097/static_and_blis_optimal_config/table.txt} \filetwo
\pgfplotstableread[col sep=comma] {plots/gmres_s_step_kernels_perf_data/horeka/dynamic_and_mkl/table.txt} \filethree
\pgfplotstableread[col sep=comma] {plots/gmres_s_step_kernels_perf_data/horeka/static_and_mkl/table.txt} \filefour

\def\barOne{default}
\def\barTwo{static}
\def\barThree{default}
\def\barFour{static}
\def\nnz{23500731}
\def\iter{970}

\begin{tikzpicture}-
	\begin{axis}[
		tick align=outside,
		tick pos=left,
		width  = 0.48\linewidth,
        height = 0.22\textheight,
		ymin=0,
		ymax=28,
		ybar,
		bar width=7pt,
		enlarge x limits=0.6,
		minor y tick num=1,
		ylabel={Perf (\GFS)},
		xtick={0,1,2},
		xticklabels={\barOne, \barTwo},
		xticklabel style={rotate=0},
		legend image post style={xscale=0.5},
		legend image code/.code={\draw[#1, draw=none] (0cm,-0.15cm) rectangle (0.26cm,0.15cm);},
		legend style={
			legend pos=north west,
			draw=none },
		ylabel shift = -3 pt
		]
		
		\pgfkeys{/pgf/fpu}

		\pgfplotstablegetelem{0}{SpMV}\of\fileone
		\pgfmathsetmacro{\spmv}{\pgfplotsretval}
		
		\pgfplotstablegetelem{0}{MPK}\of\fileone
		\pgfmathsetmacro{\mpk}{\pgfplotsretval}
		
		\pgfmathsetmacro{\matrixApply}{2*\iter*\nnz*1e-9/(\mpk+\spmv)}

		\pgfplotstablegetelem{0}{SpMV}\of\filetwo
		\pgfmathsetmacro{\spmvTwo}{\pgfplotsretval}
		
		\pgfplotstablegetelem{0}{MPK}\of\filetwo
		\pgfmathsetmacro{\mpkTwo}{\pgfplotsretval}
		
		\pgfmathsetmacro{\matrixApplyTwo}{2*\iter*\nnz*1e-9/(\mpkTwo+\spmvTwo)}

		\pgfplotstablegetelem{0}{SpMV}\of\filethree
		\pgfmathsetmacro{\spmvThree}{\pgfplotsretval}

		\pgfplotstablegetelem{0}{MPK}\of\filethree
		\pgfmathsetmacro{\mpkThree}{\pgfplotsretval}

		\pgfmathsetmacro{\matrixApplyThree}{2*\iter*\nnz*1e-9/(\mpkThree+\spmvThree)}

		\pgfplotstablegetelem{0}{SpMV}\of\filefour
		\pgfmathsetmacro{\spmvFour}{\pgfplotsretval}

		\pgfplotstablegetelem{0}{MPK}\of\filefour
		\pgfmathsetmacro{\mpkFour}{\pgfplotsretval}

		\pgfmathsetmacro{\matrixApplyFour}{2*\iter*\nnz*1e-9/(\mpkFour+\spmvFour)}
		
		\pgfkeys{/pgf/fpu=false}
		
		\addplot[ybar, fill=cyan] plot coordinates {(0,\matrixApplyThree) (1,\matrixApplyFour)};
		\addplot[ybar, fill=ao(english)] plot coordinates {(0,\matrixApply) (1,\matrixApplyTwo)};

		\legend{ICL, ROME}

	\end{axis}
\end{tikzpicture}

%% file: plots/GMRES_Transport_ortho.tex
\pgfplotstableread[col sep=comma] {plots/gmres_s_step_kernels_perf_data/tg097/static_and_mkl_w_fakeintel/table.txt} \fileone
\pgfplotstableread[col sep=comma] {plots/gmres_s_step_kernels_perf_data/tg097/static_and_blis_default_config/table.txt} \filetwo
\pgfplotstableread[col sep=comma] {plots/gmres_s_step_kernels_perf_data/tg097/static_and_blis_optimal_config/table.txt} \filethree
\pgfplotstableread[col sep=comma] {plots/gmres_s_step_kernels_perf_data/horeka/static_and_mkl/table.txt} \filefour

\def\barOne{MKL}
\def\barTwo{BLIS*}
\def\barThree{BLIS}
\def\iter{247}



\begin{tikzpicture}-
	\begin{axis}[
		tick align=outside,
		tick pos=left,
		width  = 0.48\linewidth,
		height = 0.22\textheight,
		ymin=0,
		ybar,
		bar width=7pt,
		enlarge x limits=0.3,
		legend style={at={(0.5,-0.20)},
			anchor=north,legend columns=2},
		minor y tick num=1,
		ylabel={Perf (\#Ortho/s)},
		xtick={0,1,2},
		xticklabels={\footnotesize{\barOne}, \footnotesize{\barTwo}, \footnotesize{\barThree}},
		xticklabel style={rotate=0},
		ylabel shift = -3 pt
		]
		
		\pgfkeys{/pgf/fpu}

		\pgfplotstablegetelem{0}{Bortho}\of\fileone
		\pgfmathsetmacro{\bortho}{\pgfplotsretval}
		\pgfplotstablegetelem{0}{TSQR}\of\fileone
		\pgfmathsetmacro{\tsqr}{\pgfplotsretval}
		\pgfmathsetmacro{\orthoOne}{\iter/(\bortho+\tsqr)} 

		\pgfplotstablegetelem{0}{Bortho}\of\filetwo
		\pgfmathsetmacro{\bortho}{\pgfplotsretval}
		\pgfplotstablegetelem{0}{TSQR}\of\filetwo
		\pgfmathsetmacro{\tsqr}{\pgfplotsretval}
		\pgfmathsetmacro{\orthoTwo}{\iter/(\bortho+\tsqr)}
	
		\pgfplotstablegetelem{0}{Bortho}\of\filethree
		\pgfmathsetmacro{\bortho}{\pgfplotsretval}
		\pgfplotstablegetelem{0}{TSQR}\of\filethree
		\pgfmathsetmacro{\tsqr}{\pgfplotsretval}	
		\pgfmathsetmacro{\orthoThree}{\iter/(\bortho+\tsqr)} 

		\pgfplotstablegetelem{0}{Bortho}\of\filefour
		\pgfmathsetmacro{\bortho}{\pgfplotsretval}
		\pgfplotstablegetelem{0}{TSQR}\of\filefour
		\pgfmathsetmacro{\tsqr}{\pgfplotsretval}	
		\pgfmathsetmacro{\orthoFour}{\iter/(\bortho+\tsqr)} 
		
		\pgfkeys{/pgf/fpu=false}
		
		\addplot[ybar, fill=cyan] plot coordinates {(0,\orthoFour)};
		\addplot[ybar, fill=ao(english)] plot coordinates
		{(0,\orthoOne) (1-0.15,\orthoTwo) (2-0.15,\orthoThree)}; 

		
	\end{axis}
\end{tikzpicture}